\input amstex
\documentstyle{amsppt}
\magnification =\magstep1
\pageheight{9truein}
\baselineskip=18pt
\input EPSfig.macros
\nologo
\overfullrule=0pt

\newif\ifproofmode                      
\proofmodefalse				

\newif\ifforwardreference		
\forwardreferencefalse			

\newif\ifchapternumbers			
\chapternumbersfalse			

\newif\ifcontinuousnumbering		
\continuousnumberingfalse		

\newif\iffigurechapternumbers		
\figurechapternumbersfalse		

\newif\ifcontinuousfigurenumbering	
\continuousfigurenumberingfalse		

\font\eqsixrm=cmr6			
\def\marginstyle{\eqsixrm}		

\newtoks\chapletter			
\newcount\chapno			
\newcount\eqlabelno			
\newcount\figureno			

\chapno=0
\eqlabelno=0
\figureno=0


\def\chapfolio{\ifnum \chapno>0 \the\chapno \else \the\chapletter \fi}


\def\bumpchapno{\ifnum \chapno>-1 \global \advance \chapno by 1
	\else \global \advance \chapno by -1 \setletter\chapno \fi
	\ifcontinuousnumbering \else \global\eqlabelno=0 \fi
	\ifcontinuousfigurenumbering \else \global\figureno=0 \fi}

%


%

\def\tempsetletter#1{\ifcase-#1 {}\or{} \or\chapletter={A}\or\chapletter={B}
  \or\chapletter={C} \or\chapletter={D} \or\chapletter={E}
  \or\chapletter={F} \or\chapletter={G} \or\chapletter={H}
  \or\chapletter={I} \or\chapletter={J} \or\chapletter={K}
  \or\chapletter={L} \or\chapletter={M} \or\chapletter={N}
  \or\chapletter={O} \or\chapletter={P} \or\chapletter={Q}
  \or\chapletter={R} \or\chapletter={S} \or\chapletter={T}
  \or\chapletter={U} \or\chapletter={V} \or\chapletter={W}
  \or\chapletter={X} \or\chapletter={Y} \or\chapletter={Z}\fi}

%

\def\chapshow#1{\ifnum #1>0 \relax #1%
   \else {\tempsetletter{\number#1}\chapno=#1 \chapfolio} \fi}

%
\def\today{\ifcase\month\or
January\or February\or March\or April\or May\or June\or
July\or August\or September\or October\or November\or December\fi
\space\number\day, \number\year}

%

\def\initialeqmacro{\ifproofmode
 \headline{\tenrm \today\hfill \jobname\ --- draft\hfill\folio}
     \hoffset=-1cm \immediate\openout2=allcrossreferfile \fi
 \ifforwardreference \input labelfile
     \ifproofmode \immediate\openout1=labelfile \fi \fi}


%

\def\chaplabel#1{\bumpchapno \ifproofmode \ifforwardreference
   \immediate\write1{\noexpand\expandafter\noexpand\def
   \noexpand\csname CHAPLABEL#1\endcsname{\the\chapno}}\fi\fi
   \global\expandafter\edef\csname CHAPLABEL#1\endcsname
   {\the\chapno}\ifproofmode\llap{\hbox{\marginstyle #1\ }}\fi\chapfolio}

%
\def\eqnum{\global\advance\eqlabelno by 1
   \eqno(\ifchapternumbers\chapfolio.\fi\the\eqlabelno)}

\def\eqlabel#1{\global\advance\eqlabelno by 1 \ifproofmode\ifforwardreference
 \immediate\write1{\noexpand\expandafter\noexpand\def
 \noexpand\csname EQLABEL#1\endcsname{\the\chapno.\the\eqlabelno?}}\fi\fi
 \global\expandafter\edef\csname EQLABEL#1\endcsname
 {\the\chapno.\the\eqlabelno?} \eqno(\ifchapternumbers\chapfolio.\fi
 \the\eqlabelno)\ifproofmode\rlap{\hbox{\marginstyle #1}}\fi}

\def\leqlabel#1{\global\advance\eqlabelno by 1 \ifproofmode\ifforwardreference
 \immediate\write1{\noexpand\expandafter\noexpand\def
 \noexpand\csname EQLABEL#1\endcsname{\the\chapno.\the\eqlabelno?}}\fi\fi
 \global\expandafter\edef\csname EQLABEL#1\endcsname
 {\the\chapno.\the\eqlabelno?} \leqno(\ifchapternumbers\chapfolio.\fi
 \the\eqlabelno)\ifproofmode\rlap{\hbox{\marginstyle #1}}\fi}

\def\eqalignnum{\global\advance\eqlabelno by 1
   &(\ifchapternumbers\chapfolio.\fi\the\eqlabelno)}

\def\eqalignlabel#1{\global\advance\eqlabelno by1 \ifproofmode
 \ifforwardreference\immediate\write1{\noexpand\expandafter\noexpand\def
 \noexpand\csname EQLABEL#1\endcsname
     {\the\chapno.\the\eqlabelno?}}\fi\fi
 \global\expandafter\edef\csname EQLABEL#1\endcsname
 {\the\chapno.\the\eqlabelno?}\ifchapternumbers\chapfolio.\fi
 \the\eqlabelno\ifproofmode\rlap{\hbox{\marginstyle
 #1}}\fi}

\def\eqref#1{(\ifundefined{EQLABEL#1}***\ifproofmode\ifforwardreference)%
   \else \write16{ ***Undefined Equation Reference #1*** }\fi
   \else \write16{ ***Undefined Equation Reference #1*** }\fi
   \else \edef\LABxx{\getlabel{EQLABEL#1}}%
   \def\LAByy{\expandafter\stripchap\LABxx}\ifchapternumbers
   \chapshow{\LAByy}.\expandafter\stripeq\LABxx
   \else\ifnum \number\LAByy=\chapno \relax\expandafter\stripeq\LABxx
   \else\chapshow{\LAByy}.\expandafter\stripeq\LABxx\fi\fi)\fi
   \ifproofmode\write2{Equation #1}\fi}

%

\def\fignum{\global\advance\figureno by 1 \relax
   \iffigurechapternumbers\chapfolio.\fi\the\figureno}\

\def\figlabel#1{\global\advance\figureno by 1\relax
 \ifproofmode\ifforwardreference
 \immediate\write1{\noexpand\expandafter\noexpand\def
 \noexpand\csname FIGLABEL#1\endcsname{\the\chapno.\the\figureno?}}\fi\fi
 \global\expandafter\edef\csname FIGLABEL#1\endcsname
 {\the\chapno.\the\figureno?}\iffigurechapternumbers\chapfolio.\fi
 \ifproofmode$^{\hbox{\marginstyle #1}}$\relax\fi\the\figureno}

\def\figref#1{\ifundefined{FIGLABEL#1}!!!!\ifproofmode\ifforwardreference)%
   \else \write16{ ***Undefined Equation Reference #1*** }\fi
   \else \write16{ ***Undefined Equation Reference #1*** }\fi
   \else \edef\LABxx{\getlabel{FIGLABEL#1}}%
   \def\LAByy{\expandafter\stripchap\LABxx}%
   \iffigurechapternumbers\chapshow{\LAByy}.\expandafter\stripeq\LABxx
   \else\ifnum\number\LAByy=\chapno \relax\expandafter\stripeq\LABxx
   \else\chapshow{\LAByy}.\expandafter\stripeq\LABxx\fi\fi
   \ifproofmode\write2{Figure #1}\fi\fi}

%

%

\def\getlabel#1{\csname#1\endcsname}
\def\ifundefined#1{\expandafter\ifx\csname#1\endcsname\relax}
\def\stripchap#1.#2?{#1}
\def\stripeq#1.#2?{#2}

\figurechapternumberstrue  

\chapternumberstrue        

\def\thmlbl#1{\figlabel{#1}}
\def\thmref{\figref}
\def\eqnlbl#1{\leqlabel{#1}}

\def\eqnref#1{\eqref{#1}}
\def\sectionnumber{\chapno}
\def\theoremnumber{\figureno}
\def\equationnumber{\eqlabelno}

\input amssym.def
\def\notexist{\raise1pt\hbox{$\not$}\exists}

\def\BbbR{\Bbb{R}}
\def\BbbC{\Bbb{C}}

\def\BbbA{\Bbb{A}}


\def\CalO{\Cal{O}}

\def\sgn{\text{\rm \ sgn \ }}
\def\Tr{\text{\rm \ Tr \ }}


\def\div{\roman{div}}

\def\R{\roman{Re\ }}
\def\Im{\roman{Im\ }}
\def\Span{\hbox{Span}}


\def\barU{{\bar{u}}}

\def\T3{{\widetilde{\xi}}}
\def\TPsi{{\widetilde{\Psi}}}

\def\tildeZ{{\widetilde{z}}}

\def\tildeW{{\widetilde{w}}}

\topmatter
\title
Stability of periodic solutions of conservation laws with viscosity:
Analysis of the Evans function
\endtitle
\leftheadtext{Stability of periodic solutions}
\rightheadtext{M. Oh and K. Zumbrun}
\abstract
Nonclassical conservation laws with viscosity
arising in multiphase fluid and solid mechanics
exhibit a rich variety of traveling-wave phenomena,
including homoclinic (pulse-type) and periodic solutions 
along with the standard heteroclinic (shock, or front-type) solutions.
Here, we investigate 
stability of periodic traveling
waves within the abstract Evans function framework
established by R.A. Gardner.
Our main result is to derive a useful {\it stability index}
analogous to that developed by Gardner and Zumbrun in the
traveling-front or -pulse context, giving necessary conditions for
stability with respect to initial perturbations that are periodic
on the same period $T$ as the traveling wave;
moreover, we show that the periodic stability
index has an interpretation analogous to that of the traveling-front or -pulse
index in terms of well-posedness
of an associated Riemann problem for an inviscid medium, 
now to be interpreted
as allowing a wider class of measure-valued solutions, 
or, alternatively, in terms of existence and nonsingularity of a 
local ``mass map'' from perturbation mass to 
potential time-asymptotic $T$-periodic states. 
A closely related calculation yields also a complementary 
{\it long-wave stability criterion} necessary for stability 
with respect to periodic
perturbations of arbitrarily large period $NT$, $N\to \infty$.
We augment these analytical results with numerical investigations 
analogous to those carried out by Brin in the traveling-front or -pulse case,
approximating the spectrum of the linearized operator about the wave.

The stability index and long-wave stability criterion
are explicitly evaluable in the same
planar, Hamiltonian cases as is the index of Gardner\&Zumbrun,
and together yield rigorous results of instability similar to those
obtained previously for pulse-type solutions;
this is established through a novel
{\it dichotomy} asserting that the two criteria are in certain cases 
logically exclusive.  
%
In particular, we obtain results bearing on the
nature and mechanism for formation of 
highly oscillatory  Turing-like patterns 
observed numerically by Frid\&Liu and  \v{C}ani\'c\&Peters in models
of multiphase flow.
Specifically, for the van der Waals model considered by Frid\&Liu, we show
instability of all periodic waves such that the period
increases with amplitude in the one-parameter family 
of nearby periodic orbits, and in particular of large- and 
small-amplitude waves; for the standard, double-well potential,
this yields instability of all periodic waves.
Likewise, for a quadratic-flux model like that considered by 
\v{C}ani\'c\&Peters, we show instability of large-amplitude
waves of the type lying near observed patterns, and of all
small-amplitude waves;
our numerical results give evidence that 
intermediate-amplitude waves are unstable as well.
These results give support for an alternative mechanism for pattern formation
conjectured by Azevedo, Marchesin, Plohr, and Zumbrun, 
not involving periodic waves.
\endabstract
\author
M. Oh and
K. Zumbrun 
\endauthor
\date
November 27, 2000; Revised: February 18, 2002.
\enddate
\address
Department of Mathematics,
Indiana University,
Bloomington, IN  47405-4301, 
USA
\endaddress
\email
myoh\@indiana.edu
\endemail 
\address
Department of Mathematics,
Indiana University,
Bloomington, IN  47405-4301, 
USA
\endaddress
\email
kzumbrun\@indiana.edu
\endemail 
\endtopmatter
\vfil\eject
\document
\sectionnumber=1 
\theoremnumber=0
\equationnumber=0
\smallskip
\TagsOnLeft
\specialhead
\S1.\quad Introduction
\endspecialhead
In this paper, we study stability of periodic traveling-wave
solutions of {\it conservation laws with viscosity}
$$
u_t+f(u)_x=(B(u)u_x)_x,
\eqnlbl{1}
$$
$u$, $f \in \BbbR^n$, $B\in \BbbR^{n\times n}$,
modeling flow in compressible media.
Periodic solutions do not occur
for classical systems \eqnref{1} possessing a convex entropy 
in the sense of [Sm,Ka], for example, ideal gas dynamics
or magnetohydrodynamics, being forbidden by energy considerations.
However, they arise in a natural way in certain
nonclassical systems modeling media with
multiple phases, for example van der Waals gas dynamics
and elasticity, or three-phase flow in porous media,
in conjunction with a host of other complex phenomena
not seen in classical conservation laws:
in particular, hetero- and homoclinic cycles of
traveling-front and -pulse type solutions, and associated
nonuniqueness of Riemann solutions in the corresponding
first-order system
$$
u_t+f(u)_x=(B(u)u_x)_x;
\eqnlbl{hyp}
$$
see, e.g., [Sl.1--5,Sh.1--3,AMPZ.1--4,Z].
As discussed in [AMPZ.2], these features have
the common seed of instability of constant solutions
in certain regions of state space: in the case of
phase-transitional models, regions that are ``between
phases.''
Such instability is frequently (but not always)
associated with change in type from
hyperbolic to elliptic in the first order system \eqnref{hyp}.

In the case of van der Waals gas and solid dynamics, there is
a close relation between the mechanical model \eqnref{1}
and the variational Cahn Hilliard model for phase transition;
see Section 7 for further details. 
In particular, their stationary theories agree, reducing to 
the Euler-Lagrange equations for the associated van der Waals/Cahn Hilliard
energy: a planar Hamiltonian ordinary differential equation.
Thus, it is no surprise that the mechanical model features the
same rich solution structure found in the study of Cahn Hilliard
equations $u_{tt}+ f(u)=u_{xx}$, or Allen Cahn (reaction diffusion) 
equations $u_t+ f(u)=u_{xx}$ based on the same energy functional.
Indeed, the existence of hetero-, homoclinic, and periodic cycles
is already suggested
by the planar Hamiltonian form of the stationary-wave ordinary
differential equation.

On the other hand, the dynamics of these models are quite different, and
so a priori it is not clear to what extent, or in what way, this stationary
structure will be manifested in the asymptotic behavior of solutions.
In particular, one would like to know whether the mechanical model
\eqnref{1} indeed supports persistent transition layers (fronts)
between phases, 
and, more, whether it can successfully
predict, at least at a qualitative level, the experimentally observed
phenomena of nucleation (spontaneous formation of 
transition layers in previously smooth flow)
and pattern formation (e.g., Martensitic
crystal structure in stressed elastic solids).
%
In the case of three-phase flow models commonly 
used in oil recovery simulation, with gas, water and oil
treated as separate phases of a single fluid,
both modeling and experimental data are less certain
than in the case of the van der Waals model.
Nonetheless, the same basic questions are of interest,
now for purposes of experimental prediction and validation
or invalidation of the mathematical model.
These questions must be addressed in the context of the full, {\it dynamical}
behavior of model \eqnref{1},
and ultimately involve the careful study of stability.

Stability and behavior of hetero- and homoclinic, or front- and pulse-type
solutions of phase-transitional models
has been investigated numerically and analytically in, e.g., 
[ZPM,AMPZ.1--3,LZ.1--2,GZ,Z.1--2].
For both van der Waals and three-phase flow models, the picture that
has emerged is that traveling-front solutions connecting one pure
phase to another are stable, as are simple Riemann patterns involving
two or more traveling fronts moving away from each other with
nonzero speed, whereas
traveling-pulse solutions connecting a single phase to itself,
and therefore Riemann patterns in which they appear, are unstable.
In terms of time-asymptotic dynamics, 
the latter seem to play the role of saddle points separating the
basins of attraction of the former, attracting asymptotic states [AMPZ.1].
These studies essentially answer the first two questions posed above,
showing that phase-transitional layers are indeed supported by the models,
and can form spontaneously from smooth initial data.
Moreover, the experiments of [AMPZ.1] show that a classical, stable
small-amplitude
Riemann pattern not involving phase transitions may bifurcate
under an initial perturbation with compact support to a stable
large-amplitude pattern involving two or more phase-transitional layers; 
this is a stronger sense in which nucleation is seen to occur.

The present investigation is motivated by the third question,
concerning the possibility of
{\it pattern formation} in models of form \eqnref{1}.
Particularly intriguing are numerical experiments of
Frid\&Liu [FL.1--2,LF] and \v{C}ani\'c\&Peters [CP], 
in which Riemann solutions of various phase-transitional models
were seen to exhibit highly oscillatory,
Turing-like patterns reminiscent of Martensitic crystal
structure or nucleation in phase-transitional elasticity.
A natural conjecture, by analogy
with similar phenomena in Cahn Hilliard or reaction diffusion equations,
is that these are stable patterns consisting of fronts connecting various
periodic or constant states:
the classical mechanism for pattern formation.

Support for this point of view was given by recent investigations 
of \v{C}ani\'c suggesting a connection with Hopf bifurcation
and the appearance of limit cycles [C] in the traveling-wave
equations arising from the specific Riemann data associated with these
phenomena.
A closer look [AMPZ.4--5] reveals a rich global structure in the
phase portraits of the
associated traveling-wave ordinary differential equations, which indeed 
feature large-amplitude periodic orbits resembling 
each of the observed oscillatory motifs.
However, these periodic solutions appear in 
one-parameter families of varying wavelength, 
terminating in the infinite-wavelength limit at a hetero- or homoclinic cycle
corresponding to a degenerate Riemann pattern of front- or pulse-type 
traveling-waves with common speed:
a metastable, slowly interacting pattern of waves which by themselves
would be stable.
The criterion that would select a particular wavelength from among
such a family is unclear.
Moreover, in further numerical experiments [AMPZ.4--5], the patterns,
though apparently robust, do not appear to settle down into a final,
periodic configuration; rather, as conjectured in [AMPZ.4], 
they seem to be the result of complex metastable interactions of an 
infinite pattern of slowly interacting traveling fronts between
pure phases, driven by linear instability in the limiting states at
plus and minus spatial infinity (the ``generalized Riemann data'' in
the sense of [AMPZ.1]).

More precisely, it appears that oscillations in phase originate
at spatial infinity through a linear instability mechanism.
A stationary phase calculation [AMPZ.5] reveals that the response under
perturbation of the unstable constant states at infinity 
(i.e., the Green function 
of the constant-coefficient equations obtained by linearizing 
about constant solutions with those values)
is well approximated by a sum of modulated, time-exponentially 
growing Gaussian wave-packets
$$ 
\bar\Psi_j:=
\frac {
e^{\R \lambda(k_j^*) t} 
e^{ ik_j^*(x-\alpha_j^* t) }
e^{-(x-\alpha_j^* t)^2/4\beta_j^* t}}
{\sqrt{4\pi \beta_j^* t}}
\frac{
r_j l_j^*}
{\langle r_j,l_j\rangle},
\eqnlbl{unstablediff}
$$
where $k_j^*$ denotes the frequency
for which the associated dispersion relation $\lambda_j(k)$ takes on its
maximum real part (temporal growth rate), 
$i \alpha_j^*:= -d\lambda_j/dk(k_j^*)$, $\beta_j^*:=-(1/2)d^2\lambda_j/dk^2(k_j^*)$,
and $r_j$, $l_j$ denote right and left eigenvectors of the flux Jacobian
$df(u)$ evaluated at the background constant state.
The waves $\bar\Psi_j$
might be called {\it unstable linear diffusion waves}, by analogy with
the picture given by Liu\&Zeng [LZe] in the stable case.
Their oscillations grow exponentially in amplitude until they reach
the nonlinear regime, at which time they sharpen into slowly moving
fronts connecting approximately pure phases and do not grow further.

Note that this process does not involve periodic waves, or even 
their infinite-wavelength heteroclinic-cycle limit.
The average wavelength of the pattern is, rather, closely approximated
by the characteristic wavelength of the {\it linear} oscillations
originating at spatial infinity, i.e., the rate at which new waves
are ``born.''  Likewise, the ``front'' separating an oscillatory
region from an unstable constant state is just the front edge of
the linear Gaussian wave-packet, whose speed of propagation can be 
well-approximated by tracking the point at which the Gaussian 
envelope reaches a fixed amplitude representing the transition
to nonlinear dynamics.
For more detailed discussion, and a derivation of \eqnref{unstablediff},
see Appendix B.
This description, if correct, represents a novel and
nonclassical mechanism for pattern formation, different
from that seen in Cahn Hilliard and reaction diffusion
models for phase transition.
Indeed, it seems more related to certain models for turbulence, 
in which energy from high frequency modes drives the evolution
of characteristic large-scale structures on a
lower-dimensional attracting manifold.

It seems difficult to rigorously verify this picture of 
pattern formation, which by its nature lies outside of the 
usual analytical frameworks.
However, it is possible to give indirect support by eliminating the more
usual scenario involving periodic waves.
For, recall that the patterns of [FL.1--2,LF,CP], 
consist of one or more oscillatory regions sandwiched between two 
unstable constant regions extending to plus and minus spatial infinity.
If these oscillatory regions indeed represent pieces of different
periodic solutions, then the most likely scenario is
that at least one of these periodic solutions should be {\it stable};
the alternative, a stable pattern consisting entirely of unstable
pieces, would represent a new type of dynamic stabilization for which
we know of no possible mechanism.  
Thus, it is strong evidence against the classical scenario
if we can show that there
exist no large-amplitude {stable} periodic solutions
lying nearby the observed oscillatory patterns.

\bigskip
Motivated by these considerations, we here study the general
question of stability of periodic solutions of conservation
laws with viscosity.  Of particular interest is the situation present in the
above-mentioned numerical studies, of a planar, periodic family of
solutions, originating from a nonlinear center and bounded by a
limiting homoclinic or heteroclinic cycle, i.e., a typical (global)
Hopf bifurcation [GH,HK].  

Our analysis is by spectral Evans-function
techniques, using an analytic framework developed by R. A. Gardner [G].  
Specifically, changing coordinates to a rest frame for the traveling
wave $\bar u$, we obtain the linearized equation
$$
v_t=Lv:= (Bv_x)_x - (Av)_x,
\eqnlbl{1.3}
$$
about the (now stationary) wave $u=\bar u(x)$.
We shall investigate stability of $\bar u$ by the study of the spectrum
of the linearized operator $L$ about the wave.
More precisely, we investigate the {\it weak spectral stability} of
$\bar u(\cdot)$ as a solution of \eqnref{1}, defined as
$$
\sigma(L)\subset \{\R \lambda \le 0 \},
\eqnlbl{spectralstability}
$$
i.e., nonexistence of time-exponentially unstable eigenmodes.
Note that failure of \eqnref{spectralstability} implies exponential
linearized and (appropriately defined) nonlinear instability;
hence, the study of spectral stability is appropriate for 
investigation of {\it instability phenomena}.
(Linearized {\it stability} involves different issues, and
will be discussed in a companion paper [OZ]).

Following Gardner [G], we study the spectrum of $L$ ``directly,'' via
the {\it periodic Evans function}, a determinant $D(k,\lambda)$
involving the monodromy matrix of the linearized eigenvalue equation
$$
(L-\lambda)w=0
\eqnlbl{1evalue}
$$
for $L$, whose zero set $(k, \lambda)$ corresponds to bounded
solutions $e^{ikx}w(x)$ of \eqnref{1evalue}, with $w$ periodic;
for a detailed derivation, see Section 2.
This approach has the advantage of generality,
whereas more specialized analyses based on variational 
[M.1--2,LP] or Hamiltonian [Mi] structure may take better
advantage of the specific structure of the underlying 
evolution model.
For conservation laws \eqnref{1}, there does not seem
to be any such useful structure to which we may appeal.

Evans function techniques have been used successfully to study
the stability of traveling-front and -pulse solutions; however,
computations in the periodic case become considerably more
complicated.  Indeed, to our knowledge, ours is one of only
three such explicit computations that have been carried out 
using Gardner's framework, the others being low-frequency
expansion  of the Evans function $D(k,\lambda)$
in the wave number $k$, carried out by E. Eszter
[Es] in a singularly perturbed limit 
and by B. Sandstede and A. Scheel
[SS.1] in the large-amplitude (long-period) limit.
These previous analyses were obtained by different
techniques, in the somewhat different
reaction diffusion setting, and under hypotheses 
that do not apply here:
in particular, that zero be a simple eigenvalue of the linearized 
operator about the wave. 
However, the main distinction of the present analysis is 
that it is not in principle restricted to any type of limiting
case, giving useful stability criteria for waves of arbitrary
amplitude or type. (On the other hand, as we shall see, 
such limiting situations can be extremely helpful in the analytic 
evaluation of these criteria).
This feature distinguishes our results also from those
obtained in, e.g., [Ec,Ma.1--2,BMi.1-2,Mi]
by other than Evans function techniques,
all of which concern the limiting cases of a
large-amplitude bounding cycle or a small-amplitude constant solution.

\medskip
{\it Abstract result 1.}
Our main result is the development of a stability index
analogous to that obtained by Gardner\&Zumbrun [GZ] in
the traveling-front or -pulse context,
relating evolutionary stability to the dynamics of the 
traveling-wave ordinary differential equation.   
Specifically, we show that 
$$
\Gamma= \sgn \gamma \Delta \det df(u_-)\ge 0
\eqnlbl{result1}
$$
is necessary for stability with respect to periodic perturbations
of the same period as the background wave, 
where $\gamma$ 
is a transversality coefficient for the
traveling wave ordinary differential equation
$$
u'=B^{-1}(u)(f(u)-f(u_-)-s(u-u_-)),
\eqnlbl{1.1}
$$
$$
\barU(0)=\barU(T)=u_0,
\eqnlbl{1.2}
$$
and $\Delta:=\det (\partial \bar u_m/\partial u_-)$,
where $\bar u_m$ denotes the mass over one period
of the periodic profile $\bar u(\cdot)$ and $u_-\in \BbbR^n$
is an artificial parameter.
For a detailed derivation of \eqnref{1.1}, \eqnref{1.2} 
and 
precise definitions of $\gamma$ and $\Delta$,
see Sections 2 and 5, respectively.

Similarly as in the traveling-front or -pulse case [GZ,BSZ,ZS,Z.3], 
the coefficient $\Delta$
is seen to be related to well-posedness of an associated 
inviscid Riemann problem; however, this must now be
interpreted in the context of measure-valued solutions
appropriate for oscillatory solutions.
Alternatively, $\Delta\ne 0$ may be viewed as expressing the
existence and nonsingularity of a local ``mass map'' from
perturbation mass to potential time-asymptotic $T$-periodic
states, an evident necessary condition (by conservation of mass)
for orbital stability
of $\bar u$ with respect to nonzero mass perturbations within
the class of $T$-periodic solutions of \eqnref{1};
for further discussion, see Section 5.

\medskip
{\it Abstract result 2.}
The stability index $\Gamma$ detects 
strong (time-exponentially growing) instabilities analogous in the 
front or pulse case to unstable point spectrum of the linearized operator 
about the wave.
By closely related computations, we derive also a complementary
long-wave stability criterion
detecting weak (time-algebraically growing) instabilities,
analogous in the front or pulse case to unstable essential spectrum, 
or linearized instability of the limiting, constant states connected
by the profile of the traveling wave.
Whereas the stability index concerns
instability with respect to periodic perturbations with the same 
period as the background wave, 
the long-wave stability
criterion concerns instability
with respect to periodic perturbations on a different period, in
the limit as this period goes to infinity: i.e., 
so-called ``sideband instabilities''.
Specifically, in the ``quasi-Hamiltonian'' case that 
the traveling-wave ordinary differential equation \eqnref{1.1}
has an integral of motion,
we establish an illuminating small-frequency expansion 
$$
\aligned
D(k,\lambda)&= \gamma \lambda 
\det \big( -\lambda (\partial \bar u_m/\partial u_-)df(u_-)^{-1}
- ikT  \big) \\
&+ \CalO\big((|k|+|\lambda|)^{n+2}\big),
\endaligned
\eqnlbl{lf}
$$
of the Evans function,
where $\gamma$ and $(\partial \bar u_m/\partial u_-)$ 
are as described just above.
This yields a distinguished dispersion relation
$$
\lambda_0(k)=o(k)
\eqnlbl{dist}
$$
and $n$ dispersion relations 
$$
\lambda_j(k)=-i kT\alpha_j + O(k^2), \quad j=1,\dots,n,
\eqnlbl{disp}
$$
where $\alpha_j$ denote the eigenvalues of
$( \partial \bar u_m/\partial u_-)^{-1}df(u_-)$.
Relations \eqnref{disp} generalize those seen in the 
constant-coefficient case, for which 
$\partial \bar u_m/\partial u_-$ reduces to $I$;
for further discussion see Example 3.4 and Remarks 3.5 and 5.10.
In particular, they show that the ``generalized hyperbolicity''
requirement 
$$
\sigma \Big( ( \partial \bar u_m/\partial u_-)^{-1}df(u_-)\Big)\, 
\text{\rm real}
\eqnlbl{longwave}
$$
is necessary for long-wave stability.

{\bf Remark \thmlbl{analogy}.}
There is a further analogy between
conditions \eqnref{result1} and \eqnref{longwave} and
the stability index and long-wave stability conditions
arising in the study of multi-dimensional planar viscous
shock fronts [GZ,BSZ,ZS,Z.3], where, similarly, the stability
index concerns strong instability with respect to planar perturbations
respecting the symmetry of the background solution, while
the long-wave stability condition concerns weak instability
with respect to perturbations with small but nonzero transverse frequency.
However, the origins of the long-wave dispersion relations are rather
different in the two cases; in particular, \eqnref{lf} appears to be
closely tied to the assumed quasi-Hamiltonian structure, 
whereas the corresponding expansion in the shock front case is 
completely general.

\medskip

{\it Numerical results.}
Along with these analytic tools, we develop also numerical techniques
for approximating the spectrum of the linearized operator about
the wave, analogous to
those developed by Brin [Br.1--2,BZ] in the traveling-front or -pulse case.
Here, the issues are rather different.
For example, the evaluation of the Evans function in the periodic case is
numerically straightforward, since it involves only integration of a 
well-behaved ordinary differential equation on a finite interval;
by contrast, the traveling-front or -pulse
case involves integration on an infinite interval, making the problem
stiff.
On the other hand, the location of the spectrum, that is, the zero
set $(k, \lambda(k))$ of the Evans function $D(\cdot,\cdot)$ becomes
considerably more complicated than in the traveling-front or -pulse case,
for which the Evans function depends on a single argument only.
In particular, it is quite difficult to resolve the $n+1$ curves
\eqnref{dist}, \eqnref{disp}
bifurcating from the point $(k,\lambda)=(0,0)$
Thus, our numerical techniques are effectively restricted to 
$(k,\lambda)$ bounded away from the origin, and so are complementary
to the analytical techniques described above, which include
but are not limited to 
the Taylor expansion of $D$ about the origin.

{\it Applications.}
The stability index and long-wave stability criterion
are explicitly {\it comparable} in the same
``planar Hamiltonian'' 
case for which the index of Gardner\&Zumbrun was explicitly
{\it evaluable} in the traveling-front or -pulse case:
namely,
the case for which the traveling-wave ordinary 
differential equation \eqnref{1.1} is planar Hamiltonian for 
some distinguished speed,
and together yield rigorous results of instability similar to those
obtained previously for pulse-type solutions in [GZ,Z.1]
In this case, using a Poincar\'e-Bendixon argument similar to that 
used in [GZ], we may express the sign of the  Melnikov integral $\gamma$ 
appearing in the stability index 
in terms of the sign of the derivative
$dT/da$ of period $T$ with respect to amplitude $a$
in the embedding family of nearby periodic orbits of \eqnref{1.1},
as simply
$$
\sgn \gamma = \sgn dT/da.
\eqnlbl{signgamma}
$$
This reduces both of stability conditions \eqnref{result1} 
and \eqnref{longwave}
to requirements on the derivative $\partial \bar u_m/\partial u_-$, 
a quantity which in the small-amplitude and the large-amplitude 
homoclinic limits approaches the identity.

In general (i.e., for intermediate-amplitude waves),
$\partial \bar u_m/\partial u_-$ seems difficult if not impossible
to evaluate analytically.
Nonetheless, we are able to obtain rigorous
instability results through a novel {\it dichotomy}, 
asserting that criteria \eqnref{result1} and \eqnref{longwave}
are in certain cases logically exclusive
based only on structure or symmetries of 
$\partial \bar u_m/\partial u_-$ rather than its precise value.
The basic strategy is to show the 
matrix $df(u_-)(\partial u_m/\partial u_-)^{-1}$
to be trace-free, in which case its eigenvalues are real if and only
if its determinant is negative and thus $\Gamma = -\sgn dT/da<0$
if $dT/da>0$.
(Recall that $df(u_-)(\partial u_m/\partial u_-)^{-1}$
is $2\times 2$ in the planar case we consider.)
%
Specifically, we show for 
the equations of van der Waals gas dynamics and elasticity
that there exist {\it no} stable periodic solutions, of any amplitude,
with the property that period increases with amplitude in the 
one-parameter family of nearby periodic orbits: in particular
that large-amplitude waves are unstable;
a separate argument shows that small-amplitude solutions are
unstable, due to Majda-Pego instability of the limiting 
constant state [MP].  
In the case of a standard double-well potential, it can be shown
that the period is everywhere monotone increasing (see Remark 7.2, below), 
and so {\it all periodic solutions are unstable}.
Likewise, for the special (Hamiltonian) class of planar, quadratic-flux 
systems considered in [GZ]
we establish instability of large-amplitude waves lying near the
numerically observed patterns, and (by separate, continuity argument)
of small-amplitude solutions.

These two cases are quite relevant to the 
issue of pattern formation that originally motivated our investigations.
The former equations, modeling phase-transitional fluids/solids,
are precisely those that were seen to exhibit
pattern formation in the numerical investigations of [FL.1].
The latter, 
which serve as qualitative models for multi-phase flow in porous media, 
are prototypical for the other models in which
pattern formation has been observed [FL.2,CP,AMPZ.4--5].
Indeed, taken together, our results are strongly 
suggestive of nonexistence of {\it any} stable periodic waves 
in the cases where oscillatory patterns have been observed.
Further support for this conjecture is given
by our numerical experiments, in which we find instability
of periodic solutions across their entire (numerically determined) 
region of existence. 
Our results thus appear to eliminate from consideration
the usual paradigm of multiple stable periodic states, pointing
to a {\it different mechanism} for pattern formation in
multiphase conservation laws than that typically seen in 
reaction diffusion equations: for example, the one proposed in [AMPZ.4--5].  

Besides the direct, physical interest of our conclusions,
we point out an indirect contribution of this paper that
is perhaps more important.  Namely, we provide a useful,
and {\it explicit} analogy between the Evans function
framework for the periodic case and that of the better-studied
traveling-front or -pulse case, {\it different than the
large-period limit} studied by Gardner [G.1--2].
In the present paper, we have applied this analogy to the 
results of [GZ,Z.1,BSZ,ZS,Z.3] to obtain interesting instability results 
in the context of conservation laws with viscosity, i.e., 
{\it necessary} conditions for stability.
A similar translation of the complementary 
pointwise semigroup methods developed in [ZH] yields 
interesting {\it sufficient} conditions for linearized stability;
these results will be presented in the companion paper [OZ].

\medskip

{\it Discussion and open problems.}
%
The quasi-Hamiltonian assumption made in deriving
\eqnref{lf} seems to be fundamental for the
evaluation of the long-wave viscosity criterion.
Besides appearing frequently in physical examples, quasi-Hamiltonian
systems are shown in [OZ] to be the {\it only} type that can
support periodic waves that are asymptotically stable in the usual, diffusive
sense (see further discussion, beginning of Section 6).

A related assumption on the structure of
stationary solutions is made by Bridges and Mielke
[BMi.2] in their powerful study, by quite different, center manifold techniques,
of instability of {\it multidimensional} 
but small-amplitude periodic waves in the Cahn Hilliard-Allen Cahn setting
$$
m\Psi_{tt}+ d\Psi_t=  \Psi_{xx}+ \Psi_{yy} + dF(\Psi),
\eqnlbl{BMi}
$$
$F$, $\Psi$, $x \in \BbbR$, $y\in [0,\pi]$ with periodic boundary,
$m$, $d\ge 0$.
Namely, they assume that the center manifold of small bounded solutions
of the stationary equation near some rest point 
is foliated by periodic solutions
consisting of level sets of an appropriate ``spatial'' Hamiltonian:
that is, it has the same solution structure as does the 
(planar Hamiltonian) one-dimensional stationary equation
$$
\Psi_{xx} =dF(\Psi)
\eqnlbl{BMistat}
$$
obtained formally by shrinking the $y$-dimension to zero in \eqnref{BMi}.
Likewise, this structure is present, and used in an important way
in the authors' (related) landmark study of Benjamin-Feir instability 
in [BMi.1] and in Mielke's characterization, using a Lyapunov-Schmidt type
reduction, of linearized stability of roll-up solutions of the 
Swift-Hohenberg equation [Mi].

It is interesting to compare our results to those of Bridges and
Mielke in [BMi.2].  They derive a relation (see Lemma 4.1 of the reference) 
$$
\tilde D(k,\Lambda)= b_0 (\Lambda - C k^2) + O(|\Lambda|^2+|k||\Lambda|+ |k|^3)
\eqnlbl{reduced}
$$
similar to our \eqnref{lf},
with 
$$
b_0\ne 0, \quad \sgn C= \sgn dT/da.
\eqnlbl{signs}
$$
Here, $\Lambda:=m\lambda^2+d\lambda$ is the generalized spectral parameter
arising in the eigenvalue equation associated with the linearized
stability problem, 
and $\tilde D$ is essentially an Evans function
for the reduced eigenvalue equation they obtain by center manifold reduction.
{{From \eqnref{reduced}, \eqnref{signs} there follows immediately
the geometric necessary condition
$$
dT/da < 0
\eqnlbl{BMinec}
$$
for stability that is their main result:  equivalently,
a sufficient condition $dT/da>0$ for instability.
It is to be noted that {\it this is a sideband-type instability condition,
and not a stability index} in the sense of our condition \eqnref{result1}.

This result generalizes well-known results in the one-dimensional
case (obtained by dropping $y$ in \eqref{BMi})
relating $dT/da>0$ to instability, described, e.g., in [IR].
However, as pointed out in [BMi.2, discussion below (1.13)], 
$dT/da<0$ in the one-dimensional case
{\it also implies instability}, though of a different kind:
strong (large $\lambda$) rather than sideband instability.
Thus, it might well be that such instabilities occur also in 
the multidimensional case, even though they are not detected
by the sideband instability analysis of Bridges\&Mielke.

In this regard, it is interesting to consider the stability of 
periodic solutions of the one-dimensional version 
$$
m\Psi_{tt}+ d\Psi_t=  \Psi_{xx}+ dF(\Psi)
\eqnlbl{BMi1}
$$
of \eqnref{BMi} from the point
of view of this paper.
By exactly the same calculations used to establish
\eqnref{result1}, \eqnref{lf}, and \eqnref{signgamma},
we obtain for this problem a low-frequency expansion
$$
D(k,\Lambda)= \gamma \Lambda -  T^2 k^2 + O(|\Lambda|^2+|k||\Lambda|+ |k|^3)
\eqnlbl{unreduced}
$$
and a stability index 
$$
\Gamma= \sgn \gamma= \sgn dT/da;
\eqnlbl{ind}
$$
indeed, the calculations reduce substantially in this case.
Here, $D$ as usual denotes the Evans function for the 
full (unreduced) eigenvalue equation associated with linearized stability,
and $\Lambda:=m\lambda^2+d\lambda$ as above denotes
the generalized spectral parameter arising in that equation.
Thus, we not only recover the sideband instability analysis 
\eqnref{reduced}--\eqnref{BMinec} of
[BMi], but, through the stability index, detect also strong
instabilities in case $dT/da<0$.
That is, we see a similar dichotomy
in this substantially different setting to the one used here to
show instability of solutions of \eqnref{1}; in fact, 
the absence of the ``hyperbolic factor'' 
$\Delta=\partial \bar u_m/\partial u_-$ in the formulae 
makes the argument much simpler to apply in this
case, and leads to a stronger result
(instability without regard to  $\sgn dT/da$).
Note that there is no requirement here on the amplitude of solutions.


The above discussion suggests that an Evans function analysis
based on ``spatial dynamics'' might be a possible direction for 
generalization of the multidimensional analysis of Bridges\&Mielke,
both to the case $dT/da<0$ and to the case of large-amplitude waves.
The Evans function described by Sandstede and Scheel [SS.2],
based on Galerkin approximation on finite subspaces, appears 
to be a natural candidate for such investigations.
This would be an extremely interesting direction for further study.

The analysis of \eqnref{BMi1} also suggests interesting questions
in the one-dimensional case.  Our results for the vectorial
van der Waals phase-transitional model are now seen to be natural,
though somewhat weaker,
generalizations of those for the scalar phase-transitional model \eqnref{BMi1}.
This begs the question whether one might obtain
a weakened version of these results
also for the vector version $\Psi \in \BbbR^n$ of \eqnref{BMi1}.
Indeed, one obtains for more general equations 
$mu_{tt}+ du_t=  u_{xx}+ G(u)$, $u\in \BbbR^n$,
in the quasi-Hamiltonian case, 
a low-frequency expansion
$$
D(k,\Lambda)= \gamma \Lambda -  
\delta T^2 k^2 + O(|\Lambda|^2+|k||\Lambda|+ |k|^3),
\eqnlbl{vreduced}
%
$$
and a stability index $\Gamma=\sgn \gamma$, where 
the transversality coefficient $\gamma$ no longer has
a simple geometric interpretation, and the
coefficient $\delta $ is no longer explicitly evaluable.
Thus, we do seem to obtain a result of instability when $\delta >0$,
analogous to that obtained for the van der Waals model when
$dT/da>0$.  However, it is not clear whether $\delta >0$ can
in fact occur;  this would be another interesting issue for further
investigation.
Likewise, in the conservation law setting, it would be interesting to
determine whether or not the results for the van der Waals model carry
over to the general planar Hamiltonian case, or even to the vector
quasi-Hamiltonian case, {\it without} further assumptions restricting
the structure of the matrix $\partial \bar u_m/\partial u_-$.
We suspect strongly that the answer is ``no,'' but do not so far have any
counterexamples.


Finally, we mention the related analysis carried out 
by Laugesen and Pugh [LP] using quite different, variational methods, 
of stability under periodic perturbations of periodic solutions
(with the same period) 
of the thin film evolution model 
$$
h_t = (f(h) h_{xx} + g(h)h)_{xx} + ah,
\eqnlbl{thinfilm}
$$ 
$h\in \BbbR$, in the zero-gravity case $a=0$, which as far as we know is the
only other treatment of stability of periodic waves in the
conservation law setting (they obtain results also for $a\ne0$,
but these are not relevant to our discussion).
Here again the stationary-wave equation is a second-order
nonlinear oscillator, so planar Hamiltonian, even
though the general (nonzero speed) 
traveling-wave equation is third-order.
Applying our methods to this problem, we conjecture that
one should be able to recover a partial version of their results, 
namely a stability index of $\Gamma=\sgn (dT/da)$, 
where the derivative $dT/da$ is taken with the area under $h$
(i.e., the mass) held fixed, yielding instability whenever
$dT/da<0$.\footnote{
The situation here is somewhat degenerate due to the fact that the 
traveling-wave equation reduces in order at speed zero; 
indeed, condition (H3) of section 2 is violated.
However, the computation can still be carried out, as discussed
in Remark 5.3 below.
}
(Laugesen and Pugh in fact {\it characterize} stability by this condition,
obtaining positive stability results as well.)
What is more interesting, one might equally well hope to carry out
a low-frequency expansion analogous to \eqnref{lf}, to obtain
conditions for sideband instability as well.  These were not treated
in [LP], and do not seem to be accessible by their methods.

{\bf Remark \thmlbl{SS}.}
The planar Hamiltonian structure of \eqnref{BMi1} is what makes possible
the explicit evaluation in \eqnref{unreduced}
of the coefficient of the quadratic term in $k$.
For systems of general type, this can usually be done only
in some asymptotic limit;  see for example the analyses of
[Es,SS.1] in the large-amplitude limit.
\medskip

{\bf Plan of the paper.} 
In Section 2, we frame the problem, describing the equations and 
assumptions under consideration, and deriving traveling-wave and
eigenvalue ordinary differential equations.
In Sections 3 and 4, we define the Evans function, following [G.1--2],
and recall a result of [G.2] relating the
spectra of periodic waves in the large-amplitude (large-period) limit
to that of a bounding homoclinic; we also explore briefly the more elementary
small-amplitude limit, in which periodic waves approach a limiting, 
constant solution. 
In Section 5, we define a periodic stability index analogous to that
of the traveling-front or -pulse case, and establish the
fundamental relation \eqnref{result1}.
In Section 6, restricting to the quasi-Hamiltonian case, 
we derive by similar techniques the small-frequency expansion \eqnref{lf}.
This yields, in particular, 
an appealing formula for the ``averaged,'' or effective 
constant-coefficient equation governing behavior under perturbation,
and also gives further details completing the description
of large-amplitude behavior given by Gardner in [G.2].
Finally, restricting further to the planar Hamiltonian case, and applying
a Poincar\'e-Bendixon argument similar to that used in [GZ], 
we obtain the sign of the  Melnikov integral $\gamma$ 
appearing in the stability index in terms of the sign of the derivative
$dT/da$ of period $T$ with respect to amplitude $a$
in the embedding family of nearby periodic orbits;
In Section 7, we use the analytic tools developed in
Sections 5 and 6 to establish rigorous instability results for
the two classes of example systems described above.
In Appendix A, we describe a numerical algorithm for location of
the spectrum, and carry out systematic numerical experiments 
for the same two classes of example system;
these support and in some cases extend our earlier analytical results.
Finally, in Appendix B, we describe in detail the alternative
mechanism for pattern formation proposed in [AMPZ.4--5], based
on metastable configurations of slowly interacting fronts, driven
by linear instability in the constant states at spatial infinity.
As described in the appendix, these could be thought of as
{unstable nonlinear diffusion} waves generalizing the stable
versions described in [LZe].

\bigskip

\sectionnumber=2 
\theoremnumber=0
\equationnumber=0
\smallskip
\TagsOnLeft
\specialhead
\S2.\quad 
Preliminaries
\endspecialhead

Consider a conservation law \eqnref{1}
and a periodic traveling-wave solution $u=\barU(x-st)$, of period
$T$, satisfying the traveling-wave ordinary differential equation
$$
(B(u)u')'= f(u)' -su',
\eqnlbl{second_order}
$$
with initial conditions
$$
\aligned
\bar u(0)&=\bar u(T)=:u_0,\\
\bar u'(0)&=\bar u'(T)=:u_1.\\
\endaligned
\eqnlbl{init}
$$
Here, and elsewhere, `$'$' denotes $\partial/\partial x$.
Integrating \eqnref{second_order}
from $0$ to $x$, we obtain a first-order dynamical system
$$
u'=B^{-1}(u)(f(u)-su-q),
\eqnlbl{1.1a}
$$
$$
\barU(0)=\barU(T)=u_0,
\eqnlbl{1.2a}
$$
parametrized by $(q,s)\in \Bbb{R}^{n+1}$, where the ``total flux''
$$
\aligned
q:= &B(u)u'- f(u) + su \\
\equiv &B(u_0)u_1 - f(u_0) + su_0\\
\endaligned
$$
is a constant of motion.  Notice that the map $(u_0,u_1)\to (u_0,q)$
is locally invertible so long as $\det B(u_0)\ne 0$, by the 
Inverse Function Theorem, hence we have lost no information
by this reparametrization.

To emphasize the connection with the traveling-front or -pulse case, 
we assume further that there exists some rest point $u_-$ of
\eqnref{1.1a}, i.e., $q=f(u_-)-su_-$, so that we can rewrite
\eqnref{1.1a}, \eqnref{1.2a} in the final form \eqnref{1.1},
\eqnref{1.2} given in the introduction.
We shall use this form of the equations throughout the paper.
Recall that \eqnref{1.1} is exactly the equation satisfied by
a front or pulse type traveling-wave solution
$$
u=\bar u(x-st), \quad \lim_{x\to \pm \infty}\bar u= u_\pm,
\eqnlbl{front}
$$
allowing the convenient comparison of front or pulse type
and periodic solutions within the same dynamical system framework.
Note, in the planar case of our main interest, that such a rest
point always exists within the region bounded by a periodic wave,
so there is no loss of generality in changing to the new coordinates;
moreover, provided that $u_-$ is a nondegenerate rest point,
$\det (df(u_-)-sI)\ne 0$, the map $u_-\to q$ is again locally invertible.

We make the following nondegeneracy assumptions,
analogous to those made in [GZ,ZH] for the traveling-front or -pulse case:

\quad ({H0})  $f$, $B\in C^2$.
\smallskip

\quad ({H1})  $\R \sigma(B(\bar u(x)))>0$ for all $x$.
\smallskip

\quad ({H2})  $df(u_-)-sI$ {invertible}.
\smallskip

\quad ({H3})  
The traveling-wave profile $\bar u$ is a transversal orbit
of traveling-wave equation \eqnref{1.1} 
augmented with $s'=0$, under periodic boundary conditions $u(0)=u(T)$;
in particular, for fixed $u_-$, $\barU$ 
is, locally, the unique $T$-periodic solution of 
\eqnref{1.1} up to translation, even allowing 
variation in the speed of propagation $s$.

\smallskip

\noindent As mentioned above, ({H2}) concerns 
nondegeneracy of the parametrization $(u_0,u_-,s)$ of \eqnref{1.1},
and thus is not strictly necessary for our analysis of
periodic waves below; a similar analysis without
this assumption can be carried out in $(u_0,q,s)$ 
coordinates, with only minor expositional changes.

{\it Note:}  Assumption ({H3}) {\it does not preclude}
the interesting case of a planar 
Hamiltonian ordinary differential equation possessing a family of 
nested periodic solutions, since generically the
orbits of the periodic family have distinct periods; this will in 
fact be the main source of our examples in Sections 6--7.

Now, assume without loss of generality that speed $s=0$,
i.e., $u\equiv\bar u(x)$ is a {\it stationary} solution of
\eqnref{1}.
Linearizing \eqnref{1} about $\barU(\cdot)$, we obtain the usual 
linearized equation \eqnref{1.3},
where now $A,B$ are {\it periodic}, rather than asymptotically
constant as in the traveling-front or -pulse context.  
The eigenvalue equation for $L$ is likewise
$$
(Bw')'=(Aw)'+\lambda w,
\eqnlbl{2}
$$
where again `$'$' denotes $\partial/\partial x$,
or, written as a first-order system:
$$
W'=\BbbA(\lambda ,x)W, 
\eqnlbl{4}
$$
where $W:= (w,w')^t$ and coefficient
$$
\BbbA:=
\pmatrix
0 &I\\
B^{-1}(\lambda I+A') & B^{-1}(A-B')
\endpmatrix \eqnlbl{4.1}
$$
is a periodic $2n \times 2n$ matrix.

Through a study of the eigenvalue equations,
we shall investigate the {\it weak spectral stability} of
$\bar u(\cdot)$ as a solution of \eqnref{1}, 
as defined in \eqnref{spectralstability}.
That is, we shall investigate whether or not the linearized operator $L$ 
possesses unstable eigenmodes $\lambda: \, \R \lambda >0$.
Failure of \eqnref{spectralstability} implies
{\it exponential linear instability} with respect to test-function
initial data, as measured in any norms whatsoever; this can be
seen by applying the evolution operator to $\phi \chi^M$ for $M$ 
sufficiently large, where $\phi$ is a (merely bounded) unstable mode and 
$\chi^M(x):=\chi(x/M)$, with $\chi$ a smooth cutoff function 
that is $1$ on $[-1,1]$ and vanishes off of $[-2,2]$.
By the ``almost-finite propagation speed'' property of \eqnref{1.3},
we find that the amplification of the resulting solution $v(t)$ in going
from time zero to time $t$ is of order
$$
\|v(t)\|_{L^p}/\|v(0)\|_{L^q}\sim
t^{1/p-1/q}e^{\R \lambda t} \to \infty
$$
as $t\to \infty$, precluding uniform $L^q \to L^p$ stability
for any choice of $p$ and $q$.

{\bf Remark \thmlbl{strongstab}}.
It is shown in [OZ] that {\it strong spectral stability,}
defined as \eqnref{spectralstability} augmented with appropriate
nondegeneracy conditions (conditions (D1)--(D3) of [OZ]),
implies linearized $L^1 \to L^p$ asymptotic stability
for all $p>1$, with uniform rates of decay equal to those
for the standard heat equation.

\bigskip

\sectionnumber=3
\theoremnumber=0
\equationnumber=0
\smallskip
\TagsOnLeft
\specialhead
\S3.\quad 
The Evans function and the Spectrum of $L$
\endspecialhead

A brief calculation reveals that $L$ has no point spectrum in $L^p, 
p < \infty$:
Following [G.1--2], we introduce the
{\it monodromy matrix} 
$$
M(\lambda ):= \Psi(T,\lambda),
\eqnlbl{3}
$$
where $\Psi(\cdot,\lambda)$ is the fundamental solution of \eqnref{4}, i.e.,
$$
\Psi'=\BbbA(\lambda ,x)\Psi,\quad \Psi(0,\lambda)= I
\eqnlbl{4.2}
$$
and $T$ is the period of the coefficients.  Then,
$$
W(NT)=M(\lambda )^N W(0)
\eqnlbl{5}
$$
for any integer $N$, for any solution $W$ of the eigenvalue equation
\eqnref{4}, whence $W$ can be at most bounded and not decaying
at $\pm \infty$.  Indeed, we have: 

\proclaim{Proposition \thmlbl{1}}  $\sigma(L)$ consists, for all
$L^p$, precisely of $L^\infty$ eigenvalues, i.e., $\lambda $ such that
$$
\det(M(\lambda )-\gamma )=0,
\quad |\gamma |=1.
\eqnlbl{6}
$$
\endproclaim

{\bf Proof}.  If $M(\lambda )$ has no eigenvalue of modulus $1$, then
there exist $k$ ``stable'', i.e., with modulus less than $1$, eigenvalues, and
$2n-k$ ``unstable,'' or with modulus greater than $1$, eigenvalues, 
the associated normal modes decaying
exponentially at $+\infty$ and $-\infty$, respectively,
and linearly independent.  We
can thus construct a Green function $G_\lambda $, as in the
asymptotically constant case, 
to obtain a bounded resolvent, hence $\lambda \in \rho (L)$.
For further details, see [He,Z.3--4,OZ].

On the other hand, $L^\infty$ eigenvalues can be shown to lie in any
$L^p$ spectrum, by a standard limiting argument: specifically,
by showing that 
$$
\|\chi^M\phi\|_{L^p}/ \|(L-\lambda)\chi^M \phi \|_{L^p}\to \infty
$$
as $M\to \infty$ for any bounded $\phi$ such that $(L-\lambda)\phi=0$,
where $\chi^M(x):=\chi(x/M)$ for any $C^\infty$ cutoff function $\chi$
that is one on $[-1,1]$ and vanishes off of $[-2,2]$.
\qed 
\bigskip

Loosely following [G.1--2], we define the {\it Evans function}
$$
D(k,\lambda):= \det (M(\lambda) - e^{ikT})
\eqnlbl{evans}
$$
for any $(k, \lambda)\in \BbbR\times \BbbC$.
Note that $D(k, \lambda)$ is clearly jointly analytic 
in $k$ and $\lambda$ on all of $\BbbR \times \BbbC$; thus, 
it is somewhat better behaved than the corresponding object
in the traveling-front or -pulse case, see [GZ,ZH] for further discussion.
{{From Proposition \thmref{1}, we have immediately:

\proclaim {Corollary \thmlbl{Dcor}}
The spectrum of $L$ consists of the set of all $\lambda$
such that $D(k,\lambda)$ vanishes for some $k$.
\endproclaim

In the terminology of [G.1--2], points $\lambda$ satisfying
\eqnref{6} are called $\gamma$-eigenvalues;
we will call them $k$-eigenvalues, where $\gamma=e^{ikT}$.
The parametrization by $k$ gives a more transparent analogy to
the analysis by Fourier transform of the constant-coefficient case,
see Example 3.4 just below (or, see [OZ] or [S.1--3] for
a deeper discussion of this analogy).

{\bf Remark \thmlbl{1.1}}.
More generally,
$$
M(\lambda )=\TPsi(L) \TPsi^{-1}(0),
\eqnlbl{6.1}
$$
for any matrix $\TPsi$ of solutions of \eqnref{4}.  The coordinate-independent 
representation \eqnref{6.1} is quite useful in computations.
Likewise, $D(\cdot,\cdot)$ is invariant under
linear changes of coordinates,
$$
\Psi \to P(x) \Psi,
\eqnlbl{6.2}
$$
with $P$ periodic.
({Note:}  This includes changes of coordinates at the level of
\eqnref{2}, but also more general ones at the level of phase
coordinates, \eqnref{4}).

\medskip

It is interesting to compare to the spectrum of $L$ considered as an
operator on periodic functions $L^2[0,T]$.  
Necessarily discrete, the
spectrum of $L$ acting on the bounded interval $[0,T]$ 
corresponds precisely to the set of 
($k=0$)-eigenvalues of $L$ acting on the unbounded domain;
conversely, $(k=0)$-stability corresponds to stability with respect 
to periodic perturbations of period $T$.
More generally, $(k=2\pi m/n)$-stability
corresponds to stability with respect to perturbations that are periodic 
with period $nT$; letting $n\to \infty$, we recover stability with respect
to {\it general perturbations}.
Applying the Implicit Function Theorem to \eqnref{6}, 
we find that the spectrum of $L$ 
on the real line consists typically of {\it curves} of $k$-spectra, 
given by the closure of the spectra of $L$ considered as an operator
on periodic functions over all multiple intervals $[0,nT]$, $n$ an integer.

\medskip

{\bf Example \thmlbl{cc}.}
In the constant-coefficient case, $Lv:=Bv_{xx} - Av_x$ with
$A$, $B$ identically constant, an elementary computation yields 
$$
D(k,\lambda)=
\Pi_{j=1}^{2n} (e^{\mu_j(\lambda)T}-e^{ikT}),
$$
where $\mu_j$, $j=1,\dots,2n$ denote the roots of the characteristic equation
$$
\det(B\mu^2-A\mu-\lambda)=0
\eqnlbl{char}
$$
associated with eigenvalue equation \eqnref{4}.
On the set of $\lambda$ such that $\BbbA$ is diagonalizable, this is easily seen
by changing coordinates to a basis in which it is diagonal, see Remark
\thmref{1.1}; recalling that this set is dense, we obtain the full 
result by continuity.
Thus, the zero-set of $D$ consists of all $k$, $\lambda$ such that
$$
\mu_j(\lambda)=ik \, (\text{\rm mod } 2\pi i/T)
\eqnlbl{zeroset}
$$
for some $j$.
Setting $\mu=ik$ in \eqnref{char}, we obtain the dispersion relation
$$
\det(\lambda -ikA-k^2B)=0,
\eqnlbl{dispersion}
$$
recovering the standard characterization of $\sigma(L)$ by Fourier transform.
\qed
\bigskip
\medskip

Note that $D(0,0)=0$ in the example above, for any choice of $A$, $B$,
since $\mu=0$ is an $n$-fold root of \eqnref{char} at $\lambda=0$;
we shall see later that this holds also in the general variable-coefficient case.
An obvious necessary condition for stability is thus
$$
(\partial/\partial k) \lambda _*(0)\quad \hbox{ is imaginary},
\eqnlbl{8}
$$
where $\lambda _*(k)$ is
any smooth root of $D(k,\lambda_*(k))\equiv 0$ 
bifurcating from $\lambda _*(0)=0$.  
For, otherwise there would
exist exponentially unstable modes $\lambda=\lambda_*(k)$
for $k$ sufficiently small (in fact for all $k>0$ or all $k<0$, depending
on $\sgn \R (\partial/\partial k) \lambda _*(0)$).
We shall make essential use of this simple observation later on.

{\bf Remark \thmlbl{hyp}.}
In the constant-coefficient case,  relation \eqnref{dispersion} 
yields expansions
$$
\lambda_j( k)=0 - ia_j k+ \cdots, \quad j=1,\dots, n,
\eqnlbl{ccexp}
$$
for the $n$ roots bifurcating from $\lambda(0)=0$, 
where $a_j$ denote the eigenvalues of $A$.
Thus, we obtain the necessary stability condition of {\it hyperbolicity},
$\sigma(A)$ real.
\medskip

\bigskip

\sectionnumber=4
\theoremnumber=0
\equationnumber=0
\smallskip
\TagsOnLeft
\specialhead
\S4.\quad 
A result of Gardner
\endspecialhead

Before beginning our main analysis, we recall an interesting related
result of Gardner [G.2] concerning a particular (rather typical) 
kind of large-amplitude limit for a family of periodic waves.  
Specifically, consider the situation of a family
of periodic traveling waves $\{\bar u^\varepsilon\}$, approaching as $\varepsilon\to0$
to a limiting homoclinic, or ``pulse-type,'' solution $\bar u^0$,
the vertex of which is a nondegenerate (i.e., hyperbolic) rest point
of the associated traveling-wave ordinary differential equation.
Here, $\bar u^\varepsilon$ may be taken without loss of generality to be
stationary (i.e., zero speed), and the partial differential equations
they solve to vary smoothly
with $\varepsilon$, so long as $\varepsilon \to 0$ is a regular perturbation,
with all $\varepsilon$-equations parabolic of some fixed order.
Included in this framework is the usual case of traveling-wave 
solutions of the {\it same} (parabolic) partial differential equation, 
with possibly differing speeds,
approaching a homoclinic separatrix in the $\varepsilon\to 0$ limit.

In this scenario, as $\varepsilon\to 0$,
the period $T^\varepsilon$ of $\bar u^\varepsilon$ goes to infinity, while
the profile on a single period $[0,T^\varepsilon]$ approaches uniformly
to an appropriate shift of the homoclinic profile $\bar u^0$,
without loss of generality to $\bar u^0(x-T/2)$.
Thus, it is natural to ask whether the stability properties of
periodic waves may be related in the large-period limit to
those of the limiting homoclinic wave.
The following result of Gardner shows that this is indeed correct.

Assume as is standard [GZ,ZH] that the linearized operator $L^0$ 
about the homoclinic wave $\bar u^0$ has essential spectrum contained 
in $\R \lambda \le 0$, so that point spectrum determines spectral
stability, \eqnref{spectralstability}.  
This point spectrum may be detected by the vanishing of a traveling-front-type 
Evans function
$$
D(\lambda):=\det (\phi_1^-,\dots,\phi_k^-,\phi_{k+1}^+,\dots,\phi_N^+),
\eqnlbl{fronttype}
$$
where $\{\phi_1^-,\dots,\phi_k^-\}$ and $\{\phi_{k+1}^+,\dots,\phi_N^+\}$
are appropriately chosen bases of the subspaces of decaying solutions
at $-\infty$, $+\infty$, respectively, of the eigenvalue equation \eqnref{4}.
This can be constructed in a way that is analytic in $\lambda$ on all of 
$\{\lambda:\, \R \lambda >0\}$; for details, see [GZ,ZH].

Then, we have:

\proclaim{Proposition \thmlbl{2} [G.2]}  The homoclinic limit
$\bar u^0$ is weakly spectrally stable if and only if  the family
of periodic waves $\{\bar u^\varepsilon\}$ 
is ``weakly spectrally stable in the limit,'' in the sense that
the spectra of the associated linearized operators $L^\varepsilon$
about $\bar u^\varepsilon$ satisfy
$$
\sigma (L^\varepsilon) \subset 
\{\lambda: \, \R \lambda \le \eta(\varepsilon)\},
$$
with $\eta \to 0$ as $\varepsilon \to 0$.
In particular, if the homoclinic $\bar u^0$ is spectrally unstable,
then so are all $\bar u^\varepsilon$ for $\varepsilon$ sufficiently small.
\endproclaim

More precisely:

\proclaim{Proposition \thmlbl{limit} [G.2]}  
Let $\Lambda\subset \BbbC$ be compactly contained in 
$\{ \lambda: \, \R \lambda >0\}$, with $\partial \Lambda$
contained in the resolvent set of $L^0$.
Then, for $\varepsilon$ sufficiently small, the number of zeroes in $\Lambda$
of the periodic Evans function $D(\cdot,k)$, for any fixed $k$, is equal
to the number of zeroes of the front-type Evans function $D(\cdot)$.
\endproclaim

Proposition \thmref{2} is an immediate consequence of 
Proposition \thmref{limit} together with standard sectorial
bounds [Pa,He,Z.3--4] restricting the spectra of elliptic operators 
in $\{\lambda: \, \R \lambda \ge 0\}$ to a bounded domain.
Proposition \thmref{limit} is of more general application,
applying (under suitable, mild assumptions) also to evolution 
equations that are not parabolic.
{{From Proposition \thmref{limit} we obtain the revealing 
picture that eigenvalues of $L^0$, as $\varepsilon$ increases
from zero, unfold into loops of (essential) spectrum of $L^\varepsilon$;
looking in reverse, we see that, as $\varepsilon \to 0$,
the radius of the corresponding spectral loop decreases to zero
({\it Proof.} Denoting the eigenvalue in question as $\lambda_*$,
apply Proposition \thmref{limit} on balls $B(\lambda_*,r)$,
with $\varepsilon(r)\to0$ as $r \to 0$.)

In [G.2], (an equivalent version of) Proposition \thmref{limit} 
was established using the topological index (Chern number) construction 
introduced in [AGJ].  Here, we sketch an alternative, more elementary proof.
{On an initial reading, the reader may wish to skip this proof, 
which is independent from and uses a different set of techniques
than the analysis in the rest of the paper.}
\medskip

{\bf Proof of Proposition \thmref{limit}}.  
We sketch the proof in the context of \eqnref{1},
in which case $k=n$ and $N=2n$ in \eqnref{fronttype}, where 
$n$ as in \eqnref{1} is the dimension of $u$.
Analogous computations hold in the general case.
Since both periodic- and front-type Evans functions are analytic
in $\lambda$ on $\bar \Lambda \subset \{\lambda: \, \R \lambda >0\}$
the result will follow (using Rouch\'e's Theorem) by
a winding number calculation about the contour $\Gamma=\partial \Lambda$, 
provided we establish:

{\it Assertion \thmlbl{lim}}.  
On any compact subset $\Gamma$ of the resolvent set of $L^0$
intersect $\{\lambda:\, \R\lambda >0\}$, 
$$
e^{-\alpha_-T^\varepsilon}(-e^{ikT^\varepsilon})^{n}
D(k,\lambda_*)
\to C(k,\lambda) D(\lambda_*),
$$
holds as $\varepsilon\to0$,
uniformly in $k$, $\lambda$, where $C(k,\lambda)$ is a
nonvanishing (jointly) analytic function, $\alpha_-(\lambda)$
is analytic, and $n$ as above is the dimension of $u$.

{\it Proof of Assertion}.  
By the general theory of [AGJ,GZ,ZH], the basis elements
$\phi_j^+$ and $\phi_j^-$ approach exponentially as $x\to +\infty$,
$-\infty$, respectively, to the stable and
unstable subspace of the limiting coefficient matrices
$\BbbA_\infty:=\BbbA(\pm\infty)$
of \eqnref{4} (recall that in the homoclinic
case coefficients are not periodic, but approach limiting values
at the same rate as does the background wave $\bar u^0$:
exponential, in the case assumed here that $\bar u^0(\pm \infty)$ 
is a nondegenerate rest point of the traveling-wave ordinary differential equation).
More precisely,
$$
\phi_1^- \wedge \cdots \wedge \phi_n^-
\sim
\bar\phi_1^- \wedge \cdots \wedge \bar\phi_n^-
$$
and
$$
\phi_{n+1}^+ \wedge \cdots \wedge \phi_{2n}^+
\sim
\bar\phi_{n+1}^+ \wedge \cdots \wedge \bar\phi_{2n}^+,
$$
where $\cdot \wedge \cdots\wedge \cdot$ denotes an $n$-fold exterior
algebraic product, or minor, and $\bar \phi_j^\pm$ denote basis
elements for the manifolds of decaying solutions at $\pm \infty$
of the limiting, constant coefficient equations
$$
W'=\BbbA_\infty W
$$
at infinity.
In turn, we have
$$
\aligned
\bar\phi_1^- \wedge \cdots \wedge \bar\phi_n^-
&=
e^{\alpha_- x} V_1^- \wedge \cdots \wedge V_n^-, \\
\bar\phi_{n+1}^+ \wedge \cdots \wedge \bar\phi_{2n}^+
&=
e^{\alpha_+ x}
V_{n+1}^+ \wedge \cdots \wedge V_{2n}^+,
\endaligned
$$
where vectors $V_j^-$ and $V_j^+$ respectively span the unstable and stable 
subspaces of $\BbbA_\infty$, and $\alpha_-$ and $\alpha_+$
denote the trace of $\BbbA_\infty$ on these respective subspaces.

Now, let the matrix of solutions $\tilde \Psi$ in \eqnref{6.1}, 
Remark \thmref{1.1}, be chosen by
$$
\tilde \Psi(0):=(\tilde \psi_1^-(0),\dots \tilde \psi_n^-(0),
\tilde \psi_{n+1}^+(+T/2), \dots, \tilde \psi_{2n}^+(+T/2)),
$$
where $\tilde \psi_j^\pm$ are determined by
$$
\tilde \psi_j^-(0):=\bar\phi_j^-(-T/2),
\qquad
\tilde \psi_{j}^+(T):= \bar\phi_j^+(+T/2).
$$
Then, the standard estimates of [GZ,ZH] yield that 
$\det \tilde\Psi(\lambda,T/2) \to D(\lambda)$ uniformly in $\lambda$
as $\varepsilon \to 0$.

Moreover, since $D(\lambda)$ was assumed not to vanish
on $\Gamma$, we have that
$$
\phi_{1}^-(x)\wedge \cdots\wedge \phi_{n}^-(x)
\sim
C_-(\lambda)
e^{\alpha_+ x}
V_{n+1}^+ \wedge \cdots \wedge V_{2n}^+,
$$
as $x\to +\infty$,
which says that decaying solutions at $-\infty$ must be 
growing at $+\infty$, and likewise
$$
\phi_{n+1}^+(x)\wedge \cdots\wedge \phi_{2n}^+(x)
\sim
C_+(\lambda)
e^{\alpha_- x}
V_{n+1}^+ \wedge \cdots \wedge V_{2n}^+
$$
as $x\to -\infty$; the transmission coefficients $C_\pm$ as in [ZH] 
are determined by solving appropriate systems of linear equations 
with analytic coefficients, whose respective determinants factor 
as $D(\lambda)$ times a nonvanishing analytic function.
Therefore, 
$$
\tilde \psi_{1}^-(T)\wedge \cdots\wedge \tilde \psi_{n}^-(T)
\sim
C_-(\lambda)
e^{\alpha_+ T/2}
V_{n+1}^+ \wedge \cdots \wedge V_{2n}^+,
\eqnlbl{rep-}
$$
and likewise
$$
\tilde \psi_{n+1}^+(0)\wedge \cdots\wedge \tilde \psi_{2n}^+(0)
\sim
C_+(\lambda)
e^{-\alpha_- T/2}
V_{1}^- \wedge \cdots \wedge V_{n}^-.
\eqnlbl{rep+}
$$

Using the representation 
$D(k,\lambda)=
\det \big(\tilde \Psi(T)- e^{ikT}\tilde \Psi(0)\big)/\det \tilde \Psi(0)$
afforded by Remark \thmref{1.1}, combined with \eqnref{rep-}, \eqnref{rep+}
and the fact that $\alpha_+<0<\alpha_-$ (more precisely, $\phi_1^+,\dots,\phi^+_n$
decay exponentially, while $\phi_{n+1}^-,\dots \phi_{2n}^-$ grow exponentially),
we therefore obtain
$$
\aligned
D(k,\lambda)&\sim
\det (\tilde\psi_1^-(T),\dots,\tilde \psi_n^-(T),
-e^{ikT}\tilde \psi_{n+1}^+(0),\dots,-e^{ikT}\tilde\psi_{2n}^+(0))/
\det \tilde \Psi(0)\\
&\sim
C_-C_+ (-e^{ikT})^ne^{(\alpha- -\alpha_+)T}/C_+e^{-\alpha_+T}\\
&=
C_-(-e^{ikT})^ne^{\alpha-T)}.
\endaligned
\eqnlbl{perasymptotic}
$$

On the other hand, Abel's formula and the exponential
convergence $\BbbA\to \BbbA_\infty$ as $x\to \pm \infty$ imply that
$$
\aligned
D(\lambda)
&\sim
\det \tilde \Psi(T/2)\\
&=
C_0 e^{-(\alpha_- + \alpha_+)(T/2)}\det \tilde \Psi(T)\\
&\sim
C_0 C_- \det(V_1^+,\dots,V_n^+,V^-_{n+1},\dots,V_{2n}^-),
\endaligned
\eqnlbl{frontasymptotic}
$$
where 
$$
C_0(\lambda):=
e^{\int_0^\infty (\Tr \BbbA_\infty(\lambda)
-\Tr \BbbA(\lambda,x))dx}=\CalO(1).
$$
Comparing \eqnref{perasymptotic} and \eqnref{frontasymptotic}, we 
obtain the result, with
$$
C(k,\lambda):=
\big(C_0  \det(V_1^+,\dots,V_n^+,V^-_{n+1},\dots,V_{2n}^-)\big)^{-1}.
$$
\qed
\bigskip

{\bf Remarks.}  1.  As noted above, the calculation in the proof
of Proposition \thmref{limit} also shows explicitly that the $k$-loops
of the spectrum shrink as $T\to \infty$ to single eigenvalues, 
thus illuminating the passage from continuous to point spectrum 
as the period $T$ goes to infinity.
Though we did not do it here,
the rate of shrinking could be quantified
by a more careful version of the same argument.

2. The final assertion of Proposition \thmref{2}
(the one that mainly concerns us here)
is intuitively clear from consideration of
a small multiple of an unstable eigenfunction of the homoclinic
wave given as initial perturbation of a periodic wave--- 
this will grow exponentially for some time, by the property of
``almost finite propagation speed,'' or approximate localization of behavior.

3. In the more general case that the family of orbits
corresponding to $\{\bar u^\varepsilon\}$ converges as $\varepsilon\to0$ to
a heteroclinic cycle,  one can carry out an entirely similar analysis
to show that the $k$-spectra of the periodic operators $L^\varepsilon$, 
for $k$ fixed, $\varepsilon$ sufficiently small,
correspond approximately to the {\it union} of the spectra of the 
linearized operators
about each of the heteroclinic  raveling waves in the limiting cycle.
(This computation is in the spirit of ``multi-bump'' calculations
carried out for multiple traveling-pulse solutions in models
of nerve-impulse and optical transmission).
\medskip

{\bf Applications.}
In [GZ], there was derived a {\it stability index}
suitable for the evaluation of stability of traveling-front- or 
-pulse-type solutions of conservation laws. 
Precisely, this index yields the parity of the number of unstable
eigenvalues $\lambda$ such that $\R \lambda >0$ of
the linearized operator about the wave; evidently odd parity establishes
instability, while even parity is consistent with stability
but inconclusive.
This index was explicitly evaluated in [GZ] for the
class of planar, Hamiltonian quadratic-flux models, and later
in [Z.1] for the equations of van der Waals gas dynamics:
precisely the models we will study here in connection with pattern
formation.
In both cases, homoclinic, or pulse-type, solutions were found
to have odd index, implying instability, while heteroclinic, 
or front-type, solutions
were found to have even index, consistent with stability.
Though the latter result is inconclusive,
there is substantial evidence from other quarters [AMPZ.1,LZ.1--2]
that the heteroclinic waves are indeed stable.

Combining these results with the result of Proposition \thmref{2},
we find for these models that periodic solutions approaching a 
bounding homoclinic must be {\it unstable in the large-period limit},
an observation that was already made in [GZ].
This is suggestive, but not conclusive evidence regarding the nature of
numerically observed pattern formation.
For, though numerically observed patterns do appear to lie near limiting 
separatrices [AMPZ.5], the distance relative to the requirements
of the abstract theory is difficult to quantify.
Also, in both cases, the patterns appear to lie near not only homoclinic
but also heteroclinic cycles: $3$- and $2$-cycles, respectively.
Presuming, as evidence suggests, that the component heteroclinic fronts 
of such a limiting cycle are stable,  we find from Remark 3 above that 
the approaching periodic waves are at worst {\it weakly unstable}, in
the sense that unstable eigenvalues must have vanishingly small real part
as $\varepsilon \to 0$.  Thus, the question of stability is in this case
much more sensitive. 
Indeed, our numerical experiments in Appendix A
indicate that such periodic waves are in the large-amplitude
(i.e., large-period) limit
{\it stable} with respect to period-$T$/($k=0$) perturbations, but
{\it unstable} with respect to perturbations of 
some other periods/values of $k$, whereas, recall, the methods of this
section do not distinguish between different $k$-values.
Similar considerations hold in the small-amplitude limit; 
see Remark 4.5 below.

These difficulties motivate our development in the following sections
of a more direct, and general approach to stability of periodic waves.
First, rather than using the front or pulse stability index to obtain
information in the large-period limit, we will define an analogous
stability index in the periodic case, and compute this directly,
thus obtaining, at least in principle, information about
small- and intermediate- as well as large-amplitude waves.
This corresponds essentially to determining the Taylor expansion
in $\lambda$ of the Evans function at the origin $(k,\lambda)=(0,0)$.
Next, performing a similar
but less generally applicable
Taylor expansion in the variable
$k$, we will obtain a complementary, ``long-wave'' stability criterion
differentiating between different values of $k$.
Together, these will turn out to be sufficient to obtain 
rigorous instability results relevant to the pattern formation 
phenomena discussed above: in particular, in the case of the van der
Waals equations,  for waves of {\it arbitrary} amplitude.

{\bf Remark \thmlbl{multi}.}
At least formally,\footnote{
There is a technical difficulty associated with the accumulation of the
essential spectrum at the imaginary axis of the linearized operators
about these waves;  in particular, note that the argument of Proposition
\thmref{2} makes essential use, specifically
in the derivation of the crucial estimate \eqnref{perasymptotic},
of the fact that $\phi_1^+,\dots,\phi_n^+$ are exponentially decaying for
$x$ near $+\infty$, while $\phi_{n+1}^-,\dots,\phi_{2n}^-$ are exponentially
growing for $x$ near $-\infty$ (which is to say they decay as $x\to -\infty$).
It is a standard fact (see, e.g. [He,GZ,ZH])
that this property is equivalent to the assumption that $\lambda$
lie uniformly to the right of the essential spectrum boundary of $L$.
On the other hand, the more general bundle construction of [AGJ,G.1--2]
should still go through, in conjunction with analytic continuation into
the essential spectrum of the traveling-front- or -pulse-type Evans 
function using the Gap Lemma of [GZ,KS]; see, for example, the
related analysis in [DGK].
}
the same kind of ``multi-bump'' analysis as
sketched in Remark 3, above, suggests that, for a family of
homoclinic orbits $\{\bar u^\varepsilon\}$ converging to
a $2$-cycle limit $\bar u^0$, the spectra of the associated
linearized operators $L^\varepsilon$ should converge as $\varepsilon\to0$
to approximately the union of the spectra of the component heteroclinic 
waves in the $2$-cycle.
We have conjectured that these heteroclinic waves are stable, in which
case their point spectra consist, in the set $\{\lambda:\, \R\lambda\ge 0\}$,
of single translational eigenvalues at $\lambda=0$.  Pursuing this
line of reasoning, we deduce that the approaching homoclinics 
should have, besides the obligatory translational eigenvalue at 
$\lambda=0$, a second eigenvalue of vanishingly small real part,
the sign of which determines stability or instability.
Thus, the stability of homoclinics would be quite sensitive 
in the limit as they approach a $2$-cycle, a result consistent with the
sensitive stability we have deduced for their nearby periodic
waves (which must with them approach to the same $2$-cycle).
\medskip
{\bf Remark \thmlbl{smallamp}.}
It is interesting also to consider the opposite situation from that
studied by Gardner, namely the {\it small-amplitude limit} 
as a family of periodic waves $\{\bar u^\varepsilon\}$ 
shrinks to a single point, or nonlinear center $u_c$.
Specifically, let us consider the typical situation of a Hopf
bifurcation, for which the linearization 
$$
v'=B^{-1}A(u_c) v
\eqnlbl{linearizedrest}
$$
of the traveling-wave ordinary differential equation $Bu'=Au$ about $u_c$ 
possesses a two-dimensional
center manifold corresponding to a single pair of complex eigenvalues
$\alpha_\pm:=\pm i\tau$ of the coefficient matrix $B^{-1}A(u_c)$.
In this scenario, the period $T^\varepsilon$ of $\bar u^\varepsilon$
converges as $\varepsilon \to 0$ to $T^0:= 2\pi/\tau$, while the
amplitude shrinks to zero.

Likewise, the eigenvalue equation \eqnref{4} converges as $\varepsilon\to0$
to a constant-coefficient equation as considered in Example \thmref{cc},
with $A\equiv A(u_c)$ and $B\equiv B(u_c)$.
Referring to \eqnref{char},
we find that, for $\lambda=0$, there is an $n$-fold root $\mu=0$
and $n$ remaining roots $\mu$ consisting of the eigenvalues
of $B^{-1}A$.  Consulting \eqnref{zeroset}, and noting
that roots $\pm i\tau$ are by definition equal to $\pm 2\pi i/T^0$,
we find, fixing $k=0$, that $\lambda=0$ is an $(n+2)$-fold root of $D(0,\cdot)=0$,
while all other roots are far from the origin, and may in fact be stable:
this is the case, for example, when $B=I$.
Thus, stability with respect to periodic perturbations of period $T$,
corresponding to stability/instability of the roots of the
restricted Evans function $D(0,\cdot)$,
may be quite sensitive in the small-amplitude limit.
On the other hand, waves are usually quite {\it unstable}
with respect to perturbations on different periods, corresponding
to zeros of $D(k,\cdot)$ with $k\ne 0$: for example, in case $B=I$,
we have that $A$ has imaginary eigenvalues $\pm i\tau$, and so
hyperbolicity, hence also \eqnref{8}, is violated.
\medskip

\sectionnumber=5
\theoremnumber=0
\equationnumber=0
\smallskip
\TagsOnLeft
\specialhead
\S5.\quad 
The Stability Index
\endspecialhead

Motivated by the discussion of the previous section, we now
define directly a stability index for periodic waves analogous
to the one defined in [GZ] for traveling-front- or -pulse-type solutions.
This index, being designed to identify {\it strong} instabilities,
concerns only ($k=0$)-eigenvalues, i.e., instability with respect to
periodic perturbations of the same period $T$ as the background solution
$\bar u$.
For this purpose the choice of $k$ is expected to be somewhat arbitrary
(by continuity; see also the large-period analysis of the previous section).  
And, this allows us to make a link to the previous analyses of heteroclinic waves,
for which the Evans function has only the argument $\lambda$.

As in the analysis of [GZ] in the traveling-front or -pulse case, 
our goal is to relate $\sgn D(0,\lambda)$ for $\lambda$ near 
$+\infty$ to $\sgn D(0,\lambda )$
for $\lambda$ near $0$, for $\lambda $ restricted to the real axis,
and, in turn, to relate the latter to the dynamics of the traveling-wave 
ordinary differential equation.
This approach, introduced by Evans in the pioneering papers [E.1--4],
and widely generalized in, e.g., [J,AGJ,PW],
has proven to be a powerful tool in the stability analysis of
traveling-front or -pulse type solutions.
However, until now
it does not seem to have been carried out in the periodic case. 

\medskip

{\bf Large $\lambda $ Behavior}.  
We begin with the large-$\lambda$ limit, 
which admits a particularly simple treatment.

\proclaim{Lemma \thmlbl{4}}  As $\lambda \to +\infty$ along the
real axis, $\sgn
D(0,\lambda )\to (-1)^n$, where $n$ as in \text{\rm \eqnref{1}} is the dimension
of $u$.
\endproclaim

{\bf Proof}.  
By standard G\"arding-type (i.e., sectorial) energy estimates,
$L$ has no spectrum in $\R\lambda\ge 0$ for $|\lambda|$ sufficiently
large.  Moreover, the Evans function varies continuously with respect
to continuous changes in the coefficients of $L$, by continuous
dependence with respect to initial data of solutions of ordinary differential equation.
Thus, the quantity
$$
\lim_{\lambda\to +\infty} \sgn D(0,\lambda),
\eqnlbl{limitsign}
$$
with $\lambda$ restricted to the real axis, is both well-defined and
invariant under homotopy in $L$ within the class of strictly 
elliptic operators with periodic coefficients of period $T$.
Deforming $L$ to the Laplacian $\bar L:=(\partial/\partial_x)^2$ 
via the homotopy
$$
\theta \bar L+ (1-\theta) L,
$$
$\theta$ going from $0$ to $1$,
we may thus evaluate \eqnref{limitsign} by an explicit and elementary 
computation, which we omit.
This could alternatively be carried out directly, for the original operator
$L$, using a rescaling argument as in [GZ].
A related, but much more complicated,
homotopy argument was used in [BSZ] to treat the traveling-front or -pulse case.
\qed
\bigskip

{\bf Remark.}  Note that in this periodic context, the sign of
\eqnref{limitsign} is an absolute quantity, and does not depend
on any choice of coordinates.
Likewise, in the homoclinic case, there is a natural choice of
coordinates by which the sign may be made absolute 
(namely, choosing initializing bases of the stable subspace at $+\infty$ 
and the unstable subspace at $-\infty$
that together form a basis of $\BbbR^{2n}$ 
with the standard orientation; see construction of the front- or pulse-type
Evans function in the previous section).
By contrast, a significant difficulty confronted in [GZ] for
the traveling-front case was to relate
the sign at infinity to the normalization chosen at $\lambda=0$.
Indeed, this limited the original analysis of [GZ] to the case $n=2$;
for the extension to the general case, see [BSZ,Z.3].
\medskip

{\bf Small $\lambda $ Behavior}.  
We next address the crucial small-$\lambda$ case.
By analogy with the traveling-front or -pulse case,
we seek to relate small-$\lambda $ behavior to the
dynamics of the traveling-wave ordinary differential equation \eqnref{1.1}--\eqnref{1.2}.

Notice, in the present periodic context, that \eqnref{1.1}, \eqnref{1.2} 
involves $2n+1$ parameters $(u_-,u_0,s)$ rather than the $n+1$
parameters $(u_-,s)$ of the traveling-front or -pulse case,
since the choice of initial condition $u_0$ is completely 
independent of the critical point $u_-$.
For each choice of parameters, there is a unique solution
$\barU^{(u_{-},u_{0},s)}(x)$ of \eqnref{1.1}.
We can thus define the {\it special separation function}
$$
d(u_-,u_0,s):=
\barU^{(u_{-},u_{0},s)}(x) |^T_0.
\eqnlbl{29}
$$
{\it Note:}  This is a bit different from the usual separation
function in that the vanishing of $d(\cdot)$ corresponds to existence of
a periodic solution of {\it precisely period $T$}.  
A standard Melnikov function would be based, rather, on the 
Poincar\'e return map; here, however, we are concerned only with period $T$.

Variations
$$
w_j:= \partial\barU/\partial u_0 \cdot e_j, 
\quad j=1,\cdots n,
\eqnlbl{30}
$$
satisfy the linearized traveling-wave equation
$$
Bw'=Aw,
\quad w(0)=e_j.
\eqnlbl{31}
$$
Without loss of generality, take coordinates such that
$$
\barU_x(0)=e_1;
\eqnlbl{32}
$$
hence 
$$
w_1=\barU_x.
\eqnlbl{33}
$$

Likewise, the variation
$$
\tildeZ_1 := -\partial \barU/\partial s
\eqnlbl{34}
$$
satisfies
$$
Bz'=Az+(\barU-u_-), \quad z(0)=0,
\eqnlbl{35}
$$
while
$$
w_{n+j}:=(\partial \barU/\partial u_-) \cdot e_j, \quad 
j=1,\cdots n,
\eqnlbl{36}
$$
satisfy
$$
Bw'=Aw - f'(u_-)e_j, \quad w(0)=0.
\eqnlbl{37}
$$

We thus have relations 
$$
(\partial d/\partial u_0) \cdot e_j=[w_j],
\eqnlbl{38}
$$
$$
-\partial d/\partial s = [\tildeZ_1],
\eqnlbl{39}
$$
$$
(\partial d/\partial u_-) \cdot e_j = [w_{n+j}],
\eqnlbl{40}
$$
where $e_j$ denotes the $j$th standard basis element in $\BbbR^n$.

With these definitions, we may restate (H3) in the more quantitative form:
\medskip
(H3')
\qquad  $\partial d/\partial (u_0,s)=([\tildeZ_1],[w_1],\dots, [w_n])$ 
is full rank,
\medskip
\noindent i.e, (recalling that $[w_1]=[\barU_x]=0$)
$$
\gamma := \det([\tildeZ_1],[w_2],\dots, [w_n])\not= 0.
\eqnlbl{g1}
$$
%
That is, for a fixed period $T$ and equilibrium $u_-$,
the orbit $\barU(\cdot)$ is locally unique up to translations, even
allowing variation in $s$; moreover, it corresponds to a transverse
intersection of the tangent manifolds at $\barU(0)^{(u_-,u_0,s)}$
and $\barU(T)^{(u_-,u_0,s)}$ with respect to variations in $(u_0,s)$.  

Then, by the Implicit Function Theorem, there is an $(n+1)$-dimensional
surface 
$$
(u_-,u^1_0) \to (u_-,u_0,s)
\eqnlbl{42}
$$
in parameter space 
for which $\barU^{(u_-,u_0,s)}(\cdot)$ has period
$T$, where $u^1_0:=u_0\cdot e_1$ denotes the component of $u_0$
in the $e_1=\barU_x(0)$ direction.
Fixing $u^1_0$ to factor out translation invariance and fix the
phase, we find $n$ directions in which connections persist.  
We can thus uniquely specify $n$ variations in $\barU^{(u_-,u_0,s)}$:
$$
\tildeW_{n+j}\in w_{n+j}\oplus\Span (\tildeZ_1,w_2,\cdots w_n),
\eqnlbl{43}
$$
by the requirement
$$
[\tildeW_{n+j}]=0.
\eqnlbl{44}
$$

Finally, we define the ``zero-viscosity stability coefficient:''
$$
\Delta:=\det \big(\int \tildeW_{n+1},\cdots, \int \tildeW_{2n} \big)
= \det(\partial \bar u_m/\partial u_-),
\eqnlbl{g2}
$$
where $\int \tildeW_{n+j}$ denotes $\int^T_0 \tildeW_{n+j}(y) dy$
and
$$
\bar u_m(u_-,u_0^1):=\int_0^T \bar u^{(u_-,u_0,s)}(y)dy
\eqnlbl{mean}
$$
denotes the mass over one period of the solution 
$\bar u^{(u_-,u_0,s)}$ determined by $(u_-,u_0^1)$ via map \eqnref{42};
note that this does not depend on $u_0^1$, since it is invariant 
under phase shifts.
The condition $\Delta\ne 0$ corresponds 
to the requirement that there exist no nearby periodic orbits of period $T$
that have the same mass as does $\barU(\cdot)$,
at least up to linear order in the perturbation of $(u_-,u_0,s)$.
This has an interesting heuristic interpretation in the spirit of [FL.1--2]
as linearized well-posedness within a special class of 
measure-valued solutions, in the limit of zero viscosity,
of an associated {\it Riemann problem} having left and right states both 
equal to the average value of $\bar u$ over one period: namely,
the class of measure-valued
solutions for which the associated limiting sequence
consists of periodic functions with a fixed ratio between
period and viscosity.
(Note that average value is preserved both under compact perturbations 
of $\barU(\cdot)$ and, by conservation of mass, under the 
nonlinear flow of \eqnref{1}.)  
The latter restriction comes from the fact that we are here
considering stability only within the class of periodic functions
of fixed period $T$; its somewhat awkward form
reflects the link via rescaling between the long-time
and small-viscosity limits.
For a careful description of measure-valued
solutions and their relation to asymptotic behavior of \eqnref{1},
we refer the reader to [FL.1--2].

Alternatively, $\Delta\ne 0$ may be viewed as the requirement that,
near $\bar u$, there is a unique periodic orbit of period $T$ 
having a given mass $\bar u_m$ over one period: that is, the ``mass map''
from perturbation mass to possible time-asymptotic (periodic) states
is both well-defined and nonsingular.
Since mass per period is preserved under the flow of \eqnref{1} with periodic
boundary conditions, this condition is clearly necessary for orbital
stability of $\bar u$ within the class of $T$-periodic solutions,
under perturbations with nonzero mass.

In this sense, $\Delta$ is precisely analogous to the corresponding
one-dimensional zero-viscosity stability coefficient $\Delta(0,1)$
defined in [ZS,Z.3] 
\footnote{
The relation between the stability index, linearized 
well-posedness of the Riemann problem, and the mass map was first noted in [GZ],
in slightly less explicit form; see [FreZ] for related applications.
}
for the traveling-front or -pulse case, which
has the same relations to linearized well-posedness of the Riemann problem
and nonsingularity of the mass map.
With these definitions, we have the following fundamental relation,
analogous to the one described in [GZ,ZS,BSZ,Z.3] 
for the traveling-front or -pulse case.

\proclaim{Proposition \thmlbl{5}}  
Let \text{\rm (H0)--(H3)} hold. 
Then
$$
D(0,\lambda)=\lambda ^{n+1} (-1)^n\det \big(df(u_-)^{-1}\big) 
\gamma \Delta+\CalO(\lambda ^{n+2}),
\eqnlbl{expansion}
$$
or, equivalently,
$$
(\partial/\partial \lambda )^k D(0,0)=0,
\quad 0 \le k \le n,
\eqnlbl{44.1}
$$ 
and
$$
(\partial/\partial \lambda )^{n+1}D(0,0)=(-1)^n (n+1)! 
\det df(u_-)^{-1}\gamma  \Delta,
\eqnlbl{45}
$$
where $\gamma$ is the transversality coefficient defined in
\text{\rm \eqnref{g1}}
and $\Delta=\det(\partial \bar u_m/\partial u_-)$ 
is the inviscid stability coefficient defined in 
\text{\rm \eqnref{g2}, \eqnref{mean}}.
\endproclaim

{\bf Remark \thmlbl{fullrank}.}
In case $\gamma =0$, but $\partial d/\partial(u_0,u_-,s)$
is still of full rank, a similar calculation gives
$$
(\partial/\partial \lambda )^{n+1} D(0,0)=(-1)^n (n+1)! 
\det df(u_-)^{-1}\gamma_1  \Delta,
$$
where $\gamma _1$ now denotes the determinant of an appropriately chosen
alternative transverse set of vectors, and $\Delta$ is as in
\eqnref{g2}, with $\tildeW_{n+1},\dots, \tildeW_{2n}$ chosen again to
span the tangent manifold to the family of nearby periodic orbits,
modulo translation, i.e., satisfying \eqnref{44}. (Of course,
these cannot now be defined as in \eqnref{43}.)
Thus, we see precisely the same relation to linearized well-posedness of Riemann
problems as found in the case of Lax and undercompressive shock
waves in [GZ,ZS,BSZ,Z.3]. 
\medskip
{\bf Remark \thmlbl{small_lim}.}
It may happen that $\gamma \to 0$ but $\Delta$ remains bounded 
as $(u_,u_-,s)$ and $T$ approach certain limiting values.
In this case, we may conclude that $\Gamma=0$ for the limiting periodic orbit, 
since $D$, as the uniform limit of analytic functions,
has continuous partial derivatives as well.
For example, $\gamma=0$ in the small-amplitude, constant-coefficient
limit described in Remark \thmref{smallamp}, since in this case there 
is a one-parameter family of periodic orbits with period $T=T^0$,
but also $\partial \bar{u}/\partial u_- \to I$ as this limit is approached;
see Remark 5.10, below.
Thus, we may conclude that $\Gamma=0$ in the constant-coefficient case.
This is consistent with our previous observation, obtained
by direct calculation, that $D(0,\cdot)$
must vanish in this case to order $(n+2)$ and not $(n+1)$.
\medskip

{\bf Proof of Proposition \thmref{5}.}  
Defining the fundamental set of solutions 
$W_1,\dots,W_{2n}$
of eigenvalue ordinary differential equation \eqnref{4} by initialization
$$
(W_1,\dots,W_{2n})(0)= \pmatrix I&0\\
B^{-1}A & -B^{-1}df(u_-)_{|_0}
\endpmatrix,
\eqnlbl{init}
$$
and writing $W_j=:(w_j,w_j')^t$, we obtain by Remark \thmref{1.1}
the representation
$$
\aligned
D(0,\lambda )&= \det([W_1],\dots,[W_{2n}])/\det(W_1(0),\dots,W_{2n}(0))\\
&=(-1)^n \det 
\pmatrix
[w_1], \cdots [w_{2n}]\\
[w'_1], \cdots [w'_{2n}]
\endpmatrix/\det \big( B^{-1}(u_0)df(u_-)\big).
\endaligned \eqnlbl{52}
$$

With this choice of coordinates, 
the $w_j(\lambda)$ defined here agree at $\lambda=0$ with the 
traveling-wave ordinary differential equation variations 
$w_j$ defined in \eqnref{31}, \eqnref{36}, satisfying the (second-order)
linearized traveling-wave ordinary differential equation
$$
(Bw')'=(Aw)'.
\eqnlbl{2ndevalue}
$$
Likewise, variations $z_j:={w_j}_\lambda $ are seen at $\lambda=0$
to satisfy 
$$
(Bz')'=(Az)'+w_j, \quad z(0)=z'(0)=0,
\eqnlbl{53}
$$
and $y_j:= {w_j}_{\lambda \lambda}$ to satisfy 
$$
(By')'=(Ay)'+ 2 z_j, \quad y(0)=y'(0)=0.
\eqnlbl{53y}
$$
In particular, $z_1$ satisfies
$$
(Bz')'=(Az)'+\barU_x, \quad z(0)=0,
\eqnlbl{54}
$$
hence 
$$
Bz'=Az+(\barU-u_0), \quad  z(0)=0.
\eqnlbl{55}
$$
Comparing \eqnref{55} with \eqnref{35}, \eqnref{37}, we find that
$$
z_1=\tildeZ_1\quad \hbox{modulo span}\ (w_{n+1},\cdots w_{2n}).
\eqnlbl{56}
$$

Next, integrating \eqnref{2ndevalue} from $0$ to $T$, we find that,
at $\lambda=0$,
$$
B(u_0) [w'_j]=[Bw'_j]=[Aw_j]=A(u_0)[w_j];
\eqnlbl{57}
$$
hence
$$
[w'_j]-B^{-1}A(u_0)[w_j]=0
\eqnlbl{58}
$$
for all $1 \le j \le 2n$.  Likewise, integrating \eqnref{53},
\eqnref{53y} gives
$$
[z'_j]-B^{-1}A(u_0)[z_j]=B^{-1}(u_0)\int^T_0 w_j,
\eqnlbl{59}
$$
$$
[y'_j]-B^{-1}A(u_0)[y_j]=2 B^{-1}(u_0)\int^T_0 z_j,
\eqnlbl{59y}
$$
respectively.
In the special case $j=1$, \eqnref{59} gives
$$
[z'_1]-B^{-1}A(u_0)[z_1]=B^{-1}(u_0)\int^T_0 \barU_x=\barU |^T_0 =0.
\eqnlbl{60}
$$

Now, Taylor expanding formula \eqnref{52} about $\lambda=0$, 
we obtain for small $\lambda $ that 
$$
\aligned
D(0,\lambda )&= 
(-1)^n\det \big(df(u_-)^{-1}B(u_0)\big) \\
&\times \det
\pmatrix
[w_1]+\lambda [z_1]+ \frac{1}{2}\lambda^2[y_1]+ \cdots, &\dots, & [w_{2n}]+\lambda [z_{2n}]+\cdots\\
[w'_1]+\lambda [z'_1]+\frac{1}{2}\lambda^2[y_1']+\cdots, &\dots, &[w'_{2n}]+\lambda [z'_{2n}]+\cdots
\endpmatrix.
\endaligned
\eqnlbl{61}
$$
Subtracting $B^{-1}A(u_0)$ times the first row from the second, 
and using \eqnref{59}, \eqnref{59y} and $\int w_1=[w_1]=[w'_1]=0$,
we obtain
$$
\aligned
D(0,\lambda )&= 
(-1)^n\det \big(df(u_-)^{-1}B(u_0)\big) \\
&\times \det
\pmatrix
\lambda [z_1]+\cdots, &[w_2]+\cdots, &\dots, &[w_{2n}]+\cdots\\ 
\lambda^2 B^{-1}\int z_1 +\cdots,&\lambda B^{-1}\int w_2 +\cdots, & \dots, 
&\lambda B^{-1}\int w_{2n}+\cdots
\endpmatrix\\
&=\lambda ^{n+1}
(-1)^n\det \big(df(u_-)^{-1}\big) 
\det \pmatrix
[z_1], & [w_2], &\dots, &[w_{2n}]\\
\int z_1, &\int w_2, &\dots, &\int w_{2n}
\endpmatrix \\
&\quad +\CalO(\lambda ^{n+2}),
\endaligned \eqnlbl{62}
$$
where $\int$ denotes $\int^T_0$.  Applying now appropriate column
operations, we obtain from \eqnref{43}, \eqnref{44}, \eqnref{56} that
$$
\aligned
D(0,\lambda ) &= \lambda ^{n+1}
(-1)^n\det \big(df(u_-)^{-1}\big) \\
&\times \det
\pmatrix
[\tildeZ_1], & [w_2],\dots, [w_n], &0, &\dots, &0\\
\int \tildeZ_1,&\int w_2,\dots, \int w_n, &\int \tildeW_{n+1},&\dots
&\int \tildeW_{2n}
\endpmatrix\\
&\qquad +\CalO(\lambda ^{n+2})\\
&=\lambda ^{n+1} (-1)^n\det \big(df(u_-)^{-1}\big) \\
&\times
\det([\tildeZ_1],[w_2],\dots,[w_n])\det
\big(\int\tildeW_{n+1},\dots, \int \tildeW_{2n}\big)\\
&\qquad +\CalO(\lambda ^{n+2})\\
&=\lambda ^{n+1} (-1)^n\det \big(df(u_-)^{-1}\big) 
\gamma \Delta+\CalO(\lambda ^{n+2}),
\endaligned \eqnlbl{63}
$$
as claimed.\qed
\bigskip


{\bf Remark \thmlbl{6.1}}  The solutions $w_{n+1},\dots,w_{2n}$,
satisfy the inhomogeneous equations
$$
Bw'_{n+j}=Aw_{n+j}- df(u_-)e_j, \quad w_{n+j}(0)=0.
\eqnlbl{38}
$$
Thus, they also can be expressed in terms of the fundamental solution
$\psi(x)$ of the $n$-dimensional linearized traveling-wave equation
$$
Bw'=Aw
\eqnlbl{27}
$$
(different from the fundamental solution $\Psi$ of the
$2n$-dimensional eigenvalue ordinary differential equation \eqnref{4})
via Duhamel's principle:
$$
w_{n+j}(x) =
-\Big(\int^x_0 \psi(x) \psi(y)^{-1}dy\Big) \, df(u_-)e_j,
\eqnlbl{39}
$$
and therefore
$$
\aligned
[w_{n+j}] 
&=w_{n+j}(T)\\
&=\Big(\psi(T)\int^T_0 \psi(y)^{-1}dy\Big) df(u_-)e_j.
\endaligned \eqnlbl{40}
$$
Likewise, 
$$
[w_j]= (\psi(T)-\psi(0))e_j,
\eqnlbl{41}
$$
$$
\int w_j = \Big(\int^T_0 \psi(y)dy \Big)  e_j,
\eqnlbl{42}
$$
$$
[\tildeZ_1]=\psi(T)\int^T_0 \psi(y)^{-1}(\barU(y)-u_-)dy.
\eqnlbl{43}
$$

That is, formula \eqnref{g2} for $\Delta$ can be interpreted as a sort
of Melnikov integral, like the corresponding object
in the homoclinic (undercompressive) case of the traveling-front or 
-pulse theory [GZ,Z.1].
\medskip

{\bf Remark \thmlbl{6.2}} The quantity $[z_1]$ can clearly 
be substituted for $[\tildeZ_1]$ in \eqnref{g1}, since the two quantities 
are equal, modulo $\text{\rm Span} (w_{n+1},\cdots w_{2n})$.
\medskip

{\bf The stability index.}  With these preparations, we
now define the {\it stability index}:
$$
\Gamma := \sgn(\partial/\partial
\lambda )^{n+1}D(0,0)D(0,+\infty).
\eqnlbl{Gammadef}
$$
We have immediately:

\proclaim{Proposition \thmlbl{obvious}}
Let \text{\rm (H0)--(H3)} hold.
Then, the parity of
the number of unstable eigenvalues of $L$, considered as an operator
on the space of functions periodic on $[0,T]$ is
even if $\Gamma$ is positive, and odd if $\Gamma$ is negative.
In particular, $\Gamma \ge 0$ is necessary for weak spectral stability
as defined in \text{\rm \eqnref{spectralstability}}.
\endproclaim

{\bf Proof.}
Evidently, $D(0,\lambda)$ is invariant with respect to complex
conjugation, i.e., $D(0,\bar \lambda)= \bar D(0,\lambda)$,
by the definition of $D(\cdot,\cdot)$.  This yields the familiar
fact that the eigenvalues of the real-valued operator $L$, considered
as acting on the periodic functions on $[0,T]$, are either real
or else belong to complex conjugate pairs, whence the parity of
the number of unstable eigenvalues is equal to the parity of the number
of unstable {\it real} eigenvalues.  But, this is clearly determined
by $\Gamma$ in the manner stated.
\qed
\bigskip

{\bf Remark \thmlbl{strongnec}.}
The strict inequality $\Gamma >0$
is necessary for {\it strong spectral stability}
as defined in [OZ].  In particular, it is necessary for
$L^1\to L^p$ linearized stability for any $p<\infty$, by
Proposition 1.5 and condition ($\tilde{\text{\rm D3}}$)(i) of [OZ]
combined with
\eqnref{Gammadef}, \eqnref{37} above.
This is analogous to the situation in the traveling-front or -pulse
case; see discussion, Section 11 [ZH], of the ``neutrally stable case.''
\medskip

Combining the results of Lemma \thmref{4}, Proposition \thmref{5},
and Proposition \thmref{obvious},
we obtain the main result of this section, an expression for
$\Gamma$ involving only geometry of the phase space of the 
traveling-wave ordinary differential equation:

\proclaim{Theorem \thmlbl{6}}  
Let \text{\rm(H0)--(H3)} hold.  Then, 
$$
\Gamma = \sgn \gamma  \Delta \det df(u_-),
\eqnlbl{37}
$$
where $\gamma$ and $\Delta=\det(\partial \bar u_m/\partial u_-)$ 
are as defined in \text{\rm \eqnref{g1}} and 
\text{\rm \eqnref{g2}, \eqnref{mean}}, respectively.
In particular,
$$
\sgn \gamma  \Delta \det df(u_-)\ge 0
\eqnlbl{stab1}
$$
is necessary for weak spectral stability of $\bar u$.
\endproclaim

{\bf Remark \thmlbl{sufficient}.}
In the small-amplitude limit discussed in Remark \thmref{smallamp},
recall that $D(0,\cdot)$ typically has $n+2$ zeroes lying near
the origin $\lambda=0$, and the rest lying in the strictly stable
complex half-plane $\R \lambda <0$.  By Propositions \thmref{5}
and \thmref{obvious},
$n+1$ of these in fact lie precisely at $\lambda=0$, with
the remaining small eigenvalue being stable, unstable, or zero
according as the stability index $\Gamma$ is positive, negative,
or zero.  Therefore, in this case, $\Gamma\ge 0$ is {\it necessary
and sufficient} for weak spectral stability with respect to periodic
perturbations of period $T$.
\medskip

{\bf Remark \thmlbl{Id}.}
With appropriate choice of the rest point $u_-$,
it can be shown that the key quantity 
$(\partial \bar u_m/\partial u_-)$ in \eqnref{37}
approaches the identity in both the small-amplitude and the 
large-amplitude homoclinic limits considered earlier.
For the small-amplitude limit described in Remark \thmref{smallamp}, 
we may choose $u_-$ to be the rest point lying near
the nonlinear center $u_c$ to which the family of periodic orbits converges;
for the homoclinic limit, we choose $u_-$ to be the rest point lying
near the vertex of the homoclinic.
We omit the standard, but somewhat involved calculations giving these results.

This has an interesting implication in the planar case, $n=2$.
For, in the simplest case that $\gamma$ is nonvanishing
for the entire family of periodic waves, we find then by
Theorem \thmref{6} that the sign of $\Gamma$ depends entirely
on $\sgn \det(df(u_-))$.  The latter is negative when
$u_-$ is set equal to a homoclinic vertex, since this is a saddle 
of the traveling-wave ordinary differential equation, so that $\det(B^{-1}df(u_-))$, and thus also
$\det(df(u_-))=\det B \det(B^{-1}df(u_-))$, is negative,
as the product of a positive and a negative eigenvalue.
But, by similar reasoning, it is positive when $u_-$ is set 
equal to a nonlinear center, since then the eigenvalues of 
$\det(B^{-1}df(u_-))$ are complex conjugates.
Thus, combining this observation with that of Remark \thmref{sufficient}
just above, we find that if large-amplitude waves are spectrally unstable with 
respect to period-$T$ perturbations, then
small-amplitude waves are necessarily {\it stable} in this sense.
Regarding period-$T$ perturbations, we thus have the global picture
of a single, small real eigenvalue detaching from the origin at
amplitude zero and initially moving into the stable complex half-plane,
later crossing zero again at some critical intermediate amplitude
to become unstable thereafter;
this picture is borne out by our numerical investigations in Appendix A.
On the other hand, we see that the stability index by itself is 
not likely to be effective across the range of all amplitudes,
even when it is useful in the large-amplitude limit, and for 
small amplitudes is usually {\it not} useful.
\medskip

{\bf Remark \thmlbl{gamma}.}
By calculations similar to those 
in the proof of Assertion \thmref{lim} (see proof of
Proposition \thmref{limit}),
the expression for $\gamma$ given in \eqnref{g1}
may, in the large-amplitude homoclinic limit discussed in Section 4, 
be related to a corresponding quantity arising in the
study of stability of homoclinic, or pulse-type, traveling waves.
For $n=2$, this is the quantity $\Gamma$ defined 
in (3.18), Lemma 3.4 of [GZ]; a similar quantity arises in the general
case [BSZ,Z.3].
Combining this observation with that of Remark \thmref{Id} above,
one can show that the sign of the periodic stability index converges
in the large-amplitude homoclinic limit, to the sign of the
traveling-pulse-type index for the limiting homoclinic wave.
We omit the associated calculations as lying too far
from the direction of our main interest; however, see the 
related, explicit computation (6.26) in Section 6 below.
\medskip

\sectionnumber=6
\theoremnumber=0
\equationnumber=0
\smallskip
\TagsOnLeft
\specialhead
\S6.\quad 
The quasi-Hamiltonian case
\endspecialhead

We next restrict attention to the case that the traveling-wave 
ordinary differential equation \eqnref{1.1} is ``quasi-Hamiltonian''
in the sense that it admits an integral of motion
for each $(u_-,u_0)$, with $s$ held fixed at the base value
under consideration, without loss of generality $s=0$.
This situation arises in some physical applications; see
Section 7.  Furthermore, it has a particular mathematical
interest; as shown in [OZ], this is essentially the
only circumstance under which one can expect standard ``diffusive,''
or asymptotic $L^1\to L^p$ linearized stability of periodic
solutions:  in the generic case, linearized $L^1\to L^p$ stability 
can hold only for $p=\infty$, and is bounded at best; see [OZ],
Proposition 1.5 and discussion surrounding condition ($\tilde{\text{\rm D3}}$).

With this additional structure, we may determine not only the
low-frequency behavior of the Evans function  with respect to 
$\lambda$, but the entire Taylor expansion behavior with respect 
to $(k,\lambda)$:

\proclaim{Theorem \thmlbl{lowfreq}}  
Let \text{\rm(H0)--(H3)} hold, 
and suppose in addition that the traveling-wave 
ordinary differential equation \text{\rm\eqnref{1.1}}
is quasi-Hamiltonian for the value of $s$
associated with solution $\bar u$.
Then, for $k$, $\lambda$ sufficiently small,
$$
\aligned
D(k,\lambda)&= \gamma \lambda 
\det \big( -\lambda (\partial \bar u_m/\partial u_-)df(u_-)^{-1}
- ikT  \big) \\
&+ \CalO\big((|k|+|\lambda|)^{n+2}\big),
\endaligned
\eqnlbl{lf}
$$
where $\gamma$ and $(\partial \bar u_m/\partial u_-)$ 
are as defined in \text{\rm\eqnref{g1}} and 
\text{\rm \eqnref{g2}--\eqnref{mean}}, respectively.
\endproclaim

{\bf Proof.}
By assumption, \eqnref{1.1} admits an integral of motion $H(u_-,u)$
for each $(u_-,u_0)$, i.e.,
$$
H(u_-,u(t))\equiv H(u_-,u_0)
\eqnlbl{H}
$$
for all solutions $u(t)$, where $s$ is held fixed.
Then, perturbing the periodic orbit $\bar u$ with
respect to these parameters, we find that variations
$$
[w_j]:=(\partial/\partial u_0)\big(\bar u^{(u_-,u_0,s)}(T)- 
\bar u^{(u_-,u_0,s)}(0)\big),
\eqnlbl{orth0}
$$
$$
[w_{n+j}]:=(\partial/\partial u_-)\big(\bar u^{(u_-,u_0,s)}(T)- 
\bar u^{(u_-,u_0,s)}(0)\big)
\eqnlbl{orth-}
$$
must have range lying tangent to the hypersurface $ H(u_-,u)\equiv H(u_-,u_0) $
at the base values of $(u_-,u_0)$ associated with $\bar u$, 
i.e., orthogonal to $\nabla_u H(u_-,u_0)$.
For, differentiating $H(u_-,\bar u(T))-H(u_-,\bar u(0))=0$ with
respect to $u_0$ yields 
$$
0=\nabla_{u}H(u_-,u_0) [\partial \bar u/\partial u_0],
\eqnlbl{ort0}
$$
where, as usual, $[f]$ denotes $f(T)-f(0)$; likewise,
differentiating with respect to $u_-$ yields
$$
\aligned
0&=\nabla_{u}H(u_-,u_0) [\partial \bar u/\partial u_-] +
\nabla_{u_-}H(u_-,u_0)I- \nabla_{u_-}H(u_-,u_0)I\\
&=
\nabla_{u}H(u_-,u_0) [\partial \bar u/\partial u_-].
\endaligned
\eqnlbl{ort-}
$$
Moreover, $w_1:=\bar u_x(0)$ is also orthogonal 
to $\nabla_u H(u_-,u_0)$, since 
$$
0=(d/dx)H\big(u_-,\bar u(0)\big)=\nabla_u H(u_-,u_0)\bar u_x(0).
\eqnlbl{orth1}
$$

Noting that $[w_2],\dots,[w_n]$ must be independent, by (H3'),
we thus find that they are by themselves a basis for 
$$
\text{\rm Span } \{[w_j]\} \oplus 
\text{\rm Span } (w_1),
\eqnlbl{dep}
$$
with $[\tildeZ_1]$ not required.
This simply reflects the fact that, by the Implicit Function Theorem,  
the full $(n+1)$-parameter family of nearby periodic orbits of period $T$ 
can be generated with $s$ held fixed;
for, thanks to \eqnref{H} only $n-1$ 
conditions must be satisfied to obtain existence.

Defining $W_1,\dots,W_{2n}$ again as in \eqnref{init},
but allowing $k$ to vary, we obtain in place of \eqnref{52}
the representation:
$$
\aligned
D(k,\lambda )&= 
\det \big( -df(u_-)^{-1}B(u_0)\big)\\
&\times
\det\big([W_1]-(e^{ikT}-I)W_1(0),\dots, [W_{2n}]-(e^{ikT}-I)W_{2n}(0)\big).
\endaligned
\eqnlbl{52k}
$$
Taylor expanding this about $(k,\lambda)=(0,0)$, 
we thus obtain for small $k$, $\lambda $ that
$$
\aligned
&D(0,\lambda )= 
\det \big( -df(u_-)^{-1}B(u_0)\big)\\
&\times
\det
\pmatrix
[w_1]+\lambda [z_1] -ikTw_1(0)+\cdots, &\dots, 
& [w_{2n}]+\lambda [z_{2n}]- ikTw_{2n}(0)+\cdots\\
[w'_1]+\lambda [z'_1]-ikTw_1'(0)+\cdots, &\dots, 
&[w'_{2n}]+\lambda [z'_{2n}]-ikTw_{2n}'(0)+\cdots
\endpmatrix.
\endaligned
\eqnlbl{61}
$$

Recalling \eqnref{31}, \eqnref{36}, we find that
$$
\big(w_j'-B^{-1}A w_j \big)(0)= 0,
$$
$$
\big(w_{n+j}'-B^{-1}A w_{n+j} \big)(0)=  -B^{-1}(u_0)df(u_-)e_j
$$
for $j=1,\dots, n$.
Thus, subtracting $B^{-1}A(u_0)$ times the first row from the second, 
we obtain as in \eqnref{62} the expression
$$
\aligned
&D(k,\lambda ) =
\det \big( -df(u_-)^{-1}\big)\\
&\times \det
\pmatrix
\lambda [z_1]-ikTw_1(0), &[w_2], &\dots, &[w_{2n}]\\ 
\lambda^2\int \tildeZ_1,&\lambda \int w_2 , & \dots, &\lambda \int w_{2n}
+ikTdf(u_-)e_n
\endpmatrix \\
&\qquad +\CalO\big( (|k|+|\lambda|)^{n+2}\big)\\
&= \det \big( -df(u_-)^{-1}\big)
\det
\pmatrix
\lambda [z_1], &[w_2], &\dots, &[w_{2n}]\\ 
\lambda^2\int \tildeZ_1,&\lambda \int w_2 , & \dots, &\lambda \int w_{2n}
+ikTdf(u_-)e_n
\endpmatrix \\
&\qquad +\CalO\big( (|k|+|\lambda|)^{n+2}\big),\\
\endaligned \eqnlbl{62k}
$$
where in the final equality we have used \eqnref{dep} to
eliminate term $ikTw_1(0)$.
Applying column operations as in \eqnref{63}, we thus obtain
$$
\aligned
D(k,\lambda ) &=
\det \big( -df(u_-)^{-1}\big)
\det \pmatrix
\lambda [\tildeZ_1], & [w_2],\dots, [w_n], &0 \\
\int \tildeZ_1,&\int w_2,\dots, \int w_n, & M(k,\lambda)\endpmatrix\\
&\qquad +\CalO\big( (|k|+|\lambda|)^{n+2}\big),\\
\endaligned \eqnlbl{63k}
$$
where the $n\times n$ block $M$ is given by
$$
\aligned
M(k,\lambda)&:=
\big(\int \tildeW_{n+1}
+ikTdf(u_-)e_{1},\dots,
\int \tildeW_{2n} +ikTdf(u_-)e_{n}\big)\\
&=(\partial \bar u_m/\partial u_-)+ ikTdf(u_-)I.
\endaligned
\eqnlbl{M}
$$
Finally, factoring the determinant on the righthand side of \eqnref{63k} 
using its block triangular form yields
$$
\aligned
D(k,\lambda ) 
&=
\gamma\lambda
\det \Big( \lambda (\partial \bar u_m/\partial u_-)
+ikTdf(u_-) \Big)
\det \Big( -df(u_-)^{-1}\Big)\\
&\qquad +\CalO\big( (|k|+|\lambda|)^{n+2}\big).\\
\endaligned \eqnlbl{63final}
$$
Rearranging \eqnref{63final}, we obtain \eqnref{lf}.\qed
\bigskip

Notice that result \eqnref{lf} is consistent
with our previous result \eqnref{5}.
As an immediate consequence, we obtain:

\proclaim{Corollary \thmlbl{exp}}  
Let \text{\rm(H0)--(H3)} hold, 
and suppose in addition that the traveling-wave 
ordinary differential \text{\rm \eqnref{1.1}}
is quasi-Hamiltonian for the value of $s$
associated with solution $\bar u$.
If $\Delta:=\det (\partial \bar u_m/\partial u_-)\ne 0$,
then, for $k$, $\lambda$ sufficiently
small, the spectrum of $L$ consists of $(n+1)$ smooth curves:
$$
\lambda_0(k) = o(k),
\eqnlbl{lambda0}
$$
$$
\lambda_j(k)= -i\alpha_j kT + o(k), \quad j=1,\dots, n,
\eqnlbl{lambdaj}
$$
where $\alpha_j$ denote the eigenvalues of 
$$
df(u_-)(\partial \bar u_m/\partial u_-)^{-1},
\eqnlbl{alphamatrix}
$$
and where $(\partial \bar u_m/\partial u_-)$ 
is as defined in \text{\rm \eqnref{g2}, \eqnref{mean}}.
In particular,  the ``effective hyperbolicity'' condition
$$
\sigma\big( df(u_-)(\partial \bar u_m/\partial u_-)^{-1} \big)
\, \text{\rm real}
\eqnlbl{stab2}
$$
is necessary for weak spectral stability of $\bar u$.

If $\Delta=0$, then \text{\rm \eqnref{lambda0}} holds as before.
Likewise, \text{\rm\eqnref{lambdaj}}
holds with $\alpha_j:=\beta_j^{-1}$, 
where $\beta_j$ denote the eigenvalues of 
$(\partial \bar u_m/\partial u_-) df(u_-)^{-1}$,
for all $j$ such that $\beta_j\ne 0$.
Associated with vanishing $\beta_j$, however, are nonsmooth curves
$$
\lambda(k)\sim k^{p/q},
$$
$p/q$ not an integer; in particular, 
$ \det(\partial \bar u_m/\partial u_-)=0$ implies
spectral instability.
\endproclaim

{\bf Remark.}
Note that $\Gamma=0$ in this (quasi-Hamiltonian) case implies spectral, 
i.e., {\it exponential} instability, 
since then $ \det(\partial \bar u_m/\partial u_-)=0$
by \eqnref{37} together with (H2) and (H3').
This is to be compared with the (apparently) more subtle situation of the
general case; see Remark \thmlbl{strongnec}. 
\medskip

The formulae of Corollary \thmref{exp} may be regarded as the
natural generalizations of the constant-coefficient
formulae of Remarks \thmref{hyp} and \thmref{smallamp}.
In particular, comparing \eqnref{lambdaj}, \eqnref{alphamatrix} to
\eqnref{ccexp}, we see that 
$df(u_-)(\partial \bar u_m/\partial u_-)^{-1}$
plays in the quasi-Hamiltonian case the role played in the 
constant-coefficient case $\bar u\equiv u_0=u_-$ by $A:=df(u_-)$.
That is, $df(u_-)(\partial \bar u_m/\partial u_-)^{-1}$
may be regarded as an ``averaged'' or ``effective'' convection
matrix for the variable-coefficient case (indeed, this
is shown in [OZ] to hold true in a very strong sense).
Note, further,  that $df(u_-)(\partial \bar u_m/\partial u_-)^{-1}$
converges in the small-amplitude limit to 
$df(u_-)$, since, as pointed out in Remark \thmref{Id},
$(\partial\bar u_m/\partial u_-)\to I$.

Likewise, in the large-amplitude, homoclinic limit considered
by Gardner, we find again (Remark \thmref{Id}) that 
$df(u_-)(\partial \bar u_m/\partial u_-)^{-1}$
converges to $df(u_-)$, where $u_-$ now denotes the vertex
of the homoclinic orbit, i.e., the limiting state of the
homoclinic wave as $x\to \pm \infty$.
Thus, the spectral curves described in
\eqnref{lambdaj}, \eqnref{alphamatrix} correspond to first order in $k$
with the spectral curves of this limiting, constant state, which are in
turn curves of essential spectrum for the linearized
operator about the homoclinic wave: for a detailed discussion,
see the introduction of [OZ].
The remaining curve \eqnref{lambda0} corresponds to
the translational (continuous) eigenfunction at $\lambda=0$ 
associated with the spatial derivative $\bar u^\varepsilon$ of the 
profile under consideration.  

This extends to the small-frequency regime
the picture described by Gardner (Section 4) in 
the large-amplitude homoclinic limit.
It can be shown by an argument similar to that of Section 4
(together with the observations made in the introduction of 
[OZ] regarding structure of the resolvent set of the 
limiting homoclinic profile)
that the former curves in fact {\it globally} approach the essential
spectrum curves of the limiting homoclinic wave, in the sense that
they approach on a ball of radius going to infinity with the period;
thus, they extend arbitrarily far as the homoclinic limit is approached.
We conjecture without proof that the latter curve belongs to a 
closed loop of spectra, shrinking in the large amplitude limit to the 
eigenvalue $\lambda=0$ associated with the translational eigenfunction 
$\bar u^0_x$ of the limiting, homoclinic wave.

{\bf Remark.} The Taylor expansion of $D $ can likewise be carried
out in the general (non-quasi-Hamiltonian) case, but does not seem
amenable to any such simple interpretation.
\medskip

{\bf The planar Hamiltonian case.}
Finally, we consider the planar case, $n=2$, under the assumption
that the traveling-wave ordinary differential equation
 \eqnref{1.1} is in fact Hamiltonian 
for all $u_-$, with $s$ held fixed at the base speed associated 
with $\bar u$, i.e.,
$$
0=\text{\rm Tr }d_u \big( B^{-1}(u)(f(u)-f(u_-) - s(u-u_-))\big ).
\eqnlbl{ham}
$$
This is the class from which we will take our example systems in the 
following section.  Note that, for $B\equiv \text{\rm constant}$, 
condition \eqnref{ham} reduces to
$$
0=\text{\rm Tr } \big( B^{-1}(df(u)-sI)\big),
\eqnlbl{ccham}
$$
which is evidently independent of $u_-$; with this observation,
it is straightforward to construct examples.

{{From \eqnref{ham}, we find that the traveling-wave equation may be expressed
as a Hamiltonian system
$$
u'=  B^{-1}(u)(f(u)-f(u_-) - s(u-u_-))
=\nabla_u^\perp H(u,u_-),
\eqnlbl{hamiltonian}
$$
for an appropriate Hamiltonian $H(u,u_-)$, whence
$H(u(t),u_-)$ is an integral of motion; in particular, 
we see that the planar Hamiltonian systems are a subclass
of the quasi-Hamiltonian systems defined above.

Existence of an integral of motion 
in the planar case implies considerable structure
of the phase portrait of the traveling-wave ordinary differential equation.
For example, in the situation we consider, that the a phase portrait
contain a heteroclinic poly-cycle, we find that the interior of the
cycle must be entirely filled with cycles, either periodic or heteroclinic.
In the simplest case that there is only a single nonlinear center $u_c$
enclosed, the interior must be made up entirely of periodic orbits;
see, for example, the depiction of the homoclinic case in figure 1.
In this situation, we may globally parametrize these orbits by amplitude $a$, 
defined as distance from center $u_c$ along some (fixed) curve of 
steepest descent if $u_c$ is a local maximum, or ascent
if $u_c$ is a minimum.  In more general situations, we may still
define a {\it local} parametrization in this way, in the vicinity
of any fixed periodic orbit $\bar u$.  

\midinsert
\captionwidth{5truein}
\hskip 1.5 truein
\vbox to 2.5truein{\hsize=3truein \efig newlevelset y2.2
\botcaption{Figure 1} 
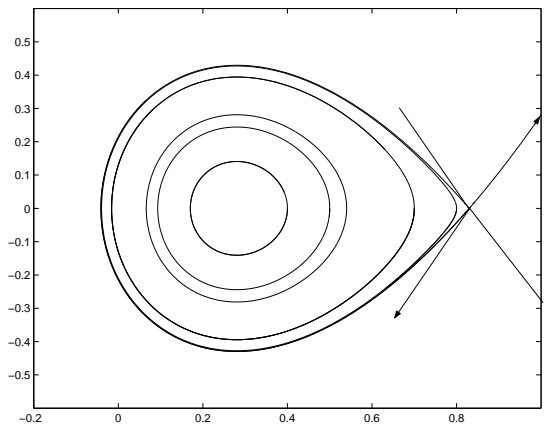
The planar Hamiltonian case: bounding homoclinic cycle.
\endcaption
}
\endinsert

%
%

Let us consider the period $T$ as a function of amplitude, $T=T(a)$.
Then, we have:

\proclaim{Lemma \thmlbl{sign}}
Let \text{\rm (H0)--(H3)}
hold, and suppose in addition that the traveling-wave
ordinary differential equation \text{\rm \eqnref{1.1}}
 is planar Hamiltonian for the value of $s$
associated with solution $\bar u$, with $B\equiv \text{\rm constant}$.
Then, 
$$
\sgn \gamma=\sgn (dT/da)
\eqnlbl{dTda}
$$
for $\gamma:=\det([\tildeZ_1],[w_2])$ as in \text{\rm \eqnref{g1}}.
\endproclaim

{\bf Proof.}
Let $w_1$, $w_2$ denote the variations defined in \eqnref{31}, \eqnref{36}, 
as usual normalizing $w_1=\bar u_x(0)=e_1$.  We may locally parametrize
amplitude by location of $u_0$ along a curve tangent to the outward 
normal $\eta=\pm e_2$ at $\bar u(0)$ to the orbit $\{\barU (x)\}$,
so that
$$
da/du_0= \eta
$$
at $\bar u(0)$, where $a(u_0)$ is the amplitude of the orbit corresponding
through $u_0$, with $(u_-,s)$ held fixed.
Differentiating the identity
$$
\bar u^{(u_-,u_0,s)}(T(a(u_0)))-\bar u^{(u_-,u_0,s)}(0)\equiv 0
$$
with respect to $u_0$ in direction $e_2$, we thus obtain
$$
[w_2] + \bar u_x(T)(dT/da)(\eta \cdot e_2)=0,
$$
or
$$
[w_2]= C\barU_x(0) = C\barU_x(T)=Ce_1,
\eqnlbl{sec2}
$$
with
$$
\sgn(C)=\sgn (dT/da)\sgn \det (\eta,\barU_x(0)).
\eqnlbl{3}
$$
The quantity 
$$
\gamma =\det([\tildeZ_1],[w_2])=C\det(\tildeZ_1(T),\barU_x(T))
\eqnlbl{3}
$$
can therefore be conveniently calculated, using inhomogeneous Abel's
equation/Duhamel principal as a Melnikov integral:
$$
\aligned
\gamma 
&=C\int^T_0
e^{\int\Tr(B^{-1}A)}\det(B^{-1}(\barU-\barU_-),\barU_x)dx\\
&= C\int^T_0 \det(B^{-1}(\barU-u_-),\barU_x)dx,
\endaligned \eqnlbl{sec4}
$$
where in the second equality we have used the Hamiltonian property
$\text{\rm Tr }(B^{-1}A)\equiv \text{\rm constant}.$
More precisely, this is found by setting 
$\gamma(t):= C\det(\tildeZ_1(t),\barU_x(t))$
and observing that $\gamma$ satisfies
$$
\gamma'= \text{\rm Tr }(B^{-1}A )\gamma + 
C\det(B^{-1}(\barU-u_-),\barU_x),
\quad \gamma(0)=0
$$
by \eqnref{31} and \eqnref{34}, \eqnref{35}; for similar calculations
in the traveling-front or -pulse case, see [GZ].

But, \eqnref{sec4}, as in [GZ], may be viewed as a contour integral
and explicitly evaluated using the Gauss-Green formula as 
$$
\aligned
\gamma &=
C\sgn \det(\eta,\barU_x)
\int_{\partial \Omega} B^{-1}(u-u_-)\cdot \eta \, ds\\
&=
C\sgn \det(\eta,\barU_x)
\int_{\Omega} \div \big( B^{-1}(u-u_-) \big) \, du\\
&=
C\sgn \det(\eta,\barU_x)
|\Omega|\text{\rm Tr }(B^{-1}),
\endaligned
\eqnlbl{5}
$$
where $\Omega$ is the region bounded by the orbit
$\barU(\cdot)$, and $\eta$ denotes the outward unit normal.
Recalling that $\text{\rm Tr }B^{-1}>0$ by (H1), and
combining with \eqnref{3}, we are done.
\qed
\bigskip

{\bf Remarks.}  1. Formulae \eqnref{sec4}, \eqnref{5} clearly show
that 
$$
\det([z_1],[w_2])=\det([\tildeZ_1],w_2]),
$$
as asserted in Remark \thmref{6.2},
since substitution of $u_0$ for $u_-$ leaves the result unchanged.

\smallskip
2.  Comparison of formula \eqnref{sec4} to that given for $\Gamma$
in (3.18), Lemma 3.4 of [GZ] shows explicitly the convergence 
described in Remark \thmref{gamma}, since the two formulae are
formally identical, differing only in the choice of profile $\bar u$.

\smallskip
3. Near a bounding poly-cycle, it is clear that $dT/da >0$,
indeed $dT/da\to +\infty$ as the boundary is approached.  Thus,
$\gamma>0$, and in combination with the observation of Remark \thmref{Id}
that $(\partial \bar u_m/\partial u_-)\to I$ and $\det df(u_-)<0$
in the homoclinic case, we find that $\sgn \Gamma<0$, recovering
the result of instability obtained by Gardner's technique together
with the homoclinic instability results of [GZ,Z.1].

It is interesting to note that our periodic instability result for the
homoclinic limit makes no requirement on the shape of the 
limiting homoclinic, whereas the general homoclinic instability results 
of [GZ] require that the orbit be {\it convex} in the vicinity of its vertex.
Indeed, in the nonconvex case, the stability index of the homoclinic
orbit is positive, 
consistent with stability, an apparent contradiction with the
Remark 2 just above.
%
In the quadratic-flux case considered in [GZ], it was shown that
all homoclinic orbits are (globally) convex, hence this 
situation does not arise.
However, as pointed out in [Z.4], it can certainly arise in the 
van der Waals models studied in [Z.1] 
(described in the following section, just below).
\footnote{
The possibility of nonconvex orbits was mistakenly not considered
in [Z.1]; however, it can be seen that they are unstable despite
their positive stability index, by the final remark of that paper.
Specifically, in the notation of the reference, 
the zero eigenfunction $\bar u_x$ 
of the Sturm-Liouville operator $M$ in the case of a locally
nonconvex orbit has two nodes, 
so must be at least the third eigenfunction of $M$.  
It follows that there exists some combination $v$ of the
first two eigenfunctions such that $\int v=0$; hence
$v=DV$ for some $V\in L^2$, and also $\langle DV,MDV\rangle<0$,
giving the result.}

A closer inspection quickly resolves this apparent paradox.
For, it is readily seen that a planar homoclinic orbit that is nonconvex
at its vertex saddle must enclose a {\it second}, interior poly-cycle,
containing the interior branches of its stable and unstable manifolds.
Periodic orbits are thus bounded by {\it both} the interior
poly-cycle and the exterior homoclinic, so that the argument for
Remark 2 does not apply.
Indeed, considering the simplest case that the interior poly-cycle
is a second, locally convex homoclinic, and applying the ``multi-bump'' 
heuristics of Remark \thmref{multi} and Remark 3 above it, 
we obtain, formally, the correct prediction of an
{\it odd} multiplicity of unstable eigenvalues for the periodic orbit, 
as the sum of even and odd multiplicities.
\medskip

Combining the results of Theorem \thmref{6}, 
Corollary \thmref{exp},
and Lemma \thmref{sign}, we obtain, finally, the composite stability test:

\proclaim{Corollary \thmlbl{hamconditions}}
Let \text{\rm (H0)--(H3)} 
hold, and suppose also that the traveling-wave 
ordinary differential equation \text{\rm \eqnref{1.1}} is planar Hamiltonian for the value of $s$
associated with solution $\bar u$, with $B\equiv \text{\rm constant}$.
Then, provided that $(\partial \bar u_m/\partial u_-)$
is invertible, both
$$
\sgn \big( dT/da\big)
 \det \big( df(u_-)(\partial \bar u_m/\partial u_-)^{-1} \big) \ge 0
\eqnlbl{stab1r}
$$
and
$$
\sigma\big( df(u_-)(\partial \bar u_m/\partial u_-)^{-1} \big)
\, \text{\rm real}
\eqnlbl{stab2r}
$$
are necessary for weak spectral stability of $\bar u$, where
$dT/da$ is the rate of change of period with respect to amplitude,
with amplitude measured in the direction of the outward normal to the
orbit $\{\bar u(x)\}$.
If $\det (\partial \bar u_m/\partial u_-)= 0$, 
on the other hand, then $\bar u$ is spectrally unstable.
(Recall, $dT/da\ne 0$ by \text{\rm (H3')} combined with Lemma 
\thmref{sign}.)
\endproclaim

\bigskip
\sectionnumber=7
\theoremnumber=0
\equationnumber=0
\smallskip
\TagsOnLeft
\specialhead
\S7.\quad 
Calculations for example systems
\endspecialhead

Using the complementary stability conditions of Corollary 
\thmref{hamconditions}, we now derive instability results for the 
two classes of model, planar Hamiltonian systems that were considered in
[Z.1] and [GZ]: van der Waals gas dynamics
with artificial dispersion-viscosity, and the class
of planar Hamiltonian models with quadratic flux functions.
These are directly relevant to the issue of oscillatory
pattern formation, being prototypes for the two kinds of
systems in which these patterns have been observed numerically.
\medskip

{\bf Van der Waals gas dynamics.}
The viscous-capillary p-system,
$$
\eqalign{
v_t - u_x &= \varepsilon_1 v_{xx} \cr
u_t + p(v)_x &= \varepsilon_2 u_{xx},
}
\eqnlbl{psys}
$$
with nonmonotone stress relation $p$, has been studied by several
authors as a model for dynamical phase transitions in compressible
van der Waals fluids/solids undergoing isothermal motion (see, for
example, [Ja,Sh.1--4,Sl.1--5,ST]).  Here,
$\varepsilon_1 > 0$ and $\varepsilon_2 >0$ are related to the
coefficients of viscosity and capillarity of the medium
($\varepsilon_1+\varepsilon_2$ and $\varepsilon_1\varepsilon_2$, respectively), 
$v$ denotes specific volume/strain, 
$\tilde u := u + \varepsilon_1 v_x$ is the velocity of the medium,
and $p(v)$ denotes pressure/stress. 
A typical stress relation is $p = -W ' (v)$, where
$$
W(v) = \frac{1}{2} (1-v^2)^2 
\eqnlbl{double}
$$
is the standard ``double-well'' potential.

In physical, $(v,\tilde u)$ coordinates, equations \eqnref{psys} 
take the form
$$
\aligned
v_t-\tilde u_x&=0,\\
\tilde u_t + p(v)_x&= (\varepsilon_1+\varepsilon_2)\tilde u_{xx}
- (\varepsilon_1 \varepsilon_2)v_{xxx},\\
\endaligned
\eqnlbl{tildeeqns}
$$
which is sometimes more natural for computations.
\medskip


\medskip
{\it Existence.}
The traveling-wave ordinary differential equation associated 
with \eqnref{psys}
(for periodic and front- or pulse-type waves alike) is
$$
\pmatrix v \\ u \endpmatrix '=
\pmatrix \frac{1}{\varepsilon_1} & 0 \\ 0 & \frac{1}{\varepsilon_2} \endpmatrix
\left[ \pmatrix -u+u_- \\ p(v)-p(v_-) \endpmatrix - s 
\pmatrix v-v_- \\ u-u_- \endpmatrix \right],
\eqnlbl{vanode}
$$
where $s$ denotes the speed of propagation of the wave.
For speed $s=0$, this becomes a Hamiltonian system with Hamiltonian
$$
\aligned
H(v,u,v_-,u_-)&:=
\frac{1}{2\varepsilon_1}(u-u_-)^2- \frac{1}{\varepsilon_2}P(v)\\
&=\frac{1}{2\varepsilon_1}(u-u_-)^2- \frac{1}{\varepsilon_2}
\big(W(v_-)-W(v) +W'(v_-)(v-v_-)\big),
\endaligned
\eqnlbl{Hamiltonian}
$$
where
$$
P(v):=\int^{v}_{v_-} \big(p(v_-)- p(z)\big) \, dz.
\eqnlbl{potential}
$$

For nonmonotone $p$, the orbits of \eqnref{vanode}, 
corresponding to level sets of $H$,
will on some range of $v_-$ include one-parameter
families of periodic orbits bounded by heteroclinic or homoclinic cycles.
For example, in the case of \eqnref{double}, there appears
for $v_-=\pm 1$ a $2$-cycle of heteroclinic orbits connecting
$(v_-,u_-)$ to $(-v_-,u_-)$ in either direction,  respecting
the vertical symmetry of $H$ about the line $u=u_-$; within
this cycle is a one-parameter family of periodic orbits converging
in the small-amplitude limit to the nonlinear center $(0,u_-)$.
For $v_-$ in the range $[0.55,1)$ or $(-1,-0.55]$, there appears 
a single homoclinic orbit with vertex at $v_-$, likewise enclosing
a one-parameter family of periodic orbits about a nonlinear
center $(v_c,u_-)$.  
Together, these constitute all periodic solutions of
\eqnref{vanode} with $s=0$ for the double-well case \eqnref{double}.

We may deduce by energy considerations that traveling periodic 
solutions of \eqnref{psys} must in fact be stationary, so that the Hamiltonian
solutions just described constitute {\it all} periodic solutions
of the van der Waals system \eqnref{psys}.
For, as noted in [Se.1--2], the flow of \eqnref{psys}, 
for periodic solutions of period $T$, serves to decrease the 
mechanical energy
$$
\tilde E(v):= E(v) + 
\int_{0}^{T} \frac{1}{2}\tilde u^2 \, dx 
= E(v) + 
\int_{0}^{T} \frac{1}{2}(u + \varepsilon_1 v_x)^2 \, dx, 
\eqnlbl{etilde}
$$
where 
$$
E(v) := \int^{T}_{0} \left( \frac{\varepsilon_1
\varepsilon_2}{2} |v'|^2+P(v)\right)dx 
\eqnlbl{energy}
$$
denotes the Cahn Hilliard/van der Waals energy for the associated
equilibrium problem.
More precisely,
$$
d\tilde E/dt= -(\varepsilon_1+\varepsilon_2)\int_0^T |\tilde u_x|^2 dx,
\eqnlbl{ederiv}
$$
as may easily be derived from \eqnref{tildeeqns}.
Thus, we may deduce that $\tilde u\equiv \text{\rm constant}$
for any periodic traveling-wave solution of \eqnref{psys}, 
from which we immediately obtain $0\equiv v_t=-sv_x$.  
It follows that either $v$ is identically
constant, in which case $u$ is as well, or else $s=0$;
in either case, we find as asserted that the wave must be stationary.

Moreover, one readily finds for $\tilde u_x=0$ that any solution
of \eqnref{tildeeqns} must satisfy the Euler-Lagrange equations
$$
\varepsilon_1\varepsilon_2 v_{xx}= p(v_-)-p(v)
\eqnlbl{eulerlagrange}
$$ 
for the equilibrium variational problem
$$
\min_v E(v),
\eqnlbl{var}
$$  
yielding back the stationary-wave ordinary differential equation
\eqnref{vanode} with $s=0$.

\medskip
{\it Stability.}
Using the result of Corollary \thmref{hamconditions}, we readily
obtain:

\proclaim{Theorem \thmlbl{vanstab}}
Periodic orbits of \text{\rm \eqnref{psys}} are unstable
whenever $dT/da> 0$, in particular, in the large-amplitude 
limit as they approach either a bounding homoclinic or two-cycle.
Orbits for which $dT/da<0$ are unstable if the spectrum
of $(\partial \bar u_m/\partial u_-)^{-1}df(u_-)$ is nonreal.
In the small-amplitude limit as
orbits approach a nonlinear center, they are unstable 
regardless of the value of $\sgn dT/da$.
\endproclaim

{\bf Proof.} Clearly, (H0)--(H2), hold, whence the
results of Corollary \thmref{hamconditions} apply whenever
(H3') holds, or $dT/da\ne0$.
By the final assertion of Corollary \thmref{hamconditions}, 
we may without loss of generality take 
$\det( \partial \bar u_m/\partial u_-)\ne 0$.
The second assertion then follows by either of \eqnref{stab1r}, \eqnref{stab2r},
the former because in this case the determinant of the real matrix
$(\partial \bar u_m/\partial u_-)^{-1}df(u_-)$ 
must be positive, as the product of complex conjugates.
To establish the first assertion, we make use of the
structure of equations \eqnref{psys}.

By direct computation, we obtain
$$
df(u_-,v_-)=\pmatrix 0 & -1 \\ p'(v_-) & 0 \endpmatrix.
\eqnlbl{A}
$$
Moreover, from translational invariance of \eqnref{psys} with 
respect to $u$, in combination with the reflective symmetry 
$u\to -u$, $x\to -x$, we find without any computation
that $\partial \bar u_m/\partial u_-$ must be of diagonal form
$$
\pmatrix \alpha &0\\
0 & 1\\
\endpmatrix,
\eqnlbl{alph}
$$
with a single undetermined quantity $\alpha\ne 0$. 
Thus, assuming that
$\det (\partial \bar u_m/\partial u_-)\ne0$, or $\alpha\ne 0$,
we find that the matrix
$$
(\partial \bar u_m/\partial u_-)^{-1}df(v_-,u_-)
=
\pmatrix 0&  -1/\alpha \\
p'(v_-)& 0\\
\endpmatrix
\eqnlbl{struct}
$$
is trace-free, whence we may conclude that its spectrum is
real if and only if its determinant is negative.
But, this implies that \eqnref{stab1r} and \eqnref{stab2r} are 
mutually exclusive when $dT/da>0$, yielding spectral instability.

Finally, the third assertion follows simply by continuity of the
Evans function, without invoking (H3') or Corollary \thmref{hamconditions},
once we observe that the limiting, constant-coefficient equations
at a nonlinear center are always unstable, since the trace-free
real matrix $df(v_c,u_c)$ must necessarily have pure imaginary eigenvalues
(or else $(v_c,u_c)$ would instead be a saddle).
\qed
\bigskip

{\bf Remark \thmlbl{doublewell}.}
In the special case of the double-well potential \eqnref{double},
monotonicity $dT/da>0$ holds for {\it all} periodic van 
der Waals orbits. 
This follows as a straightforward application of a monotonicity theorem 
of Schaaf (see [Sc, p. 102], or [LP, Theorem 7.1, p. 331]),
which asserts that monotonicity holds for the nonlinear oscillator 
$$
v''+\mu(v)=0,
\eqnlbl{nonlin}
$$
provided: 
(i) $\mu'>0$ implies $5(\mu'')^2 - 3\mu' \mu'''>0$; and,
(ii) $\mu'=0$ implies $\mu \mu''<0$.
Noting that the van der Waals traveling-wave ordinary differential equation
is of form
\eqnref{nonlin}, with $\mu(v)=p(v)-p(v_-)$, we readily
verify (i) and (ii) in the double-well case $p(v)=2v(1-v^2)$,
by direct computation.  In particular, note that
(i) follows simply from the fact that $p'''(v)=-12v<0$.
This finding is quite significant,
since it shows that there is at least one model exhibiting oscillatory
pattern formation for which no stable periodic waves exist, thus
indicating the presence of some alternative mechanism for the formation
of oscillatory patterns.

\medskip
{\bf Remark \thmlbl{gurtin}.}  
By Remark \thmref{sufficient}, in conjunction with Lemma \thmref{sign},
we find that small-amplitude periodic waves
are {\it stable with respect to periodic perturbations}
(of the same period $T$) if and only if $dT/da>0$;
in particular, they are stable for the double-well
potential, as discussed just above.
Likewise,
numerical evaluation of $\partial\bar u_m/\partial u_-$ 
for large-amplitude
waves in the two-cycle configuration yields $\Gamma>0$, suggestive of
stability under $T$-periodic perturbation, a conjecture that is further
supported by numerical approximation of the $(k=0)$-spectrum of $L$;
see Appendix A.
On the other hand, we have already seen analytically that 
large-amplitude waves in the homoclinic case are {\it unstable} with
respect to $T$-periodic perturbations.

These results are quite interesting from the point of view of the 
associated variational problem \eqnref{var}.  For, they imply that 
periodic solutions may or may not be stable critical points of the 
energy $E$, depending on the details of their structure.
This is in sharp contrast to the results of [CGS.1--2,GM]
for the same problem with Neumann in place of periodic boundary conditions,
in which critical points are seen to be stable if and only if they
are {\it monotone}.
In both cases, stability is interpreted with respect to perturbations
conserving mass. Note, in the conservation law setting, that 
stability is at best {\it orbital}, in the sense of convergence to
the three-parameter family of nearby periodic solutions of the same
period $T$, with the two conserved quantities given by the mass of 
the perturbed solution determining the limiting solution up to phase 
shift (translation).

To put things slightly differently, stability in the periodic 
setting is not determinable as in the Neumann setting by 
Sturm-Liouville-type considerations alone, but must indeed
be estimated by some such computation as we have carried out.
Similarly, stability for the variational problem on the whole
line is more subtle than that for the (Neumann) problem on
the bounded interval, requiring further analysis
as in [Z.1].

\medskip
{\bf Quadratic-flux models.}
Planar systems \eqnref{1} with quadratic flux $f$ have been studied
as qualitative models for multiphase flow in porous media
near an ``elliptic boundary'' where characteristics of the
flux Jacobian $df(u)$ coalesce; see, e.g., [SSh,MP]. 
Here, following [GZ], we consider the restricted
class of $2 \times 2$ systems with $f$ quadratic
and $B$ constant for which, additionally, the traveling-wave ordinary
differential equation 
$$
u'=B(u)^{-1}[f(u)-f(u_-)-s(u-u_-)]
$$
is {\it Hamiltonian}, i.e.,
$$
tr(B^{-1}(f'(u)-sI)) \equiv 0
$$
for some choice of s. This is a codimension-2 subclass of the
quadratic-flux models, as parametrized by the coefficients of
$B$ and $f$. 


Generically, such systems can be reduced by 
an affine change of variables to the canonical form
$$
\pmatrix u \\ v \endpmatrix_t + B\pmatrix \frac{\varepsilon}{2}
v^2-\frac{1}{2}u^2+v \\ uv \endpmatrix_x
=B \pmatrix u \\v \endpmatrix_{xx}
\eqnlbl{canon}
$$
where $\varepsilon=\pm1$. For $\varepsilon=1$ and $s=0$, the 
traveling-wave ordinary
differential equation then reduces to the Hamiltonian system
$$
\pmatrix u \\v \endpmatrix'=\pmatrix 
\frac{1}{2}v^2-\frac{1}{2}u^2+v \\ uv \endpmatrix - \pmatrix
\frac{1}{2}v_{-}^2 - \frac{1}{2}u_{-}^2+v_- \\ u_-v_- \endpmatrix,
\eqnlbl{flux}
$$
with Hamiltonian
$$
H(u,v)=\frac{1}{2}(\frac{1}{3}v^3-u^2v+v^2)-(\frac{1}{2}v_-^2v
-\frac{1}{2}u_{-}^2v+v_-v-u_-v_-u)\equiv H(u_0,v_0)
\eqnlbl{hamiltonian}
$$
preserved along orbits [GZ]. 

Explicit computation of the phase portrait of \eqnref{flux}
using \eqnref{hamiltonian} reveals the
existence of a unique three-cycle configuration for 
the special parameter values
$$
(u_-,v_-)=
(\pm \sqrt{3}/2,0), \,
(0, -1.5), \,
(0,-.5), 
\eqnlbl{params}
$$
corresponding, respectively to the three saddle equilibria at
the vertices of the three-cycle, and the single nonlinear center
contained within.  
Enclosed within the three-cycle is a nested one-parameter
family of periodic orbits.
The heteroclinic orbits making up the three-cycle
lie on straight lines, forming an equilateral triangle; in fact
the entire phase portrait has the same triple symmetry about
the center as does the bounding triangle.
The three states \eqnref{params} at the vertices of the triangle 
may be regarded as analogous to {\it pure phases} in 
physical models for three-phase flow, and the interior of the
triangle to the physical state space of volume fractions thereof [AMPZ.3].

This configuration bifurcates as $(u_-,v_-)$ is varied into various
one- and two-cycle configurations, each cycle likewise enclosing
a nested one-parameter family of periodic orbits.
Where $dT/da\ne0$, these are the only nearby periodic orbits to be
found in parameter space,  where $s$ is now allowed to vary;
however, we have not so far been able to rule out the possibility
of further periodic orbits with large speed $s$.

In numerical experiments carried out in [AMPZ.5], systems
\eqnref{canon} were seen to exhibit the same sort of oscillatory
pattern formation discussed in [FL.1--2,CP], with patterns
composed of motifs with $(u,v)$ profiles lying near either
the  three-cycle, or the various two-, and one-cycle
configurations that bifurcate from it as $(u_-,v_-)$ is varied
from \eqnref{params}.  Thus, the relevant issue from the
point of view of pattern formation seems to be the stability
or instability of large-amplitude periodic waves for parameters
near the three-cycle values \eqnref{params}.
Regarding this problem, we have a complete solution;
indeed, the additional symmetries of the three-cycle case
allow us to obtain results for {\it general amplitude} waves: 

\proclaim{Theorem \thmlbl{quadstab}}
For $B$ constant, symmetric, and positive-definite,
periodic orbits of \eqnref{canon} are unstable for 
$(u_-,v_-)$ lying near the three-cycle values \eqnref{params},
independent of amplitude, provided that $dT/da\ne 0$.
Likewise, so long as $dT/da\ne 0$, large-amplitude waves
are unstable for $(u_-,v_-)$ corresponding to a homoclinic 
configuration in the phase portrait of \text{\rm\eqnref{flux}};
moreover, small-amplitude waves are unstable for any
phase configuration.
\endproclaim

{\bf Proof.} Clearly, (H0)--(H2), hold, whence, again, the
results of Corollary \thmref{hamconditions} apply whenever
(H3') holds, or $dT/da\ne0$.
Taking $(u_-,v_-)$ to be the nonlinear center $(u_c,v_c)=(0,-0.5)$,
we find using the triple symmetry of the phase portrait that
$\partial (\bar u_m, \bar v_m)/\partial (u_-,v_-)$
must be a multiple of the identity, for any periodic orbit whatsoever.
More precisely, we note that, for $(u_-,v_-)=(u_c,v_c)=(0,-0.5)$, 
traveling-wave ordinary differential equation
 \eqnref{flux} has reflective symmetry about each
of the lines from $(u_c,v_c)$ through the vertices
of the surrounding three-cycle; moreover, each individual
symmetry is preserved as $(u_-,v_-)$ is varied along the 
respective line of symmetry.  It follows that, under
perturbations of $(u_-,v_-)$ in direction $\nu$,
where $\nu$ is the direction vector of a line of symmetry,
must be confined to the same line of symmetry, i.e., the 
resulting perturbation of the mean $(\bar u_m,\bar v_m)$ 
must lie in the same direction $\nu$, for any choice of initial 
periodic orbit $(\bar u,\bar v)$. 
In other words, each of the direction vectors of the three
lines of symmetry is an eigenvector of 
$\partial (\bar u_m,\bar v_m)/\partial (u_-,v_-)$,
whence 
$\partial (\bar u_m,\bar v_m)/\partial (u_-,v_-)$ must be a 
(real) multiple of the identity.

Thus, 
$$
df(u_-,v_-)
\partial (\bar u_m, \bar v_m)/\partial (u_-,v_-)^{-1}
$$
is simply a real multiple of 
$$
df(u_-,v_-)= B\pmatrix
0 & 0.5\\
-0.5 & 0\\
\endpmatrix.
\eqnlbl{val}
$$
But, this matrix is similar to the anti-symmetric matrix
$$
B^{1/2}
\pmatrix
0 & 0.5\\
-0.5 & 0\\
\endpmatrix
B^{1/2},
$$
hence must have pure imaginary spectrum.
(Note: in particular, it is again trace-free, as
in the proof of Theorem \thmref{vanstab}.)
We thus obtain instability by \eqnref{stab2r},
for waves of arbitrary amplitude. 
The result for $(u_-,v_-)$ near $(0,-0.5)$
follows by the same argument, together with
continuity of $\partial (\bar u_m, \bar v_m)/\partial (u_-,v_-)$.

The large-amplitude homoclinic result has been established already 
by several different techniques.
Likewise, the small-amplitude result follows from the
instability of the constant coefficient limit, as discussed
previously in Remark \thmref{smallamp}.
\qed
\bigskip

{\bf Remark \thmlbl{quadper}.}  
The proof of Theorem \thmref{quadstab} shows also that, for $(u_-,v_-)$
near the $3$-cycle values \eqnref{params} and $dT/da>0\:$-- in particular,
for large-amplitude periodic orbits approaching the $3$-cycle,
the stability index $\Gamma$ is {\it positive}, suggestive of stability
with respect to $T$-periodic perturbations.
This conjecture is supported by numerical approximation
of the $(k=0)$-spectrum, which also indicates stability 
with respect to $T$-periodic perturbations of large-amplitude
periodic solutions in the $2$-cycle limit; see Appendix A.
This is in contrast to the $T$-periodic {\it instability} observed
in the large-amplitude homoclinic limit for both van der Waals and
quadratic flux models.

\medskip

{\bf Remark \thmlbl{3cycle}.}  The general features of the
global phase portrait structure described above for the Hamiltonian
models considered here have been shown in [AMPZ.4] 
to hold, more generally, for all $2\times 2$ quadratic-flux models 
of type I under the viscous classification scheme of Hurley [HP].
In particular, there exist unique $(u_-,v_-,s)$ such that
the phase portrait of the traveling-wave ordinary differential equation
contains a three-cycle,
of which the connecting heteroclinic orbits must lie along straight
lines;  moreover, it can be shown by affine reduction to canonical
form that the corresponding ordinary differential equation
is quasi-Hamiltonian [CL].
It follows that the three-cycle encloses a nested one-parameter
family of periodic orbits, as in the ``pure'' Hamiltonian case
considered above.  Likewise, formation of oscillatory patterns
has been observed numerically for such systems [AMPZ.5], with
patterns again lying near the special three-cycle configuration.
\medskip

\bigskip
%
%
%
%
\bigskip
\sectionnumber=8
\theoremnumber=0
\equationnumber=0
\smallskip
\TagsOnLeft
\specialhead


Appendix A.\quad Numerical experiments
\endspecialhead                       

In this appendix, we describe numerical methods for 
the location of the zero set of the periodic Evans function $D(k,\lambda)$,
or equivalently the spectrum of a linear operator $L$ with
periodic coefficients, and, in particular, 
the well-conditioned {\it numerical verification of instability}.  
We demonstrate these methods for the example systems considered in 
Section 7, amplifying and (in the case of the quadratic flux models)
extending the analytical results there obtained.
The real advantage of the numerical approach, of course, is that
it is applicable to general models.
\medskip

{\bf The algorithm.} We use the following basic strategy
for numerical verification of instability:
\smallskip

1. Determine the range of existence of periodic orbits, and 
(numerically) solve for a representative sampling of profiles.

\smallskip
For each profile:

\smallskip
2. Find spectral curves $\lambda(k)$, i.e., zero-level sets of
$D(\lambda, k)$; 

3. Choose $k=k^*$ giving a maximum (positive) value of 
$\R(\lambda(k))$ and a numerically advantageous closed contour $\Gamma$
contained in the strictly positive half-plane $\R \lambda>0$ and
containing $\lambda(k^*)$.

4. Perform a winding number calculation in $\lambda$, with $k=k^*$
held fixed, to establish that there indeed
lies a root of $D(k^*, \cdot)$ in the vicinity of the original
numerical approximation of $\lambda(k^*)$ in steps 2 and 3.

\bigskip
The point of steps 3 and 4 is that the curves found in step 2
are ``graphical'' approximations of the zero-level sets, which,
though extremely illuminating, 
do not give conclusive information about existence of zeroes.
Step 4 by contrast is both accurate and well-conditioned,
and could be used as the basis for numerical proof; however,
it requires an a prior guess for $k^*$, $\lambda(k^*)$.
\medskip

{\it Choosing $k^*$ and $\Gamma$.}
By a standard G\"arding type energy estimate on \eqnref{2}, the zeros of
the Evans function must lie in the truncated wedge 
$$
V=\{\lambda \in \Bbb{C} : \Im(\lambda)+\R(\lambda) \le r, 
\Im(\lambda)-\R(\lambda) \ge -r, \R(\lambda) \le r/4 \}
$$
where $r=\|f(\xi) \|^{2}_{\infty} / \|B\|_{2}, f(\xi)=\|A(\xi)\|_2$. 
Obviously, this is not
a bounded region but it does give us a bound on the unstable
zeros of $D$. 
We find spectral curves/zero sets of $D$ numerically 
and observe that those are contained in $V$. We can see also agreement
between our numerics and the analysis in Sections 1--7 (see discussion under
{\it Applications}, below).

It is more tricky to find the zero sets of the
two-parameter periodic-type Evans function than of one-parameter
front-type Evans function. Adding to the difficulty is the
fact that $D(\lambda,k)$ vanishes to $(n+1)$st order in $\lambda$
at $(0,0)$, where $n$ is the dimension of variable $u$,
so that there are $n+1$ spectral curves bifurcating through
the origin. The apparently straightforward step $2$ in the algorithm
is thus in fact a tricky problem in graphical display. 
To find spectral curves/zero sets of $D$, we use M\"{u}ller's method
which can approximate complex roots and converges to the root for 
any initial approximation choice. We also use deflation 
to achieve a reasonable result near the bifurcation point 
$(k, \lambda)=(0, 0)$. 

Since the argument principle will be used to count the zeros
inside some contour $\Gamma$, care must be taken that $\Gamma$ 
itself does not contain any zeros.
We avoid the 
origin and high multiplicities by looking at $k \ne 0$. 
A natural choice is to choose the value
of $k$ which gives the maximum value of the positive 
$\R(\lambda)$
from the spectral curves, 
and for simplicity choose $\Gamma$ as the rectangle 
whose center is that
maximum positive value of $\R(\lambda)$ in the right
half plane. In this process, we find that a practical upper
bound on the size of the contour is enforced by the need to obtain good
error estimates, which degrade rapidly for high frequencies.
This particular $k$ allows good numerics. 
When the computed winding number exceeds zero, the
wave being tested is unstable; otherwise the wave is stable.

\medskip
{\it Applying the Argument Principle.}
Here we follow the approach of [Br.1--2].
The first problem in calculating the winding number of $D[\Gamma]$
lies in evaluating $D(\lambda, k)$ for a fixed value of $k$. 
We use the maximal value $k=k^*$ described above.  
The second problem lies in partitioning $\Gamma$ into subcontours
$\Gamma_{i}$ such that $D[\Gamma_{i}]$ lies in a slit plane; 
a slit plane being a subset of $\Bbb{C}$ of the form $\Bbb{C}
\backslash
\{r\alpha : 0 \ne \alpha \in \Bbb{C}, 0\le r \in \Bbb{R} \}$. This 
will ensure that the calculated winding number is correct. We
partition $\Gamma$ so that $D[\Gamma_{i}]$ is contained in a half
plane whose boundary passes through the origin. 

For $w, z \in \Bbb{C}$, define $\langle w,z \rangle$ to
 be the dot product between
the real 2-vectors $(\R(w), \Im(w))$ and $(\R(z), \Im(z))$ :
$$
\langle w,z\rangle=\R(w)\R(z) + \Im(w)\Im(z).
\eqnlbl{inner}
$$
If $\Gamma(s_{i})$ is the starting point of $\Gamma_{i}$, 
it suffices to find $s_{i+1}$ and $\beta \in [0,1)$ such that
$$
\langle D(s_{i}+\Delta s), D(s_{i})\rangle  > \beta |D(s_{i})|
\eqnlbl{cond}
$$
for $0 \le \Delta s \le s_{i+1}-s_{i}$. In other words, the 
component of $D(s_{i}+\Delta s)$ in the direction of $D(s_{i})$
should be a positive multiple of $D(s_{i})$. 
Strictly speaking, the right side of inequality \eqnref{cond} 
could simply be $0$ in practice. However, to allow for numerical
error,
a positive value of $\beta$ should be chosen. We choose $\beta=1/2$
over all.
This condition allows checking whether a given partition is
appropriate.

\medskip

{\it Error Estimates.}
In contrast to the front case [Br.1--2], 
computation of $D$ for moderate values
of $| \lambda|$, 
given exact values for coefficients of $L$, is straightforward and
numerically well-conditioned, and for practical purposes may be
regarded as exact. The main approximation error comes rather from
the computation of coefficients $A$ and $B$ in \eqnref{1.3}
through solution of the traveling-wave ordinary differential
equation \eqnref{1.1}, which is likewise a standard and relatively
well-conditioned problem; here we use a fourth order Runge-Kutta
method, for which analytical error estimates may be found, for
example, in [SB].
For large $|\lambda|$, it is rather the $\CalO(|\lambda|)$
exponential growth rate in the eigenvalue ordinary differential
equation that dominates
the growth of errors; indeed, it is this consideration that
effectively limits the size of our contours in the winding number
calculations described above.
However, again, the relevant analytic error bounds are standard
and straightforward.

We therefore omit the discussion of analytical error estimates,
instead carrying out a numerical convergence study.
Specifically, we use Richardson extrapolation to compute ``exact values" 
to which we compare when we calculate relative errors and to check the 
convergence of $D$ as the step size $h$ becomes small.
The resulting estimated errors are seen to be quite small, 
consistent with the above discussion.

\medskip

{\bf Numerical experiments.} We now describe the results
of our numerical experiments for the two classes of systems discussed 
in Section 7.

{\it Van der Waals gas or solid mechanics.}
Our first set of experiments is for 
the viscous-capillary p-system \eqnref{psys} with standard
double-well potential \eqnref{double}.
As discussed in Section 7, periodic solutions of this model
are of necessity stationary, satisfying the Hamiltonian
system \eqnref{vanode} with $s=0$, and this system features
for $v_-=\pm 1$
a one parameter family of heteroclinic two-cycles, 
indexed by $u_-$ due to Galilean invariance in the speed $u$, 
in which there lie nested families of periodic orbits.
Likewise, there exists a two parameter family of 
homoclinic cycles for $v_-$ in range $[0.55, 1)$, and $u_-$ again
arbitrary by Galilean invariance,
in which there lie nested families of periodic orbits.
Without loss of generality (using translation
invariance), we may
fix $u_-=0$, reducing our study to a single $2$-cycle and a
one-parameter
family of homoclinic orbits, and their  enclosed periodic orbits, see
Figures
6A-B below.

First, we use a fourth order  Runge-Kutta
method to get a profile $(v,u)$ by solving the Hamiltonian system.
By plotting $v$ versus $u$, we can see clearly 
there exist periodic orbits 
(Fig. 2A) inside the two cycle. 
{{From these orbits, we observe that a monotone relationship
between the periods and the amplitudes (Fig. 2B).
Second, we use the method again to
solve \eqnref{4}.
Finally, we use the LU decomposition method
to get a value of $D$. We also find periodic
orbits and the monotone relationship 
inside the homoclinic orbit numerically.   
We study several periodic orbits in each cycle. 

\midinsert
\captionwidth{2.5truein} 
\hbox{
\vbox to 2.7truein{\hsize=2.5truein\efig van3 y2.2
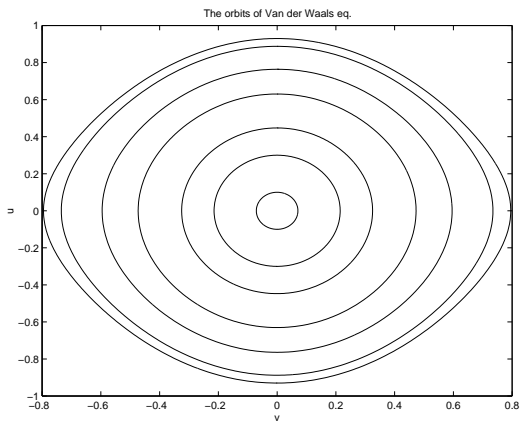
\botcaption{Figure 2a} 
Existence of periodic orbits.
\endcaption
}
\hskip .8truein 
\vbox to 2.7truein{\hsize=2.5truein\efig van35 y2.2
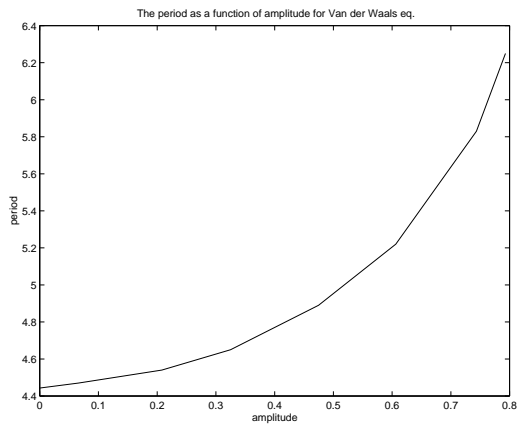
\botcaption{Figure 2b} 
Period as a function of amplitude.
\endcaption
}
} 
\endinsert

\topinsert
\hbox{
\vbox to 2.2truein{\hsize=2.5truein\efig van4 y2.2
}
\hskip .7truein 
\vbox to 2.2truein{\hsize=2.5truein\efig van48 y2.2
}
}
\hbox{
\vbox to 2.2truein{\hsize=2.5truein\efig van5 y2.2
}
\hskip .7truein 
\vbox to 2.2truein{\hsize=2.5truein\efig van6 y2.2
}
}
\hbox{
\vbox to 2.2truein{\hsize=2.5truein\efig van7 y2.2
}
\hskip .7truein 
\vbox to 2.2truein{\hsize=2.5truein\efig van75 y2.2
}
}       
\captionwidth{6truein}
\botcaption{Figure 3}
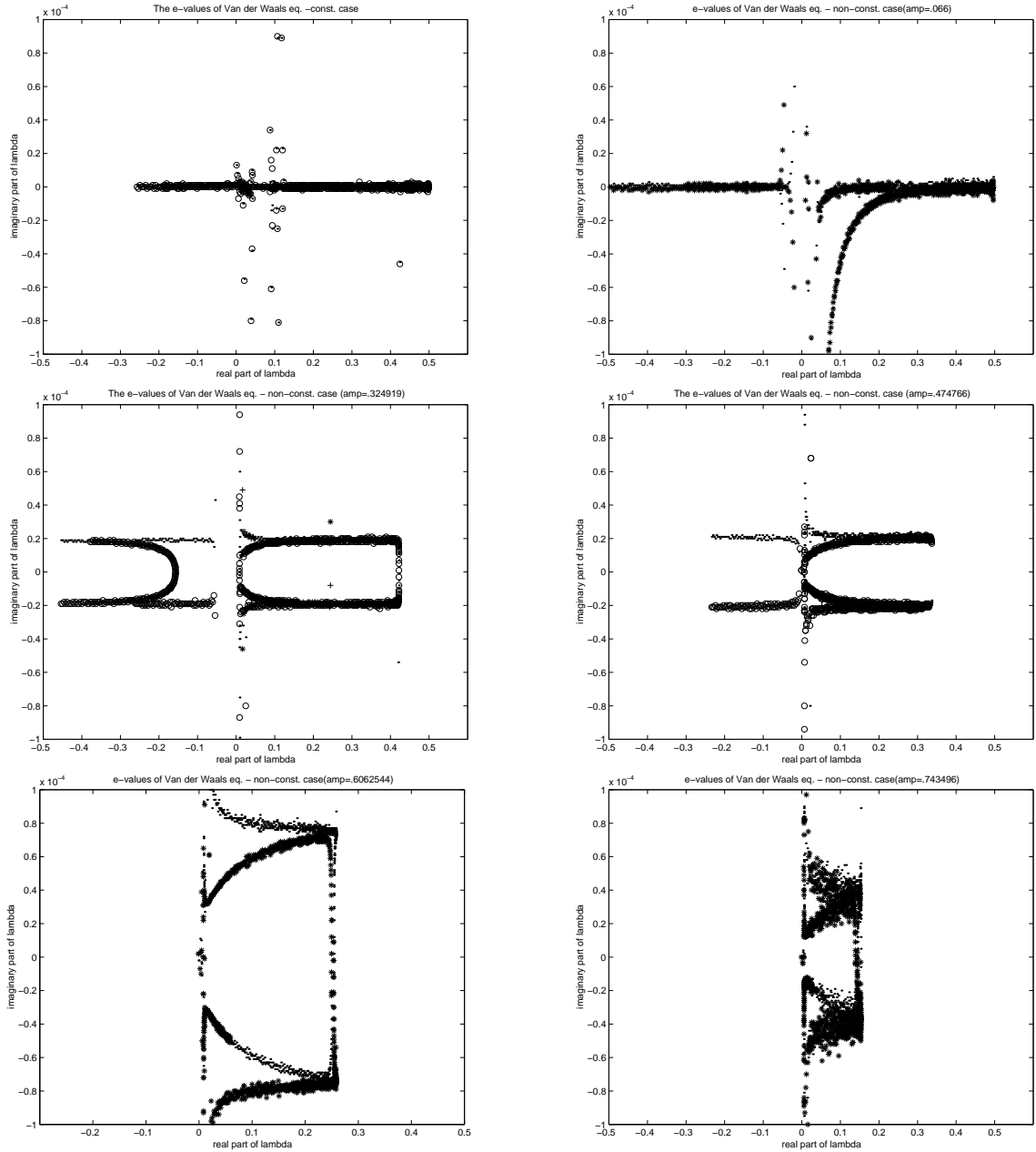
Spectral curves of the Evans function $D$.
\endcaption    
\endinsert

Figure 3 shows spectral
curves of $L$, or zero sets of $D$, for various periodic orbits
within the $2$-cycle depicted in figure 2A,
starting with the (zero amplitude) constant coefficient
case and shows how the spectral curves are changed
as the amplitude becomes large.
We see what differences are made by change of amplitudes.     
These figures are extremely illuminating, showing in a way that
calculations can not the detailed evolution of the spectrum as
the amplitude of the underlying periodic orbit is varied from 
zero (constant coefficient case) to the limiting amplitude as 
orbits approach their bounding polycycle. They both verify and 
amplify the observations obtained earlier by Evans function 
calculations. 
In particular, we immediately see that
instability is rather dramatic, whereas the analytical results
via Taylor expansion about the origin imply only that there 
exist infinitesimally unstable spectra; moreover, this instability
is not confined to small or large amplitude limits, but persists
across all amplitudes. 

This suggests the strategy described above
for numerical verification of instability; namely, to fix $k$ at
the value $k^*$ corresponding to the maximally unstable spectrum
$\lambda^* = \lambda(k^*)$ obtained by graphical examination, then
perform a winding number calculation on a well-chosen contour about
$\lambda^*$. 

As for stability results, Figures 4 and 5 show $\Gamma$ and $D[\Gamma]$
for periodic waves inside a two cycle and a homoclinic orbit
respectively. Here, $D[\Gamma]$ is normalized by the mapping
$D \longmapsto 0.1*(i+1)*D/|D|$ for $i=0,1,...,39$ (Here, $i$ is the
index of $\lambda$ on the contour $\Gamma$.
The winding numbers are greater than $1$
both cases so the periodic waves are unstable.

\midinsert
\hskip 1.5truein
\vbox to 3.5truein{\hsize=2.5truein \efig wind32 y3.3
\captionwidth{5truein}
\botcaption{Figure 4} 
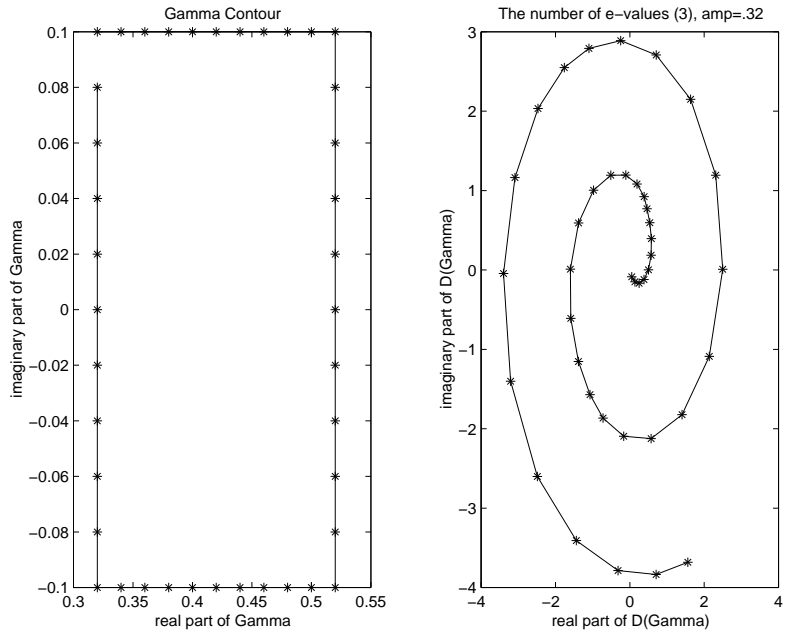
The number of unstable eigenvalues in the two-cycle case.
\endcaption
}
\endinsert
\midinsert
\hskip 1.5truein
\vbox to 3.5truein{\hsize=2.5truein \efig windhomo83 y3.3
\botcaption{Figure 5} 
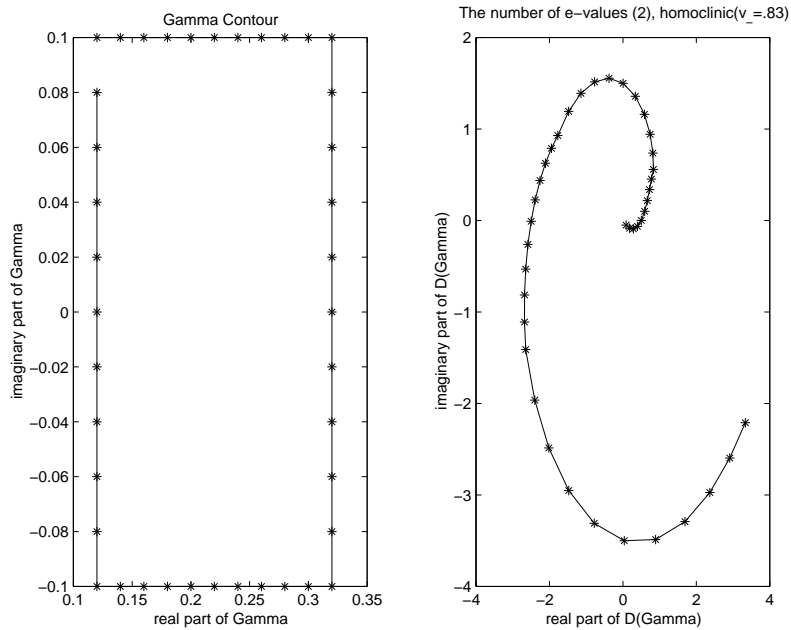
The number of unstable eigenvalues in the homoclinic case.
\endcaption
}
\endinsert

We next perform the error analysis that was mentioned earlier. 
There is no analytic solution of the eigenvalue 
ordinary differential equation \eqnref{4} for the van der Waals equations.
We have to rely entirely on numerical estimates to gauge the
accuracy of the Evans function method.  
In the tables below, third order polynomial extrapolation was
used to compute the values for $\tau=0$ where $\tau$ is
a local truncation error depending on the step size $h$, 
i.e., $\tau := Ch^4$ for some positive constant $C$ (recall that 
the error to lowest order comes from solution of the profile 
equation by the fourth order Runge-Kutta scheme).

\bigskip
$$
\matrix
\tau & D(0.32) & Extrapolated \ Relative \ error \\
10^{-5} & 3.358392 & 0.0121955 \\
10^{-6} & 3.358958 & 0.0120290 \\
10^{-7} & 3.359097 & 0.0119882 \\
10^{-8} & 3.399841 & 0.0000041 \\
. & .  & \\
. & .  & \\
. & .  & \\
0 & 3.399855 &
\endmatrix
$$
\botcaption{Table 1}
$D(0.32)$ for the van der Waals equations
\endcaption
\bigskip

$$
\matrix
\tau & D(0.52-0.1i) & Extrapolated \ Relative
\ error \\
10^{-5} & -4.636345-12.13182i & 0.0082597282 \\
10^{-6} & -4.636959-12.13066i & 0.0082772553 \\
10^{-7} & -4.637134-12.13046i & 0.0006327496 \\
10^{-8} & -4.722007-12.19752i & 0.0000185195 \\
. & . & \\
. & . & \\
. & . & \\
0 & -4.721931-12.19775i &
\endmatrix
$$
\botcaption{Table 2}
$D(0.52-0.1i)$ for the van der Waals equations
\endcaption
\bigskip
\bigskip

Starting with the
extrapolation results, examine the above Tables which describe
an illustrative calculation when the amplitude of periodic 
orbit is $0.32$.
The values chosen
for $\lambda$ are the points on $\Gamma$ closest to and furthest from
the center of $\Gamma$; in practice, these are found to give
the extrema of errors.  The extrapolated values      
are used as exact values in computing extrapolated relative errors.
Note that $\tau=10^{-6}$ is sufficient for calculating $D(\lambda,
k)$ with
approximate estimated accuracy $0.05$.  

Similar calculations were carried out for all mesh points $\lambda$
on the contours for each periodic orbit under consideration. 
For example, in the two-cycle case, we tested all periodic orbits
indicated in Figure 6A. The points in the Figure represent the amplitudes
of each periodic orbit.
{\it Note: translation invariance greatly reduces the dimension of
the problem.} 
Likewise we consider homoclinic profiles 
with vertices $v_- \in [0.55, 1)$ with vertex mesh $0.05$ (see Fig. 6B). 
For 
each homoclinic case, we then computed  
periodics with amplitude mesh $O(0.01)$ where amplitude  
$ \in [0, O(0.1)]$. 
In each case, the winding number was found to be greater than
$1$.
The bounds for the winding number calculations
are found conservatively to be order $10^{-2}$, i.e.,
the worst case Rouch\'{e}'s error over all calculation was
$$
max \frac{|\triangle D|}{|D|} \le C 10^{-2} < < 1/2,
$$
for some constant C (Table). It is comfortably below
the bound $1$ needed to apply Rouch\'{e}s Theorem.
Thus, we have numerical instability across
the whole parameter range of existence.     

\midinsert
\captionwidth{2.5truein}
\hbox{
\vbox to 3.2truein{\hsize=2.5truein \efig g y2.2
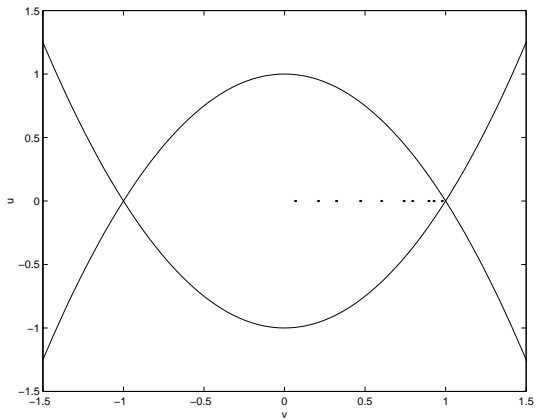
\botcaption{Figure 6a} 
The plotted points represent the amplitudes for which the experiments
were carried out in the two-cycle case.
\endcaption
}
\hskip .8truein
\vbox to 3.2truein{\hsize=2.5truein \efig g1 y2.2
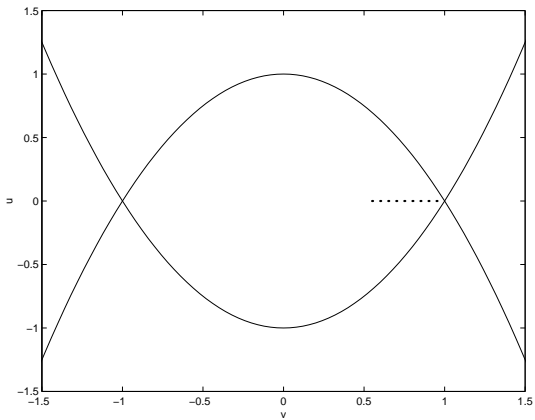
\botcaption{Figure 6b} 
The plotted points represent the amplitudes for which the experiments
were carried out in the homoclinic case.
\endcaption
}           
}
\endinsert

As discussed in the main body of the paper,
there are various interesting transitions of the number of unstable
0-eigenvalues
of the periodic waves with respect to different amplitude.
Recall that the small amplitude limit is the constant coefficient
case. For fixed $k=0$,
the Evans function has four roots at the origin if the amplitude is
zero, and three roots at zero otherwise when $n=2$;
see Remark \thmref{smallamp}.
But, we can track a fourth small root of $D$ as the 
amplitude changes. When 
the amplitude is small, the fourth root is typically in the 
left half plane of $\R (\lambda)$ (Remark \thmref{smallamp}). 
In the homoclinic case, as the amplitude
becomes large, the fourth root moves from the left to the right
half plane crossing the origin (Fig. 7); see Remark
\thmref{Id}. 
For the two-cycle case, no such transition occurs.

\midinsert
\captionwidth{5truein}
\hskip 1.5truein
\vbox to 4.7truein{\hsize=2.5truein \efig trans y4.0
\botcaption{Figure 7} 
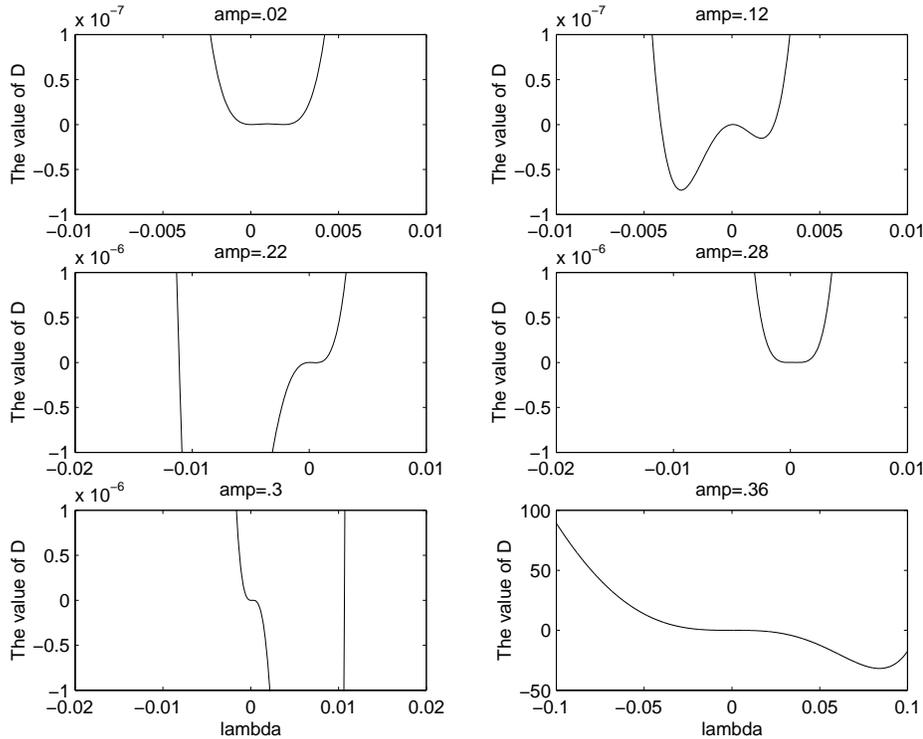
Transition of the number of unstable
0-eigenvalues for the van der Waals equations, homoclinic case.
\endcaption
}
\endinsert

\bigskip
\bigskip
{\it The quadratic-flux models}
We next consider the class of quadratic flux models
\eqnref{canon} with identity viscosity matrix $B=I$,
applying the same methods that we used for the van der Waals equations. 
First, we get a profile $(u,v)$ by solving \eqnref{flux}.
By plotting $u$ versus $v$, we can see clearly that
there exist periodic orbits
(Fig. 8A). From these orbits, we observe that there is a monotone
relationship
between the periods and the amplitudes (Fig. 8B).
Finally, we use the LU decomposition method
to get a value of $D$. We have spectral curves of $\lambda$
with respect to various amplitudes (Fig. 9).                       

\midinsert
\captionwidth{2.5truein}
\hbox{
\vbox to 2.7truein{\hsize=2.5truein\efig flux3 y2.2
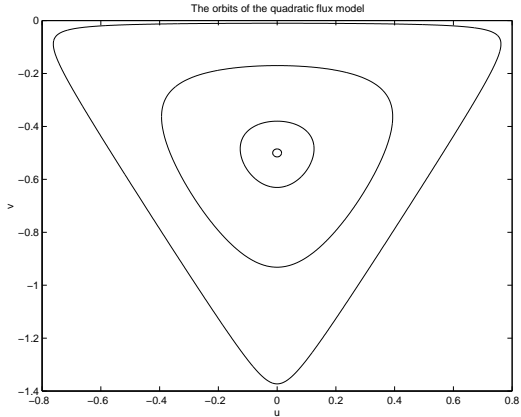
\botcaption{Figure 8a} 
Existence of periodic orbits.
\endcaption
}
\hskip .8truein 
\vbox to 2.7truein{\hsize=2.5truein\efig flux35 y2.2
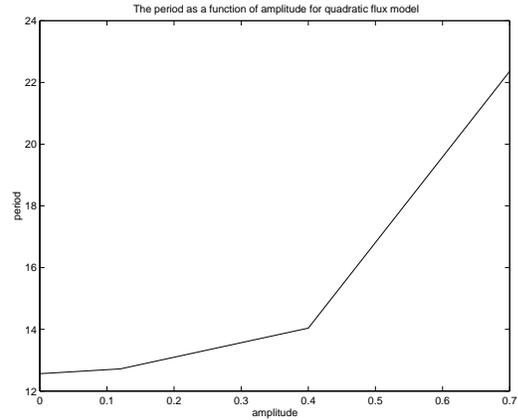
\botcaption{Figure 8b} 
Period as a function of amplitude.
\endcaption
}
}           
\endinsert

\midinsert
\hbox{    
\vbox to 2.2truein{\hsize=2.5truein\efig flux4 y2.2
}
\hskip .7truein 
\vbox to 2.2truein{\hsize=2.5truein\efig flux45 y2.2
}
}
\hbox{
\vbox to 2.2truein{\hsize=2.5truein\efig flux5 y2.2
}
\hskip .7truein
\vbox to 2.2truein{\hsize=2.5truein\efig flux6 y2.2
}
}                                                     
\captionwidth{6truein}
\botcaption{Figure 9}
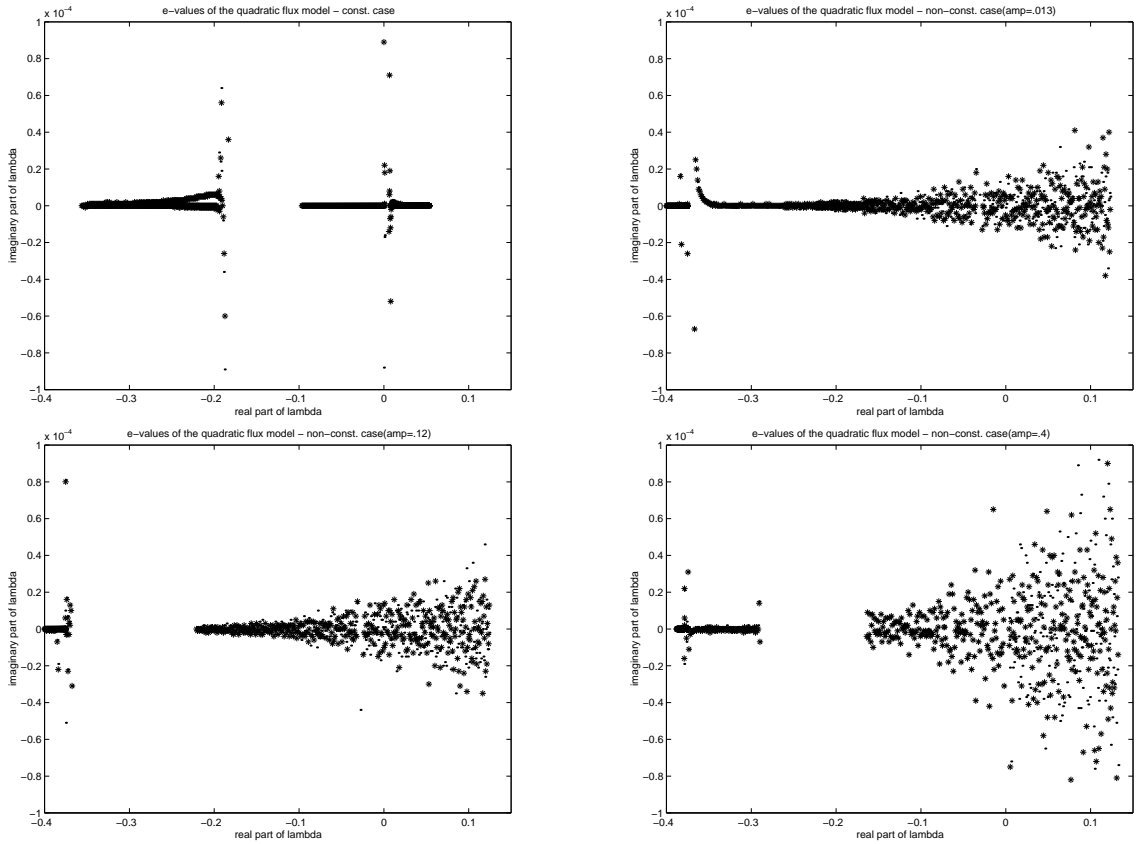
The spectral curves of Evans function $D$.
\endcaption
\endinsert

\midinsert
\hbox{
\vbox to 2.1truein{\hsize=2.5truein\efig homo85 y2.2
}
\hskip .8truein 
\vbox to 2.1truein{\hsize=2.5truein\efig homo86 y2.2
}
}        
\captionwidth{5truein}
\botcaption{Figure 10}
The break-up of the homoclinic orbit.
\endcaption
\endinsert

\midinsert
\captionwidth{2.5truein}
\hbox{                
\vbox to 3.2truein{\hsize=2.5truein \efig b y2.2
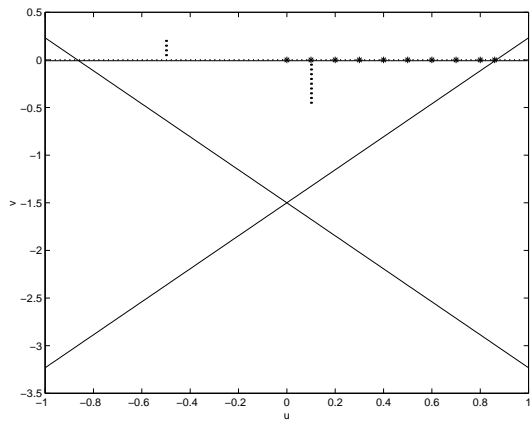
\botcaption{Figure 11a} 
The plotted points represent the vertices of one- and two-cycle for which the
experiments were carried out.
\endcaption
}
\hskip .8truein
\vbox to 3.2truein{\hsize=2.5truein \efig b1 y2.2
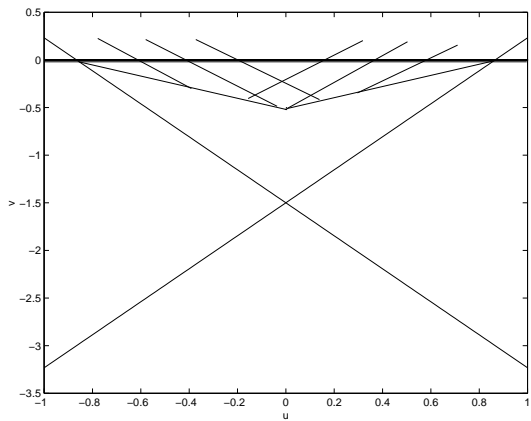
\botcaption{Figure 11b} 
The region where vertices of homoclinics lie.
\endcaption
}                         
}
\endinsert

The quadratic-flux model admits, for $s=0$, a single $3$-cycle
solution, a $1$-parameter family of $2$-cycle solutions, and 
a $2$-parameter family of homoclinic solutions, each containing
a $1$-parameter family of nested periodic orbits.
We restrict to the case s near zero, for which these are the only
periodic solutions. The location of $3$- and $2$-cycles may be 
determined analytically; however, a new issue as compared to the 
van der Waals case is that the
region of existence for homoclinics is no longer analytically
obtainable. {\it Note:} The regions of existence
for two-cycles and three-cycles are analytically obtainable,
as discussed in [AMPZ.3] and [GZ].
We determine the region of existence of homoclinic orbits numerically
by varying $(u_-, v_-)$ where $u_- =(-0.86, 0.86),
v_- \in [-0.5, 0.5)$. Figure 10 shows the break-up of the homoclinic orbit of
this model.                        
Figure 11A shows the tested points $(u_-, v_-)$
which are vertices of homoclinics ($\cdot$)
and two-cycles (*), respectively. The vertices of the (unique) three-cycle
are the three vertices of the surrounding triangle, the vertices
of two-cycles lie on the edges, and the vertices of homoclinics
lie within shaded region (Fig. 11B). For homoclinics, we only need 
to treat the region of $(u_-, v_-)$ inside the triangle since every
outside point has a corresponding inside point yielding the same
dynamical system (i.e., stationary points occur in pairs, one inside,
one outside).
 
Inside each polycycle, we find
several periodic orbits. To these, we apply 
the error analysis and the winding
number calculations as for the previous model.
The same error analysis applies here as in the previous section so
we shall not discuss this aspect.
The maximum Rouch\'{e}'s error over the entire region of
testing was $C 10^{-3}$ for some positive constant $C$, 
where in each case the winding number is 
greater than $1$.  
Figures 12, 13, and 14 show
the instability of each representative periodic orbit within
three-cycle, two-cycle and homoclinic. 
We see that small amplitude orbits and also large amplitude
orbits lying near the three-cycle
are unstable, as shown in Theorem \thmref{quadstab}, 
and also intermediate amplitude
orbits are unstable.
{\it Thus, for this model system, too,
we have instability across the whole region of existence.}

Finally, similarly as in the case of the van der Waals equations,
we may track the number of unstable $0$-eigenvalues of the
periodic waves with respect to amplitude in the homoclinic,
two-, and three-cycle cases.
The result in the homoclinic case is that there is a transition from
zero to one unstable root as amplitude increases.  In the two-
and three-cycle cases, there is no such transition, indicating that
waves of all amplitudes are stable under $T$-periodic perturbations,
even though they
are unstable under  perturbations of general period.
(Recall, $T$-periodic stability can be shown analytically for
small-amplitude waves
by consideration of the constant-coefficient case;
see Remark \thmref{smallamp}.)

\midinsert
\bigskip
\hskip 1.5truein
\vbox to 3.5truein{\hsize=2.5truein \efig wind12-3 y3.3
\botcaption{Figure 12} 
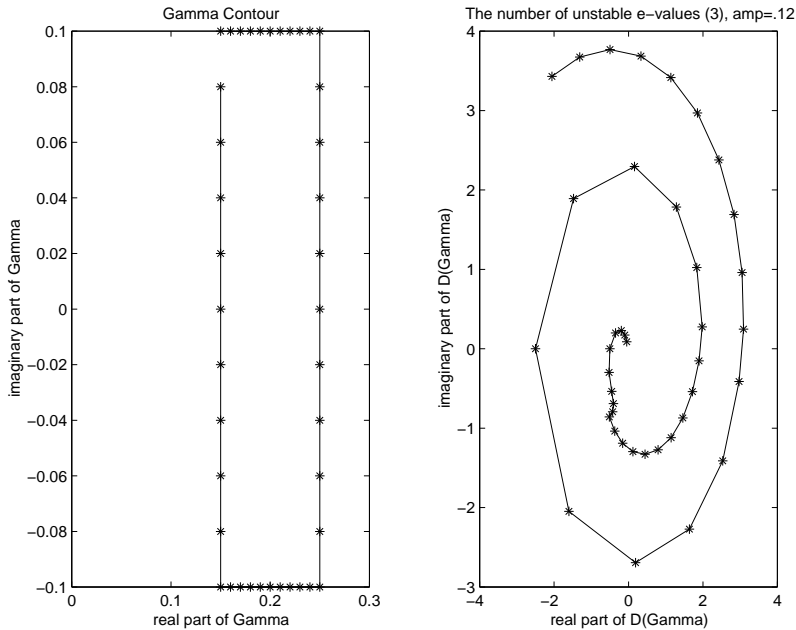
The number of unstable eigenvalues, three-cycle case.
\endcaption
}
\bigskip
\endinsert

\midinsert
\bigskip
\hskip 1.5truein 
\vbox to 3.5truein{\hsize=2.5truein \efig windtwo-3 y3.3
\botcaption{Figure 13} 
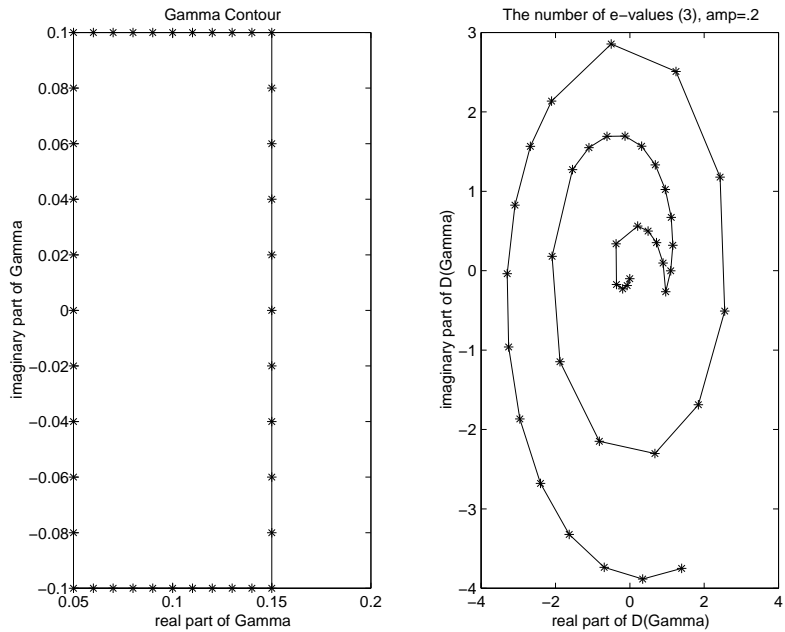
The number of unstable eigenvalues, two-cycle case.
\endcaption
}
\bigskip
\endinsert

\midinsert
\hskip 1.5truein  
\vbox to 3.5truein{\hsize=2.5truein \efig windhomo-3 y3.3
\botcaption{Figure 14} 
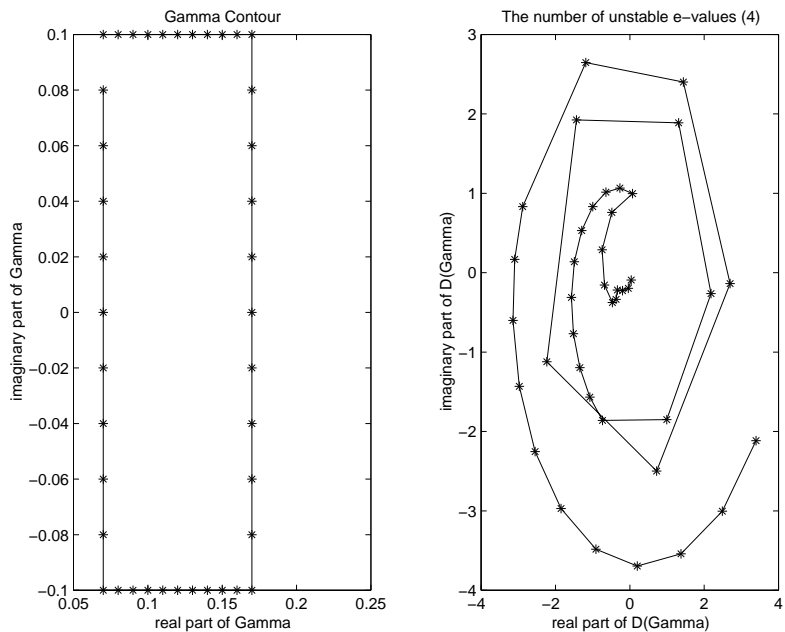
The number of unstable eigenvalues, homoclinic case.
\endcaption
}
\bigskip
\endinsert

\medskip
{\bf Discussion}

The instability of stationary periodic solutions of
the example systems of Section 7
has been shown numerically.
The code is general and easy to manipulate.
It would be quite interesting for more general systems 
to either establish analytically
a corresponding result of uniform instability
or find a counter-example of a stable waves.
Note that the code can be used also to investigate stability 
(for all $k$, not one $k$),
though it has some difficulties near $k=0$
due to numerical sensitivity associated with the multiple root
in $\lambda$ occurring there.
This allows the systematic exploration of periodic wave stability
for families of conservation laws, and also
for Cahn Hilliard and reaction diffusion equations.

\bigskip
\sectionnumber=10
\theoremnumber=0
\equationnumber=0
\smallskip
\TagsOnLeft
\specialhead
Appendix B.\quad Stationary phase approximation
\endspecialhead     

In this appendix, we reproduce 
some unpublished computations of [AMPZ.5] referred to in the introduction,
concerning the linear response under perturbation of a linearly unstable
constant solution $u\equiv u_0$ of a conservation law with viscosity
\eqnref{1}.  Consider the linearized equation \eqnref{1.3},
where $A$ and $B$ are now constant.
Taking the Fourier transform, we obtain
$$
\hat v_t= P(ik)\hat v:= (-ikA-k^2B)\hat v.
\eqnlbl{ft}
$$
The fundamental solution of \eqnref{ft} is clearly $\hat\Psi(k,t)=e^{P(ik)t}$,
whence the fundamental solution of \eqnref{1.3}, i.e., the solution
with initial data given by a Dirac delta function centered at the origin,
is given by the inverse Fourier transform
$$
\Psi(x,t)= \frac{1}{2\pi}
\int_{-\infty}^{+\infty} e^{ikx}e^{P(ik)t} \, dk,
\eqnlbl{ift}
$$
and the Green function for \eqnref{1.3} by $G(x,t;y)=\Psi(x-y,t)$.

The $n$ eigenvalues of the symbol $P(ik)$ determine $n$ dispersion
relations
$$
\lambda= \lambda_j(k),
\quad j=1,\dots, n,
\eqnlbl{disp}
$$
with linearized stability corresponding to the condition $\R \lambda_j \le 0$ 
for all $k\in \BbbR$.
Denote by $k_j^*$ the critical frequency at which $\R\lambda_j$
takes on its maximum value.
For simplicity, we restrict to the generic situation that $\lambda_j$ 
is isolated from the rest of the spectrum of $P(ik)$ at $k=k_j^*$, 
and $\lambda_j(k_j^*)$ is a nondegenerate maximum.
Thus, $\lambda_j$ admits a second-order Taylor expansion
$$
\lambda_j(k)= \lambda_j(k_j^*) - i\alpha_j (k-k_j^*) - \beta_j (k-k_j^*)^2
+ \dots
\eqnlbl{taylor2}
$$
about $k_j^*$, where $\alpha_j$ is real, and $\R \beta_j >0$.
Let $R_j(k)=R_j(k_*) + \dots$ and 
$L_j(k)=L_j(k_*) + \dots$ denote right and left
eigenvectors and $\Pi_j:= R_j L_j^*/\langle R_j,L_j\rangle$ 
the spectral projection of $P$ associated with $\lambda_j$.

Decomposing \eqnref{ift} by spectral resolution as
$\Psi(x,t)=  \sum_j \Psi_j(x,t)$, where
$$
\Psi_j(x,t) :=\sum_j \frac{1}{2\pi}
\int_{-\infty}^{+\infty} e^{ikx}e^{\lambda_j(k)t}\Pi_j(k) \, dk,
\eqnlbl{sift}
$$
we may estimate each of the $\Psi_j$ by stationary phase approximation
about $k_j^*$ as
$$
\aligned
\Psi_j(x,t) &\sim \bar \Psi_j:=
\frac{1}{2\pi}
\int_{-\infty}^{+\infty} e^{ikx}e^{(
\lambda_j(k_j^*) - i\alpha_j (k-k_j^*) -  \beta_j (k-k_j^*)^2
)t}\Pi_j(k_j^*) \, dk\\
&=\frac{1}{2\pi}
e^{ik_j^*x + \lambda(k_j^*)t} \Pi_j(k_j^*) 
\int_{-\infty}^{+\infty} e^{i(k-k_j^*)x}e^{(
- i\alpha_j (k-k_j^*) -  \beta_j (k-k_j^*)^2
)t}
\, dk.\\
\endaligned
\eqnlbl{stat}
$$

Let us now restrict to the simplest case $B=I$.
In this situation, the dispersion relations are exactly
$\lambda_j(k)= -ik a_j - k^2$, where $a_j$ are the eigenvalues of $A$,
and the associated eigenprojections $\Pi_j$ are just the 
eigenprojections $r_j l_j^{*}/\langle r_j,l_j\rangle$ of $A$ associated
with $a_j$, where $r_j$ and $l_j$ denote associated right and left
eigenvectors of $A$.
Thus, linear stability is equivalent to reality of $a_j$, or
hyperbolicity of the first-order coefficient $A$.
By direct calculation, we find that $k_j^*= \Im a_j/2$,
$\alpha_j= \R a_j$, 
$\lambda(k_j^*)=(\Im a_j/2)^2 - i \R a_j \Im a_j/2$,
and so
$$
\alpha_j = - \Im \lambda(k_j^*)/k_j^*.
\eqnlbl{cancel}
$$
Likewise, noting that $\beta_j$ is real, we may explicitly carry out
the Fourier inversion in the final factor of \eqnref{stat} to 
obtain
$$
\aligned
\bar \Psi_j&=
e^{ik_j^*x + \lambda(k_j^*)t} 
\Pi_j(k_j^*) \times
\frac{
e^{-(x-\alpha_j t)^2/4\beta_j t}  }
{\sqrt{4\pi \beta_j t}} \\
&=e^{\R \lambda(k_j^*)t} 
\frac{r_j l_j^{*}}{\langle r_j,l_j\rangle}
\times
\frac{e^{ik_j^*(x-\alpha_jt)}
e^{-(x-\alpha_j t)^2/4t}  }
{\sqrt{4\pi t}}\\
&=e^{(\Im a_j/2)^2 t} 
\frac{r_j l_j^{*}}{\langle r_j,l_j\rangle}
\times
\frac{e^{i(\Im a_j/2)(x-\R a_j t)}
e^{-(x-\R a_j t)^2/4t}  }
{\sqrt{4\pi t}},\\
\endaligned
\eqnlbl{finalform}
$$
where we have used \eqnref{cancel} in the second equality.
That is, the behavior of $\Phi$ in the $j$th mode is, to lowest order,
that of a time-exponentially growing modulated Gaussian wave packet
traveling with speed $\alpha_j$,
with exponential growth rate equal to the maximum real part
$\R \lambda_j(k_j^*)$ of $\lambda_j$, and modulating spatial
frequency equal to the critical frequency $k_j^*$.

In the case of a linearly stable mode $\Im a_j=0$, \eqnref{finalform}
reduces to the usual Gaussian ``linear diffusion wave''
$\frac{r_j l_j^{*}}{\langle r_j,l_j\rangle}
\times
\frac{ e^{-(x-\R a_j t)^2/4t}  }
{\sqrt{4\pi t}}$
of [LZe].
By analogy, we might call \eqnref{finalform} an ``unstable 
linear diffusion wave.''
Note that unstable modes $\Im a_j\ne 0$ group in conjugate pairs, with
equal propagation speeds $\alpha_j=\R a_j$,
and thus their diffusion waves combine to give a single real, sinusoidal 
diffusion wave of dimension two.
The proper analogy is, therefore, to the ``generalized diffusion
wave'' defined in [LZe] in the case of repeated real eigenvalues of $A$
(the boundary case between strict hyperbolicity and ellipticity
of $A$).

{\bf Example \thmlbl{vangas}.}
For the equations of van der Waal gas dynamics or elasticity,
we have $B=I$ and 
$
A=df(u_0,v_0)= \pmatrix 0 & -1 \\ p'(v_0) & 0 \endpmatrix$,
has eigenvalues $a_j= \pm \sqrt { -p'(v_0)}$. 
Thus, linearized stability of the constant solution
$(u,v)\equiv (u_0,v_0)$ is equivalent to $p'(v_0)<0$,
and failure of stability entails a pair of pure imaginary
eigenvalues $a_j=\pm i\tau$.
In this case, therefore, there is a single conjugate pair
of {\it standing} unstable linear diffusion waves, $\alpha_j\equiv 0$,
with time-exponential growth rate $e^{\tau^2 t/4}$ and modulating
spatial frequency $e^{i\tau x/2}$.
\medskip
{\bf Application to pattern formation.}
Now, let us consider the implications of the linear estimate
\eqnref{finalform} on the behavior under perturbation of solutions of 
the {nonlinear}
equations, for initial data approaching a linearly unstable 
constant state $u_{\infty}$ as $x\to +\infty$.
In view of the global nonlinear structure of our models,
let us make the crude approximation that nonlinear effects
take over at some prescribed amplitude $A$,
turning off the time-exponential growth and
converting the sinusoidal oscillations of the linearized
approximation \eqnref{finalform} to a series of traveling fronts
{\it moving with approximately zero speed} (the last assumption
implies some normalization to achieve a rest frame).
Looking at a single unstable mode $j$,
the picture that emerges from this crude model on the half line
$[0,+\infty)$ is an oscillatory pattern of phase transitions
separated from the constant state by a modulated front
whose location $x(t)$ is determined by the property that 
$$
|\bar\Psi(x_j(t), t)|\sim A,
\eqnlbl{frontcrit}
$$
i.e., ignoring the algebraic factor $t^{-1/2}$ in \eqnref{finalform},
travels with approximate speed $\sigma_j$ determined by
$ (\Im a_j/2)^2t -(\sigma_jt  - \alpha_jt)^2/4t =0 $.
Solving, and recalling that $\alpha_j=\R a_j$,
$\R \lambda(k_j^*)=(\Im a_j/2)^2$, we obtain a value of
$$
\sigma_j= \R a_j + |\Im a_j|
\eqnlbl{sigma}
$$
for the approximate front speed.

Since waves under this model are ``born'' in the linearized regime,
we can estimate the approximate wave-length inherited in the nonlinear
regime by dividing the distance $d=\sigma_j t$ swept out by the front
$x(t)$ from time $0$ to time $t$,
divided by the number of peaks of the linearized unstable
diffusion wave \eqnref{finalform} that it has
crossed during that time.  Working in the rest frame $x=\alpha_j t$ 
of the diffusion wave, we find that the number of peaks is
$(\sigma_j-\alpha_j)t$ divided by the wavelength $2\pi/k_j^*=2\pi/\Im a_j/2$
of its spatial modulation, yielding an estimated wavelength of
$$
T_j=
\frac{\sigma_jt}
{\frac{(\sigma_j-\alpha_j)t}{(2\pi/k_j^*)}}\,.
$$
This predicts an average (Doppler-shifted)
frequency for the nonlinear pattern of
$$
\aligned
\tilde k_j=2\pi/T_j&= k_j^* \Big( \frac{\sigma_j-\alpha_j}{ \sigma_j}\Big)\\
&=k_j^* \Big( \frac{|\Im a_j|}{ \R a_j + |\Im a_j|}\Big)\\
\endaligned
\eqnlbl{freq}
$$
In the (standing) case $\alpha_j=\R a_j=0$ of Example \thmref{vangas}, this is
just the critical frequency $k_j^*=\Im a_j/2$.
Note: here and above we are implicitly assuming that the front speed $\sigma_j$
is positive, i.e., the pattern is expansive.  This is in fact a necessary
condition for pattern formation to occur.

Numerical experiments in [AMPZ.4] showed close agreement between the
predictions of \eqnref{sigma} and \eqnref{freq} and behavior of
actual patterns in solutions of \eqnref{1}.
Likewise, the condition $\sigma_j\ge 0$ was seen to well predict
{\it nonlinear} instability, and the appearance of pattern formation,
as model parameters (e.g., Riemann endstates) were varied.

{\bf Remark \thmlbl{finalrmk}.}
In the linearly stable case, Liu and Zeng [LZe] derive also
a self-similar {\it nonlinear} diffusion wave refining the
linear estimate \eqnref{finalform},  by 
judiciously appending to an approximately diagonalized
version of the linearized equations \eqnref{1.3}
the ``diagonal'' quadratic-order terms in the Taylor expansion
of the nonlinear flux $f$.
This approximation makes essential use of the small-amplitude 
nature of solutions. 
In the large-amplitude linearly unstable case, we have
no such simple description of a nonlinear diffusion wave;
indeed, the relevant entity seems to be the entire metastable
pattern described above.
We might neatly summarize the conclusions of this paper (together with those
of [AMPZ.4--5]) as suggesting that 
patterns observed in dynamical phase-transitional
models represent {\it metastable nonlinear diffusion waves}
rather than {\it stable periodic states}.
\bigskip
\noindent{\bf Acknowledgment.}
K. Zumbrun thanks Rob Gardner for several illuminating conversations
regarding his foundational work on periodic stability.
Thanks also to Richard Laugesen for pointing out the result
of [Sc,LP] used in Section 7 to establish monotonicity of the period 
of a double-well oscillator, and to the referee for pointing out
the references [BM,M].
Finally, thanks to Eva Marie Elliot of Indiana University
for her skillful and efficient typing of this manuscript.
Research of both authors was supported in part by the 
National Science Foundation under Grants No. DMS-9706842
and DMS-007065.

\Refs
\medskip\noindent
[AGJ] J. Alexander, R. Gardner and C.K.R.T. Jones,
{\it A topological invariant arising in the analysis of
traveling waves}, J. Reine Angew. Math. 410 (1990) 167--212.
\medskip\noindent
[AMPZ.1] A. Azevedo, D. Marchesin, B. Plohr, and K. Zumbrun,
{\it Nonuniqueness of solutions of Riemann problems,}
 Z. Angew. Math. Phys. 47 (1996) 977--998. 
\medskip\noindent
[AMPZ.2] A. Azevedo, D. Marchesin, B. Plohr and K. Zumbrun,
{\it Bifurcation from the constant state of nonclassical
viscous shock waves,} Comm. Math. Phys. 202 (1999) 267--290. 
\medskip\noindent
[AMPZ.3] A. Azevedo, D. Marchesin, B. Plohr and K. Zumbrun,
{\it Capillary instability in multiphase flow models,} preprint (1999).
\medskip\noindent
[AMPZ.4] A. Azevedo, D. Marchesin, B. Plohr and K. Zumbrun,
{\it Long-lasting diffusive solutions for systems of conservation laws,} VI
Workshop on Partial Differential Equations, Part I (Rio de Janeiro, 1999). 
Mat. Contemp. 18 (2000), 1--29.
\medskip\noindent
[AMPZ.5] A. Azevedo, D. Marchesin, B. Plohr and K. Zumbrun,
{\it private communication}, unpublished numerical computations (1996).
\medskip\noindent
[BSZ] S. Benzoni-Gavage, D. Serre, and K. Zumbrun,
{\it Alternate Evans functions and viscous shock waves,}  
to appear, SIAM J. Math. Anal. 32 (2001) 929--962.
\medskip\noindent
[BMi.1] T. Bridges and A. Mielke, 
{\it A proof of the Benjamin-Feir instability,}
Arch. Rat. Mech. Anal. 133 (1995) 145--198.
\medskip\noindent
[BMi.2] T. Bridges and A. Mielke, 
{\it Instability of spatially-periodic states
for a family of semilinear PDE's on an infinite strip,}
Math. Nachr. 179 (1996) 5--25.
\medskip \noindent
[Br.1] L. Q. Brin, {\it Numerical testing of the stability of viscous
shock waves}, Ph.D. dissertation, Indiana University, May 1998.

\medskip\noindent
[Br.2] L. Q. Brin,
{\it Numerical testing of the stability of viscous shock waves,}
Math. Comp. 70 (2001) 1071--1088. 
\medskip \noindent
[BZ] L. Q. Brin and K. Zumbrun, {\it Analytically varying eigenvectors
and the stability of viscous shock waves,} to appear, Mat. Contemp.
\medskip \noindent
[BF] R. L. Burden, and J. D. Faires, {\it Numerical Analysis},
PWS Publishing Company, Boston, Massachusetts, fifth edition, 1993.
\medskip\noindent
[C] S. \v{C}ani\'c, {\it Nonexistence of Riemann solutions for a quadratic model
deriving from petroleum engineering,} preprint (1998).
\medskip\noindent
[CGS.1] J. Carr, M. Gurtin, and M. Slemrod,
{\it Structured phase transitions on a finite interval,}
Arch. Rational Mech. Anal. 86 (1984) 317--351.
\medskip\noindent
[CGS.2] J. Carr, M. Gurtin, and M. Slemrod,
{\it One-dimensional structured phase transformations under prescribed loads,}
 J. Elasticity 15 (1985) 133--142. 
\medskip\noindent
[CL] L. Cairo and J. Llibre,
{\it Integrability and algebraic solutions for planar polynomial
differential systems with emphasis on the quadratic systems,}
IME/USP (Universidade de S\~ao Paulo, Brazil) Report No. 2, (1999).
\medskip\noindent
[CP] S. \v{C}ani\'c and G. Peters, 
{\it On the oscillatory solutions in hyperbolic conservation laws,}
Nonlinear Anal. Real World Appl. 1 (2000) 287--314.
\medskip\noindent
[DGK] A. Doelman, R.A. Gardner, and T. Kaper,
{\it A stability index analysis of 1-D patterns of the
Gray-Scott Model,} to appear, AMS Memoirs.
\medskip\noindent
[E.1] J.W. Evans, 
{\it Nerve axon equations: I. Linear approximations,}
Ind. Univ. Math. J. 21 (1972) 877--885.
\medskip\noindent
[E.2] J.W. Evans, 
{\it Nerve axon equations: II. Stability at rest,}
Ind. Univ. Math. J. 22 (1972) 75--90.
\medskip\noindent
[E.3] J.W. Evans, 
{\it Nerve axon equations: III. Stability of the nerve impulse,}
Ind. Univ. Math. J. 22 (1972) 577--593.
\medskip\noindent
[E.4] J.W. Evans, 
{\it Nerve axon equations: IV. The stable and the unstable impulse,}
Ind. Univ. Math. J. 24 (1975) 1169--1190. 
\medskip\noindent
[Es] E. Eszter,  
{\it Evans function analysis of periodic travelling
wave solutions of the FitzHugh-Nagumo system}, 
Doctoral thesis, University of Massachussetts, Amherst (1999).
\medskip\noindent
[FreZ] H. Freist\"uhler and K. Zumbrun,
{\it Examples of unstable viscous shock waves,}
unpublished note, Institut f\"ur Mathematik, RWTH Aachen, February 1998.
\medskip\noindent
[FL.1] H. Frid and I.-S. Liu,
 {\it Oscillation waves in Riemann problems for
phase transitions}, Quart. Appl. Math. 56 (1998) 115--135.
\medskip\noindent
[FL.2] H. Frid and I.-S. Liu,
{\it Oscillation waves in Riemann problems inside
elliptic regions for conservation laws of mixed type}, 
Z. Angew. Math. Phys. 46 (1995) 913--931. 
\medskip\noindent
[G.1] R.A. Gardner, 
{\it On the structure of the spectra of periodic travelling waves}, 
J. Math. Pures Appl. (9) 72 (1993) 415--439.
\medskip\noindent
[G.2] R.A. Gardner, 
{\it Spectral analysis of long wavelength periodic waves and applications}, 
J. Reine Angew. Math. 491 (1997), 149--181.
\medskip\noindent
[G.3]  R.A. Gardner,
{\it Instability of oscillatory shock profile solutions of the generalized 
Burgers-KdV equation},
Phys. D 90 (1996) 366--386.  
\medskip\noindent
[GZ] R. Gardner and K. Zumbrun, 
{\it The Gap Lemma and geometric criteria for instability
of viscous shock profiles}, 
Comm. Pure Appl. Math. 51 (1998) 797--855.
\medskip\noindent
[GH] J. Guckenheimer and P. Holmes,
{\it Nonlinear oscillations, dynamical systems, and bifurcations of
vector fields},
(Revised and corrected reprint of the 1983 original),
Springer--Verlag, New York (1990), xvi+459 pp.
\medskip\noindent
[GM] M.E. Gurtin and H. Matano, 
{\it On the structure of equilibrium phase
transitions within the gradient theory of fluids,}
Quart. Appl. Math. 46 (1988) 301--317.
\medskip\noindent
[HK] J. Hale and H. Kocak,
{\it Dynamics and bifurcations,}
Texts in Applied Mathematics, 3. Springer-Verlag, New York, 1991.
xiv+568 pp. ISBN: 0-387-97141-6. 
\medskip\noindent
[He] D. Henry, 
{\it Geometric theory of semilinear parabolic equations},
Lecture Notes in Mathematics, Springer--Verlag, Berlin (1981),
iv + 348 pp.
\medskip\noindent
[HP] J. Hurley and B. Plohr,
{\it Some effects of viscous terms on solutions of Riemann
problems,} Mat. Contemp. 8 (2995) 203--204.
\medskip\noindent
[IR] E. Infeld and G. Rowlands,
{\it Nonlinear waves, solitons, and chaos,}
Second edition, Cambridge University Press, Cambridge (2000)
xiv+391 pp. ISBN: 0-521-63212-9; 0-521-63557-8. 
\medskip\noindent
[Ja] R.D. James, 
{\it The propagation of phase boundaries in elastic bars,}
 Arch.  Rational Mech. Anal. 73 (1980) 125--158. 
\medskip\noindent
[J] C.K.R.T. Jones, 
{\it Stability of the travelling wave solution of the FitzHugh--Nagumo system},
 Trans. Amer. Math. Soc.  286 (1984) 431--469.
\medskip\noindent
[KS] T. Kapitula and B. Sandstede,
{\it Stability of bright solitary-wave solutions 
to perturbed nonlinear Schrödinger equations,} Phys. D 124
(1998) 58--103.
\medskip\noindent
[Ka] S. Kawashima,
{\it Systems of a hyperbolic--parabolic composite type,
 with applications to the equations of magnetohydrodynamics},
thesis, Kyoto University (1983).
\medskip\noindent
[LP] R. Laugesen and M. Pugh,
{\it Properties of steady states for thin film equations,}
Euro. J. Appl. Math. 11 (2000) 293--351.
\medskip\noindent
[LF] I.-S. Liu and H.Frid,
{\it Phase mixtures in dynamics of pseudoelasticity},
Contin. Mech. Thermodyn. 7 (1995) 475--487. 
\medskip\noindent
[LZe] T.-P. Liu and Y. Zeng, 
{\it Large time behavior of solutions for general 
quasilinear hyperbolic--parabolic systems of conservation laws},
AMS memoirs 599 (1997).
\medskip\noindent
[LZ.1] T.P. Liu and K. Zumbrun,
{\it Nonlinear stability of an undercompressive shock for complex
Burgers equation,} Comm. Math. Phys. 168 (1995) 163--186.
\medskip\noindent
[LZ.2] T.P. Liu and K. Zumbrun,
{\it On nonlinear stability of general undercompressive viscous shock waves,}
Comm.  Math. Phys. 174 (1995) 319--345.
\medskip \noindent
[M.1] K. Maginu, {\it Existence and stability of periodic traveling
wave
solutions to Nagumo's nerve equation}, J. MATH. Biol. 10 (1980),
133-153.

\medskip \noindent
[M.2] K. Maginu, {\it Stability of periodic traveling wave solutions
with large spatial periods in reaction-diffusion systems}, J. Diff.
Equation,
39 (1981), 73-99.

\medskip\noindent
[MP] A. Majda and R. Pego,
{\it Stable viscosity matrices for systems of conservation laws},
J. Diff. Eqs. 56 (1985) 229--262.   
\medskip\noindent
[Mi] A. Mielke, {\it Instability and stability of rolls in the
Swift-Hohenberg equation,}
Comm. Math. Phys. 189 (1997) 829--853.
\medskip\noindent
[OZ] M. Oh and K. Zumbrun,
{\it Stability of periodic solutions of conservation laws with viscosity:
Pointwise bounds on the Green function,}
Arch. Rational Mech. Anal.
\medskip\noindent
[PW] R. L. Pego and M.I. Weinstein, 
{\it  Eigenvalues, and instabilities of solitary waves},
Philos. Trans. Roy. Soc. London Ser. A 340 (1992) 47--94.
\medskip \noindent
[SS.1] B. Sandstede and A. Scheel, 
{\it On the stability of periodic travelling waves with large spatial period,}
J. Differential Equations 172 (2001) 134--188. 
\medskip \noindent
[SS.2] B. Sandstede and A. Scheel, 
{\it On the structure of spectra of modulated traveling waves,}
Math. Nachr. 232 (2001) 39--93.
\medskip\noindent
[Sc] R. Schaaf,
{\it A class of Hamiltonian systems with increasing periods,}
J. Reine Angew. Math. 363 (1985) 96--109.
\medskip\noindent
[SSh] S. Schecter and M. Shearer,
{\it Undercompressive shocks for non-strictly
hyperbolic conservation laws,}
J. Dynamics Differential Equations 3 (1991), no. 2, 199--271.
\medskip\noindent
[S.1] G. Schneider, {\it Nonlinear diffusive stability
of spatially periodic solutions-Abstract theorem and higher space
dimensions}, preprint (1998).

\medskip \noindent
[S.2] G. Schneider, 
{\it Diffusive stability of spatial periodic solutions of the 
Swift-Hohenberg equation,} (English. English summary) 
Comm. Math. Phys. 178 (1996) 679--702. 

\medskip \noindent
[S.3] G. Schneider, 
{\it Nonlinear stability of Taylor vortices in infinite cylinders,}
Arch. Rational Mech. Anal. 144 (1998) 121--200.  

\medskip\noindent
[Se.1] D. Serre, 
{\it Entropie du m\'elange liquide-vapeur d'un fluide 
thermo-capillaire}, 
Arch.  Rational Mech. Anal. 128 (1994) 33--73.
\medskip\noindent
[Se.2] D. Serre, 
{\it Transitions de phase et oscillations de grande amplitude},
Qualitative aspects and applications
of nonlinear evolution equations (Trieste, 1990), 172--184, World Sci. Publishing, River Edge,
NJ, 1991. 
\medskip \noindent
[Sh.1] M. Shearer,
{\it The Riemann problem for a class of conservation laws of mixed type},
J. Diff. Eqs. 46 (1982) 426--443.
\medskip\noindent
[Sh.2] M. Shearer,
{\it Nonuniqueness of admissible solutions of Riemann inintial value
problems for a system of conservation laws of mixed type},
Arch. Rat. Mech. Anal. 93 (1986) 45--59.
\medskip\noindent
[Sh.3] M. Shearer,
{\it Dynamic phase transitions in a van der Waals gas},
Quart. Appl. Math. 46 (1988) 631--636.
\medskip\noindent
[SSMP] M. Shearer, D. Schaeffer, D. Marchesin, and P.L. Paes-Leme, 
{\it Solution of the Riemann problem for a prototype $2\times 2$ system of
nonstrictly hyperbolic conservation laws,}
Arch. Rational Mech. Anal. 97 (1987) 299--320. 
\medskip\noindent
[Sl.1] M. Slemrod, 
{\it The viscosity-capillarity approach to phase transitions,}
in: 
{\it PDEs and continuum models of phase transitions} (Nice, 1988), 201--206,
Lecture Notes in Phys., 344,
Springer, Berlin-New York, 1989.
\medskip\noindent
[Sl.2] M. Slemrod,
{\it The vanishing viscosity-capillarity approach 
to the Riemann problem for a van der Waals fluid,}
in: {\it Nonclassical continuum mechanics} (Durham, 1986), 325--335,
London Math. Soc. Lecture Note Ser., 122,
Cambridge Univ. Press, Cambridge-New York, 1987.
\medskip\noindent
[Sl.3] M. Slemrod, 
{\it A limiting ``viscosity'' approach to the 
Riemann problem for materials exhibiting change of phase,}
Arch. Rational Mech. Anal. 105 (1989) 327--365.
\medskip\noindent
[Sl.4] 
M. Slemrod, 
{\it Dynamic phase transitions in a van der Waals fluid,}
J. Diff. Eqs. 52 (1984) 1--23.
\medskip\noindent
[Sl.5] M. Slemrod,
{\it Admissibility criteria for propagating phase 
boundaries in a van der Waals fluid,}
Arch. Rational Mech. Anal. 81 (1983) 4, 301--315.
\medskip\noindent
[ST] M. Slemrod and A.E.Tzavaras, 
{\it Shock profiles and self-similar fluid dynamic limits,} 
Proceedings of the Second International Workshop on Nonlinear Kinetic Theories
and Mathematical Aspects of Hyperbolic Systems (Sanremo, 1994). 
Transport Theory Stat.  Phys. 25 (1996) 531--541. 

\medskip \noindent
[Sm] J. Smoller, 
{\it Shock waves and reaction--diffusion equations,} 
Second edition, Grundlehren der Mathematischen Wissenschaften
[Fundamental Principles of Mathematical Sciences], 258. Springer-Verlag, 
New York, 1994. xxiv+632 pp. ISBN: 0-387-94259-9. 
\medskip \noindent
[SB] J. Stoer, and R. Bulirsch, {\it Introduction to Numerical
Analysis},
Springer-Verlag, New York, second edition, 1992.
\medskip \noindent
[Z.1] K. Zumbrun, {\it Dynamical stability of phase
transitions in the p-system with viscosity-capillarity},
SIAM J. Appl. Math. 60 (2000), 1913-1929. 

\medskip\noindent
[Z.2] K. Zumbrun,
{\it Computation of scattering coefficients for 
phase transitional shock patterns,} 
unpublished note.

\medskip \noindent
[Z.3] K. Zumbrun,
{\it Multidimensional stability of planar viscous shock waves,}
Advances in the theory of shock waves, 307--516, Progr. Nonlinear
Differential Equations Appl., 47, Birkhäuser Boston, Boston, MA, 2001.
\medskip \noindent
[Z.4] K. Zumbrun, {\it Stability of viscous shock waves}, 
Lecture notes, Indiana University (1998).
\medskip\noindent
[ZH] K. Zumbrun and P. Howard, 
{\it Pointwise semigroup methods and stability of viscous shock waves,}
Indiana Math. J. V47 (1998) 741--871.
\medskip\noindent
[ZPM] K. Zumbrun,  B. Plohr, and D. Marchesin, 
{\it Scattering behavior of transitional waves}, 
Mat. Contemp. 3 (1992) 191--209.
\medskip\noindent
[ZS] K. Zumbrun and D. Serre,
{\it Viscous and inviscid stability of multidimensional 
planar shock fronts,} Indiana Univ. Math. J. 48 (1999) 937--992.
\endRefs
\enddocument
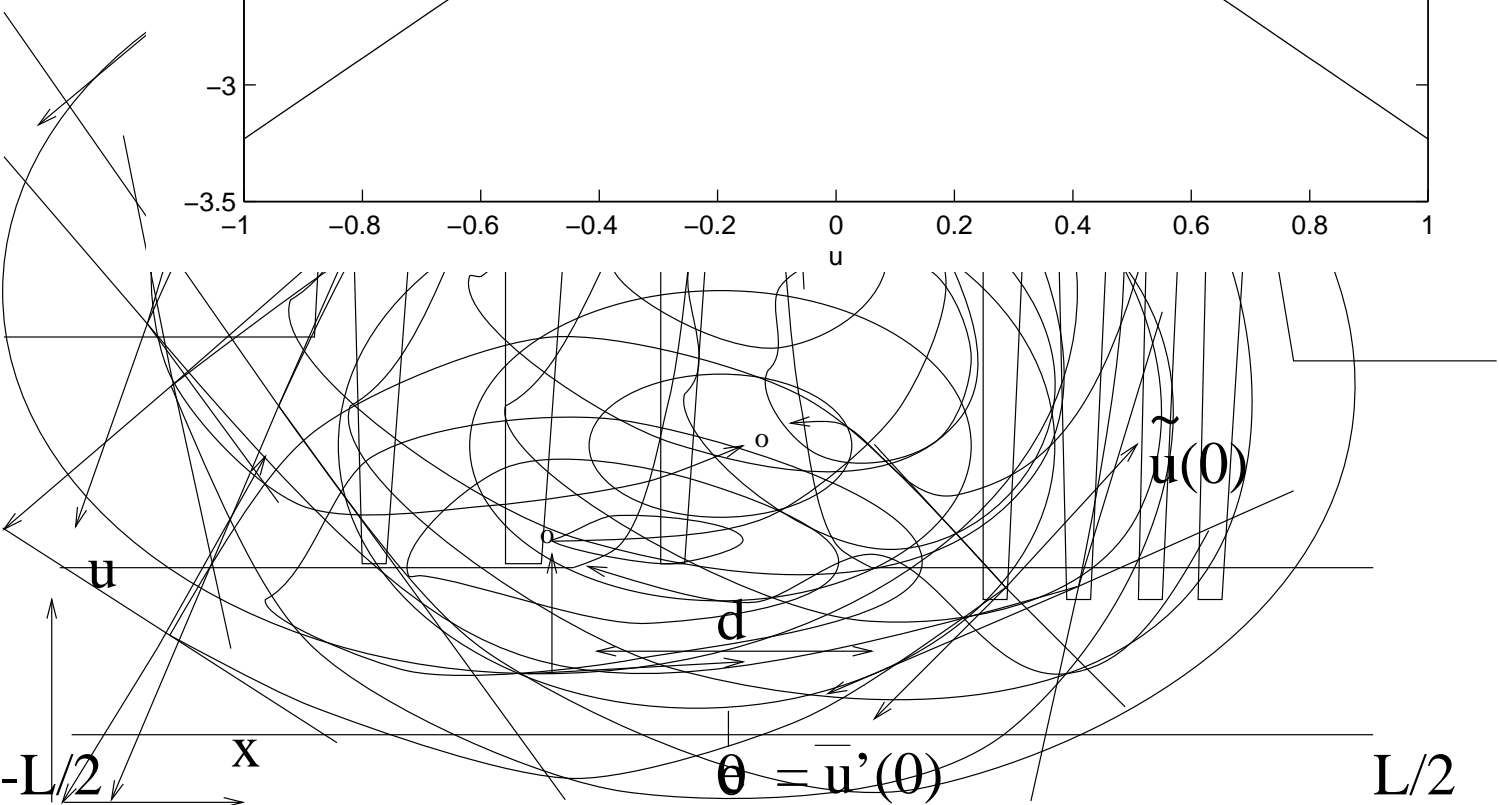
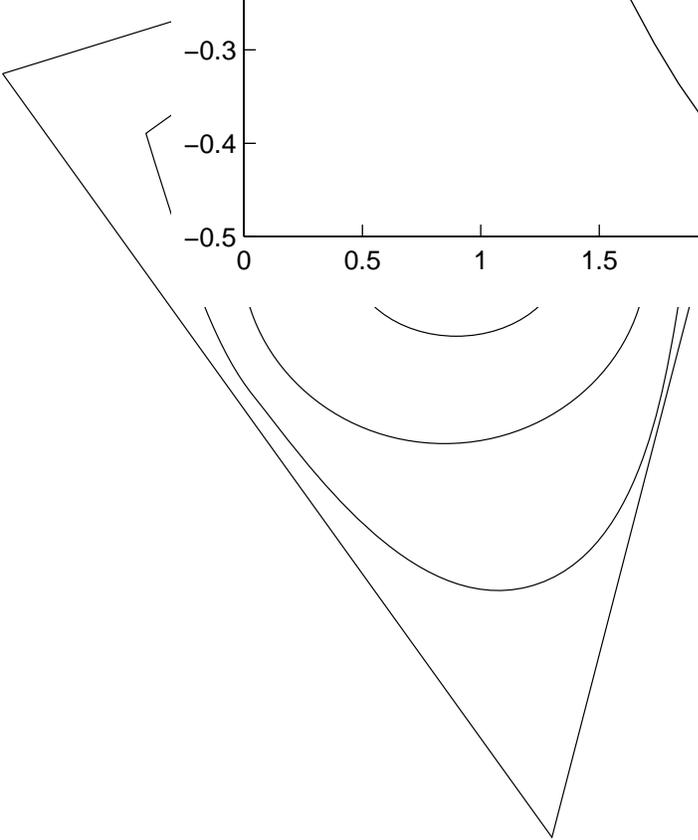
